\documentclass{article}

\usepackage{microtype}
\usepackage{graphicx}
\usepackage{subfigure}
\usepackage{booktabs} 

\usepackage{hyperref}

\usepackage[accepted]{icml2020}

\usepackage{nicefrac}
\usepackage{color}
\usepackage{makecell}
\usepackage{colortbl}
\definecolor{bgcolor}{rgb}{0.66,0.88,1.00}
\usepackage{amssymb}
\usepackage{amsmath}
\usepackage{mathrsfs}
\usepackage{bm}

\newtheorem{theorem}{Theorem}
\newtheorem{lemma}{Lemma}

\newtheorem{definition}{Definition}

\usepackage[colorinlistoftodos,bordercolor=orange,backgroundcolor=orange!20,linecolor=orange,textsize=scriptsize]{todonotes}

\renewcommand{\algorithmicrequire}{ \textbf{Input:}} 
\renewcommand{\algorithmicensure}{ \textbf{Output:}} 

\newcommand{\beq}{\begin{equation}}
\newcommand{\eeq}{\end{equation}}
\newcommand{\R}{{\mathbb R}}

\newcommand{\ns}[1]{\| #1 \|^2}

\newcommand{\nsB}[1]{\left\| #1 \right\|^2}

\newcommand{\n}[1]{\| #1 \|}
\newcommand{\inner}[2]{\langle #1, #2 \rangle}
\newcommand{\innerb}[2]{\left\langle #1, #2 \right\rangle}
\newcommand{\innerB}[2]{\left\langle #1, #2 \right\rangle}
\newcommand{\calC}{{\mathcal C}}
\newcommand{\prox}{\mathrm{prox}}
\newcommand{\head}[1]{\noindent{{\bf #1:}}}
\newcommand{\eat}[1]{}
\newenvironment{proof}{\noindent {\em Proof: }\ignorespaces}%
{\hspace*{\fill}$\Box$\par}
{\hspace*{\fill}$\Box$\par\vspace{4mm}}
\newenvironment{proofof}[1]{\smallskip\noindent{\em Proof of #1.}}%
{\hspace*{\fill}$\Box$\par}

\newcommand{\EB}[1]{\mathbb{E}\left[ #1 \right]}

\newcommand{\eqdef}{:=}

\newcommand{\alglabel}{%
	\addtocounter{ALC@line}{-1}
	\refstepcounter{ALC@line}
	\label
}

\newcommand{\xt}{y^{k}}
\newcommand{\yt}{x^{k}}
\newcommand{\vt}{z^{k}}
\newcommand{\xp}{y^{k+1}}

\newcommand{\vp}{z^{k+1}}

\newcommand{\gt}{g^{k}}
\newcommand{\xs}{x^{*}}

\newcommand{\et}{\eta_{k}}
\newcommand{\at}{\theta_{k}}
\newcommand{\bt}{\beta_{k}}
\newcommand{\st}{\gamma_{k}}

\newcommand{\norm}[1]{\left\| #1 \right\|}
\newcommand{\dotprod}[2]{\langle #1,#2 \rangle}

\newcommand{\E}[1]{\mathbb{E}\left[#1\right]}

\newcommand{\EE}{\mathbb{E}}

\def\<#1,#2>{\langle #1,#2\rangle}

\newcommand{\cC}{{\cal C}}

\newcommand{\cH}{{\cal H}}

\newcommand{\cO}{{\cal O}}

\newcommand{\cW}{{\cal W}}

\newcommand{\cY}{{\cal Y}}
\newcommand{\cZ}{{\cal Z}}

\icmltitlerunning{Acceleration for Compressed Gradient Descent  in Distributed and Federated Optimization}

\begin{document}

\twocolumn[
\icmltitle{Acceleration for Compressed Gradient Descent  in Distributed and \\ Federated Optimization}




\begin{icmlauthorlist}
\icmlauthor{Zhize Li}{kaust}
\icmlauthor{Dmitry Kovalev}{kaust}
\icmlauthor{Xun Qian}{kaust}
\icmlauthor{Peter Richt\'{a}rik}{kaust}
\end{icmlauthorlist}

\icmlaffiliation{kaust}{King Abdullah University of Science and Technology, Thuwal, Kingdom of Saudi Arabia}

\icmlcorrespondingauthor{Zhize Li}{zhize.li@kaust.edu.sa}

\icmlkeywords{Acceleration, Compressed Gradient Descent, Distributed Optimization, Federated Optimization}

\vskip 0.3in
]



\printAffiliationsAndNotice{}  

\begin{abstract}
Due to the high communication cost in distributed and federated learning problems, methods relying on compression of communicated messages are becoming increasingly popular. While in other contexts the best performing gradient-type  methods invariably rely on some form of acceleration/momentum to reduce the number of iterations, there are no methods which combine the benefits of both gradient compression and acceleration. 
In this paper, we remedy this situation and propose the first {\em accelerated compressed gradient descent (ACGD)} methods. In the single machine regime, we prove that ACGD enjoys the rate $O\Big((1+\omega)\sqrt{\frac{L}{\mu}}\log \frac{1}{\epsilon}\Big)$ for $\mu$-strongly convex problems and $O\Big((1+\omega)\sqrt{\frac{L}{\epsilon}}\Big)$ for convex problems, respectively, where $\omega$ is the compression parameter. 
Our results improve upon the existing non-accelerated rates $O\Big((1+\omega)\frac{L}{\mu}\log \frac{1}{\epsilon}\Big)$ and $O\Big((1+\omega)\frac{L}{\epsilon}\Big)$, respectively, and recover the optimal rates of accelerated gradient descent as a special case when no compression ($\omega=0$) is applied. 
We further propose a distributed variant of ACGD (called ADIANA) and prove the convergence rate $\widetilde{O}\Big(\omega+\sqrt{\frac{L}{\mu}}+\sqrt{\big(\frac{\omega}{n}+\sqrt{\frac{\omega}{n}}\big)\frac{\omega L}{\mu}}\Big)$, where $n$ is the number of devices/workers and $\widetilde{O}$ hides the logarithmic factor $\log \frac{1}{\epsilon}$. This improves upon the previous best result $\widetilde{O}\Big(\omega + \frac{L}{\mu}+\frac{\omega L}{n\mu} \Big)$ achieved by the DIANA method of \citet{DIANA}. Finally, we conduct several experiments on real-world datasets which corroborate  our theoretical results and confirm the practical superiority of our accelerated methods.
\end{abstract}

\section{Introduction}
\label{sec:intro_new}

With  the proliferation of edge devices such as mobile phones, wearables and smart home devices comes an increase in the amount of data rich in potential information which can be mined for the benefit of the users. One of the approaches of turning the raw data into information is via federated learning \citep{FEDLEARN, FL2017-AISTATS}, where typically a single global supervised model is trained in a massively distributed manner over a network of heterogeneous devices.

Training supervised federated learning models is typically performed by solving an optimization problem of the form
\begin{equation}\label{eq:prob}
\min_{x\in\R^d} \Big\{P(x)\eqdef\frac{1}{n}\sum \limits_{i=1}^{n}f_i(x) + \psi(x)\Big\},
\end{equation}
where $f_i:\R^d\rightarrow \R$ is smooth loss  associated with data stored on device $i$ and  $\psi: \R^d \rightarrow \R\cup \{+\infty\}$ is a relatively simple but possibly nonsmooth regularizer.

In distributed learning in general, and federated learning in particular, communication of messages across a network forms the bottleneck of the training system. It is thus very important to devise novel  strategies for reducing the number of communication rounds. Two of the most common strategies are i) local computations~\citep{COCOA+journal, localSGD-Stich, localSGD-AISTATS2020} and 
ii)  communication compression~\citep{Seide2015:1bit, alistarh2017qsgd, tonko, Cnat}. The former is used to perform more 
local computations on each device  before communication and subsequent model averaging, hoping that this will reduce the total number of communications. The latter is used to reduce the size of communicated messages, saving precious time spent in each communication round, and hoping that  this will not increase the total number of communications. 

\subsection{Theoretical inefficiency of local methods} Despite their practical success, local methods are poorly understood  and there is much to be discovered. For instance, there exist no theoretical results which would suggest that any local method (e.g., local gradient descent (GD) or local SGD) can achieve better communication complexity than its standard non-local variant (e.g., GD, SGD). In fact, until recently, no complexity results existed for local SGD in environments with  heterogeneous data~\citep{localGD, localSGD-AISTATS2020}, a key regime in federated learning settings~\citep{FL_survey_2019}. In the important regime when all participating devices compute full gradients based on their local data, the recently proposed stochastic controlled averaging (SCAFFOLD) method~\citep{SCAFFOLD}  offers no improvement on the number of communication as the number of local steps grows despite the fact that this is a rather elaborate method combining local stochastic gradient descent  with control variates for reducing the model drift among clients.

\subsection{Methods with compressed communication}
However, the situation is much brighter with methods employing communication compression. Indeed, several recent theoretical results suggest that by combining an appropriate (typically randomized) compression operator with a suitably designed gradient-type method, one can obtain  improvement in the the total communication complexity over comparable baselines not performing any compression. For instance, this is the case for distributed compressed gradient descent (CGD)~\citep{alistarh2017qsgd, khirirat2018distributed, Cnat, li2020unified} and distributed CGD methods which employ variance reduction to tame the variance introduced by compression~\citep{SEGA, DIANA, DIANA2, GJS, li2020unified}. 

While in the case of CGD compression leads to a decrease in the size of communicated messages per communication round, it leads to  an increase in the number of communications. Yet, certain compression operators, such as natural dithering~\citep{Cnat}, were shown to be better than no compression in terms of the overall communication complexity. 

The variance-reduced CGD method DIANA~\citep{DIANA, DIANA2} enjoys even better behavior: the number of communication rounds for this method is unaffected up to a certain level of compression when the variance induced by compression is smaller than a certain threshold. This threshold can be very large in practice, which means that massive reduction is often possible in the number of communicated bits without this having any adverse effect on the number of communication rounds. 

Recall that a function $f:\R^d\to \R$ is $L$-smooth or has $L$-Lipschitz continuous gradient (for $L> 0$) if 
\begin{equation}\label{def_smoothness}
\|\nabla f(x)-\nabla f(y)\| \le L \|x-y\|,
\end{equation}
and $\mu$-strongly convex (for $\mu\geq 0$) if
\begin{equation}\label{def_strongconvex}
f(x) - f(y) - \inner{\nabla f(y)}{x-y} \geq \frac{\mu}{2}\ns{x-y}
\end{equation}
for all $x,y \in \R^d$. The $\mu=0$ case corresponds to the standard convexity.

In particular, for $L$-smooth and $\mu$-strongly convex
$f$ with $n$ machines, DIANA enjoys the  iteration bound $O\left( \big(\omega + \tfrac{L}{\mu} +\tfrac{\omega }{n} \tfrac{L}{\mu}\big) \log \tfrac{1}{\epsilon}\right),$
where $\frac{L}{\mu}$ is the condition number and $\omega$ is the compression parameter (see Definition \ref{def_compression}).  If $\omega=0$, which corresponds to no compression, DIANA recovers the $\cO\left(\frac{L}{\mu} \log \frac{1}{\epsilon} \right)$ rate of gradient descent. On the other hand, as long as  $\omega = O\left(\min \big\{\frac{L}{\mu},n \big\}\right)$, the rate is still  $\cO\left(\frac{L}{\mu} \log \frac{1}{\epsilon} \right)$, which shows that DIANA is able to retain the same number of communication rounds as gradient descent and yet save on bit transmission in each round. The higher $\omega$ is allowed to be, the more compression can be applied.

\section{Contributions}

{\footnotesize
	
	\begin{table*}[t]
		\centering
		\caption{Convergence results for the special case with $n=1$ device (i.e., problem \eqref{eq:min_f})}\label{table:1}
		\renewcommand{\arraystretch}{1.7}
		\begin{tabular}{|c|c|c|}
			\hline
			Algorithm & $\mu$-strongly convex $f$ & convex $f$\\
			\hline
			\makecell{ Compressed Gradient Descent \\ 
				(CGD \citep{khirirat2018distributed})}
			& $O\left((1+\omega)\frac{L}{\mu}\log \frac{1}{\epsilon}\right)$ 
			&  $O\left((1+\omega)\frac{L}{\epsilon}\right)$ \\ \hline
			\rowcolor{bgcolor}
			ACGD (this paper) 
			& $O\left((1+\omega)\sqrt{\frac{L}{\mu}}\log \frac{1}{\epsilon}\right)$ 
			& $O\left((1+\omega)\sqrt{\frac{L}{\epsilon}}\right)$ \\ \hline
		\end{tabular}

		\vspace{1mm}
		\centering
		\caption{Convergence results for the general case with $n$ devices (i.e., problem \eqref{eq:prob}). Our  results are always better than previous results. }\label{table:2}
		\renewcommand{\arraystretch}{1.7}
		\begin{tabular}{|c|c|c|}
			\hline
			Algorithm 
			& \makecell{$n \leq  \omega$ \\(few devices or high compression)}
			& \makecell{$n>\omega$ \\(lots of devices or low compression)}\\
			\hline
			\makecell{ Distributed CGD \\ 
				(DIANA \citep{DIANA})}
			& $O\left(\omega \big(1+\frac{L}{n\mu}\big)\log \frac{1}{\epsilon}\right)$ 
			& $O\left(\big(\omega+\frac{L}{\mu}\big)\log \frac{1}{\epsilon}\right)$  \\ \hline 
			\rowcolor{bgcolor} 
			ADIANA (this paper) 
			& $O\left(\omega\Big(1+\sqrt{\frac{L}{n\mu}}\Big)\log \frac{1}{\epsilon}\right)$
			& $O\left(\Big(\omega+\sqrt{\frac{L}{\mu}}+\sqrt{\sqrt{\frac{\omega}{n}}\frac{\omega L}{\mu}}\Big)\log \frac{1}{\epsilon}\right)$ \\ \hline
		\end{tabular}
		
	\end{table*}
}

Discouraged by the lack of theoretical results suggesting that local methods indeed help to reduce the number of communications, and encouraged by the theoretical success of CGD methods, in this paper we seek to enhance CGD methods with a mechanism which, unlike local updating, can provably lead to a decrease of the number of communication rounds.  

	{\em What mechanism could achieve further improvements?}

In the world of deterministic gradient methods, one technique for  such a reduction is well known: {\em Nesterov acceleration / momentum}~\citep{nesterov1, NesterovBook}. In case of stochastic gradient methods, the accelerated method Katyusha~\citep{allen2017katyusha} achieves the optimal rate for strongly convex problems, and the unified accelerated method Varag~\citep{Zhize2019unified} achieves the optimal rates for convex problems regardless of the strong convexity. See also  \citep{L-SVRG,  L-SVRG-AS} for some enhancements. 
Essentially all state-of-the-art methods for training deep learning models, including Adam~\citep{ADAM}, rely on the use of momentum/acceleration in one form or another, albeit lacking in  theoretical support. 

	{\em However, the successful combination of gradient compression and acceleration/momentum has so far remained elusive, and to the best of our knowledge, no algorithms supported with theoretical results exist in this space. Given the omnipresence of momentum in modern machine learning, this is surprising. }

We now summarize our key contributions:

\subsection{First combination of gradient compression and acceleration}

We develop the first gradient-type optimization methods provably  combining the  benefits of gradient compression and acceleration: i) ACGD (Algorithm~\ref{alg:acgd}) in the single device case, and ii) ADIANA (Algorithm~\ref{alg:ADIANA}) in the distributed case.

\subsection{Single device setting}

We first study the single-device setting, and design an accelerated CGD method (ACGD - Algorithm~\ref{alg:acgd}) for solving the unconstrained smooth minimization problem
\begin{equation}\label{eq:min_f} \min_{x\in \R^d} f(x)\end{equation}
in the regimes when $f$ is $L$-smooth and i) $\mu$-strongly convex, and ii)  convex. Our theoretical results are summarized in Table~\ref{table:1}.  In the strongly convex case, we improve the complexity of CGD~\citep{khirirat2018distributed} from $O\Big((1+\omega)\frac{L}{\mu}\log \frac{1}{\epsilon}\Big)$ to 
$O\Big((1+\omega)\sqrt{\frac{L}{\mu}}\log \frac{1}{\epsilon}\Big).$ 
In the convex case, the improvement is from $O\Big((1+\omega)\frac{L}{\epsilon}\Big)$ to 
$O\Big((1+\omega)\sqrt{\frac{L}{\epsilon}}\Big),$
where $\omega\geq 0$ is the compression parameter (see Definition \ref{def_compression}).

\subsection{Distributed setting}

We further study the distributed setting with $n$ devices/nodes and focus on problem
\eqref{eq:prob} in its full generality, i.e.,
\begin{equation*}
\min_{x\in\R^d} \Big\{P(x)\eqdef\frac{1}{n}\sum \limits_{i=1}^{n}f_i(x) + \psi(x)\Big\},
\end{equation*}

The presence of multiple nodes ($n>1$) and of the regularizer $\psi$ poses additional challenges. In order to address them, we need to not only combine acceleration and compression, but also introduce a DIANA-like variance reduction mechanism to remove  the variance introduced by the compression operators.

In particular, we have developed an accelerated variant of the DIANA method for solving the general problem \eqref{eq:prob}, which we call ADIANA (Algorithm~\ref{alg:ADIANA}). 
The comparison of complexity results between ADIANA and DIANA is summarized in Table~\ref{table:2}. 

Note that our results always improve upon the non-accelerated DIANA method.  
Indeed, in the regime when the compression parameter $\omega$ is larger than the number of nodes $n$, we improve the DIANA rate $O\Big(\omega\big(1+\frac{L}{n\mu}\big)\log \frac{1}{\epsilon}\Big)$ to 
$ O\left(\omega\Big(1+\sqrt{\frac{L}{n\mu}}\Big)\log \frac{1}{\epsilon}\right).$
On the other hand, in the regime when $\omega <n$, we improve the DIANA rate $O\left( \big(\omega+\frac{L}{\mu} \big)\log \frac{1}{\epsilon}\right)$ to 
$ O\left(\Big(\omega+\sqrt{\frac{L}{\mu}}+\sqrt{\sqrt{\frac{\omega}{n}}\frac{\omega L}{\mu}}\Big)\log \frac{1}{\epsilon}\right).$ Our rate is better since $\omega + \frac{L}{\mu} \geq 2\sqrt{\frac{\omega L}{\mu}}$ and $\sqrt{\frac{\omega}{n}} <1$ (note that $\omega < n$). 

Note that if $\omega \leq n^{1/3}$, which is more often true in federated learning as the number of devices in federated learning is typically very large,  our ADIANA result reduces to $O\left(\Big(\omega+\sqrt{\frac{L}{\mu}}\Big)\log \frac{1}{\epsilon}\right).$ In particular, if $\omega = O\Big(\min \Big\{n^{1/3} , \sqrt{\frac{L}{\mu}} \Big\}\Big),$ then the communication round is $O\left(\sqrt{\tfrac{L}{\mu}} \log \frac{1}{\epsilon}\right)$; the same as that of non-compressed accelerated gradient descent (AGD) \citep{NesterovBook}. It means that ADIANA benefits from cheaper communication due to compression \emph{for free} without hurting the convergence rate (i.e., the communication rounds are the same), and is therefore better suited for federated optimization.

\section{Randomized Compression Operators}
\label{sec:compr}

We now introduce the notion of a randomized compression operator which is used to compress the gradients.

\begin{definition}[Compression operator]\label{def_compression}
	A randomized map $\calC: \R^d\mapsto \R^d$ is an $\omega$-compression operator  if  
	\begin{equation}\label{eq:comp}
	\EE[\calC(x)]=x,~~ \EE[\ns{\calC(x)-x}] \leq \omega\ns{x},~~ \forall x\in \R^d.
	\end{equation}
	In particular, no compression ($\calC(x)\equiv x$) implies $\omega=0$.
\end{definition}

Note that the conditions \eqref{eq:comp} require the compression operator to be unbiased and its variance uniformly bounded by a relative magnitude of the vector which we are compressing. 

\subsection{Examples}

We now give a few examples of randomized compression operators without attempting to be exhaustive.

\head{Example 1 (Random sparsification)}
Given $x\in \R^d$, the random-$k$ sparsification operator is defined by $$\calC(x)\eqdef \frac{d}{k} (\xi_k \odot x),$$ where 
$\odot$ denotes the Hadamard (element-wise) product and $\xi_k\in \{0,1\}^d$ is a uniformly random binary vector with $k$ nonzero entries ($\n{\xi_k}_0=k$).
This random-$k$ sparsification operator $\calC$ satisfies \eqref{eq:comp} with $\omega =\frac{d}{k}-1$.
Indeed, no compression $k=d$ implies $\omega=0$.

\head{Example 2 (Quantization)}
Given $x\in \R^d$, the ($p,s$)-quantization operator  is defined by
$$\calC(x)\eqdef \text{sign}(x)\cdot \n{x}_p\cdot \frac{1}{s} \cdot \xi_s,$$ 
where $\xi_s\in\R^d$ is a random vector with 
$i$-th element $$\xi_s(i)\eqdef \begin{cases}
l+1, &\text {with probability } \frac{|x_i|}{\n{x}_p}s-l\\
l, &\text {otherwise}
\end{cases},$$ where the level $l$ satisfies $\frac{|x_i|}{\n{x}_p}\in [\frac{l}{s}, \frac{l+1}{s}]$. The probability is chosen so that $\EE[\xi_s(i)]=\frac{|x_i|}{\n{x}_p}s$.
This ($p,s$)-quantization operator $\calC$ satisfies \eqref{eq:comp} with $\omega = 2+\frac{d^{1/p}+d^{1/2}}{s}$.
In particular, QSGD \citep{alistarh2017qsgd} used $p=2$ (i.e., ($2,s$)-quantization) and proved that the expected sparsity of $\cC(x)$ is $\EE[\n{\calC(x)}_0] = O\big(s(s+\sqrt{d})\big)$.

\section{Accelerated CGD: Single Machine}
\label{sec:acc-cgd}
In this section, we study the special case of problem \eqref{eq:prob} with a single machine ($n=1$) and no regularizer ($\psi(x)\equiv 0$), i.e., problem \eqref{eq:min_f}:   \[\min_{x \in \R^d} f(x).\]

\subsection{The CGD algorithm}
First, we recall the update step in compressed gradient descent (CGD) method, i.e., \[x^{k+1} = x^k - \eta \calC(\nabla f(x^k)),\] where $\calC$ is a kind of $\omega$-compression operator defined in Definition~\ref{def_compression}.

As mentioned earlier, convergence results of CGD are $O\left((1+\omega)\frac{L}{\mu}\log \frac{1}{\epsilon}\right)$ for strongly convex problems and $O\left((1+\omega)\frac{L}{\epsilon}\right)$ for convex problems (see Table \ref{table:1}). The convergence proof for strongly convex problems, i.e., $O\left((1+\omega)\frac{L}{\mu}\log \frac{1}{\epsilon}\right)$,  can be found in \citep{khirirat2018distributed}. For completeness, we now establish a convergence result for convex problems, i.e.,  $O\left((1+\omega)\frac{L}{\epsilon}\right)$ since we did not find it in the literature.

\begin{theorem}\label{thm:cgdconvex}
	Suppose $f$ is convex with $L$-Lipschitz continuous gradient and the compression operator $\calC$ satisfies \eqref{eq:comp}.  Fixing the step size $\eta=\frac{1}{(1+\omega)L}$, the number of iterations performed by CGD to find an $\epsilon$-solution such that $$\EE[f(x^k)]-f(x^*)\leq \epsilon$$ is at most  $$ k=O\left(\frac{(1+\omega)L}{\epsilon}\right).$$
\end{theorem}

\subsection{The ACGD algorithm}

Note that in the non-compressed case $\omega=0$ (i.e., CGD is reduced to  standard GD), there exists   methods for obtaining accelerated convergence rates of $O\left(\sqrt{\frac{L}{\mu}} \log \frac{1}{\epsilon}\right)$ and $O\left( \sqrt{\frac{L}{\epsilon}}\right)$ for strongly convex and convex problems, respectively. However, no accelerated convergence  results exist for CGD methods. Inspired by Nesterov's accelerated gradient descent (AGD) method \citep{NesterovBook} and FISTA \citep{beck2009fast}, we propose the first accelerated compressed gradient descent (ACGD) method, described in Algorithm~\ref{alg:acgd}.

\begin{algorithm}[h]
	\caption{Accelerated CGD (ACGD)}
	\label{alg:acgd}
	\begin{algorithmic}[1]
		\REQUIRE 
		initial point $x^0$, $\{\et\}, \{\at\}, \{\bt\}, \{\st\}$, $p$\\
		\STATE $z^0=y^0=x^0$
		\FOR {$k=0,1,2,\ldots$}
		\STATE $\yt=\at \xt +(1-\at)\vt$ \alglabel{line:y}
		\STATE Compress gradient $\gt=\calC(\nabla f(\yt))$ \alglabel{line:compress}
		\STATE $\xp =\yt -\frac{\et}{p} \gt$  \alglabel{line:x}
		\STATE $\vp=\frac{1}{\st}\xp + \left(\frac{1}{p}-\frac{1}{\st}\right)\xt$ \\
					\qquad\qquad$+ \left(1-\frac{1}{p}\right)(1-\bt)\vt
									+\left(1-\frac{1}{p}\right)\bt \yt $  \alglabel{line:v}
		\ENDFOR
	\end{algorithmic}
\end{algorithm}

\subsection{Convergence theory}
Our accelerated convergence results for ACGD (Algorithm \ref{alg:acgd}) are stated in Theorems~\ref{thm:acgdconvex} and \ref{thm:acgdstrongconvex}, formulated next.

\begin{theorem}[ACGD: convex case]\label{thm:acgdconvex}
	Let $f$ be convex with $L$-Lipschitz continuous gradient and let the compression operator $\calC$ satisfy \eqref{eq:comp}. Choose the parameters in ACGD (Algorithm~\ref{alg:acgd}) as follows: 
	\begin{align*}
		&\et \equiv \frac{1}{L}, \quad p=1+\omega,\\
		&\at =\frac{k}{k+2}, \quad  \bt \equiv 0, \quad \st = \frac{2p}{k+2}.
	\end{align*}
	 Then the number of iterations performed by ACGD  to find an $\epsilon$-solution such that \[\EE[f(x^k)]-f(x^*)\leq \epsilon\] is at most \[k=O\left((1+\omega)\sqrt{\frac{L}{\epsilon}}\right).\]
\end{theorem}

\begin{theorem}[ACGD: strongly convex case]\label{thm:acgdstrongconvex}
	Let $f$ be $\mu$-strongly convex with $L$-Lipschitz continuous gradient and let the compression operator $\calC$ satisfy \eqref{eq:comp}.Choose the parameters in ACGD (Algorithm~\ref{alg:acgd}) as follows: 
	\begin{align*}
	&\et \equiv \frac{1}{L}, \quad p=1+\omega,\\
	&\at \equiv \frac{p}{p+\sqrt{\mu/L}}, \quad  \bt \equiv \frac{\sqrt{\mu/L}}{p}, \quad \st \equiv \sqrt{\frac{\mu}{L}}.
	\end{align*}
 	Then the number of iterations performed by ACGD to find an $\epsilon$-solution such that \[\EE[f(x^k)]-f(x^*)\leq \epsilon\] (or $\EE[\ns{x^k-x^*}]\leq \epsilon$)  is at most \[k=O\left((1+\omega)\sqrt{\frac{L}{\mu}}\log \frac{1}{\epsilon}\right).\]
\end{theorem}

In the non-compressed case $\omega=0$ (i.e., $\calC(x)\equiv x$), our results recover the standard optimal rates of accelerated gradient descent. Further, if we consider the random-$k$ sparsification compression operator,  ACGD can be seen as a variant of accelerated randomized coordinate descent~\citep{RCDM}. Our results  recover the optimal results of accelerated randomized  coordinate descent method \citep{NUACDM, ACD2019} under the same standard smoothness assumptions. 

\subsection{Proof sketch} The following  lemma which demonstrates  improvement in one iteration plays a key role in our analysis.

\begin{lemma}\label{lem:key}
	If parameters $\{\et\}, \{\at\}, \{\bt\}$, $\{\st\}$ and $p$ satisfy 
	$\at = \frac{1 - \st/p}{1-\bt\st/p}$,
	$\bt \leq \min\{\frac{\mu\et}{\st p},1\}$, $p\geq \frac{(1+L\et)(1+\omega)}{2}$ and the compression operator $\calC^k$ satisfies \eqref{eq:comp},
	then we have for any iteration $k$ of ACGD, and for all $ x \in \R^d$, 
	\begin{align*}
	&\frac{2\et}{\st^2}\EE[f(\xp) -f(x)]+ \EE[\ns{\vp-x}] \\
	&\leq \left(1-\frac{\st}{p}\right)\frac{2\et}{\st^2}\left(f(\xt) -f(x)\right) +(1-\bt)\ns{\vt-x},  
	\end{align*}
	where the expectation is with respect to the randomness of compression operator sampled at iteration $k$.
\end{lemma}

The proof of Theorems~\ref{thm:acgdconvex} and \ref{thm:acgdstrongconvex} can be derived (i.e., plug into the specified parameters ($\{\et\}, \{\at\}, \{\bt\}$, $\{\st\}$ and $p$) and collect all iterations) from Lemma \ref{lem:key}. The detailed proofs can be found in the appendix.

\section{Accelerated CGD: Distributed Setting}
\label{sec:acc-dcgd}

We now turn our attention to the general distributed case, i.e., problem \eqref{eq:prob}:
\begin{equation*}
\min_{x\in\R^d} \Big\{P(x)\eqdef\frac{1}{n}\sum_{i=1}^{n}f_i(x) + \psi(x)\Big\}.
\end{equation*}
The presence of multiple nodes ($n>1$) and of the regularizer $\psi$ poses additional challenges.

\subsection{The ADIANA algorithm}

We now propose an \emph{accelerated} algorithm for solving problem 
\eqref{eq:prob}. 
Our method combines both acceleration and variance reduction, and hence can be seen as an accelerated version of DIANA \citep{DIANA,DIANA2}. Therefore, we  call our method ADIANA (Algorithm~\ref{alg:ADIANA}).
In this case, each machine/agent computes its local gradient (e.g., $\nabla f_i(x^k)$) and a shifted version thereof is compressed and sent to the server. The server subsequently aggregates all received messages, to form a stochastic gradient estimator $g^k$ of $\frac{1}{n}\sum_i \nabla f_i(x^k)$, and then performs a proximal step. 
The shift terms $h_i^k$ are adaptively changing throughout the iterative process, and have the role of reducing the variance introduced by compression. If no compression is used, we may simply set the shift terms to be $h_i^k=0$ for all $i,k$. 

Our method was inspired by \citet{DIANA}, who first studied variance reduction for CGD methods for a specific ternary compression operator, and  \citet{DIANA2} who studied the general class of $\omega$-compression operators we also study here. 
However, we had to make certain modifications to make variance-reduced compression work in the accelerated case since both of them were studied in the non-accelerated case.
Besides, our method adopts a randomized update rule for the auxiliary vectors $w^k$ which simplifies the algorithm and analysis, resembling the workings of the loopless SVRG method proposed by \citet{L-SVRG}.

\begin{algorithm}[th]
	\caption{Accelerated DIANA (ADIANA)}
	\label{alg:ADIANA}
	\begin{algorithmic}[1]
		\REQUIRE 
		initial point $x^0$, $\{h_i^0\}_{i=1}^n$, $h^0=\frac{1}{n}\sum_{i=1}^{n}h_i^0$, parameters $\eta, \theta_1, \theta_2, \alpha, \beta, \gamma, p$ \\
		\STATE $z^0=y^0=w^0=x^0$
		\FOR{$k = 0,1,2,\ldots$}
		\STATE $x^k = \theta_1 z^k + \theta_2 w^k + (1 - \theta_1 - \theta_2)y^k$ \alglabel{line:xk}
		\STATE {\bf{for all machines $i= 1,2,\ldots,n$ do in parallel}}
		\STATE Compress shifted local gradient $\calC_i^k(\nabla f_i(x^k) - h_i^k)$ and send to the server
		\STATE Update local shift $h_i^{k+1}=h_i^k+\alpha \cC_i^k(\nabla f_i(w^k) - h_i^k)$
		\STATE {\bf{end for}}
		\STATE Aggregate received compressed gradient information\\
		\qquad $g^k = \frac{1}{n}\sum \limits_{i=1}^n \cC_i^k(\nabla f_i(x^k) - h_i^k) + h^k$\\
		\qquad $h^{k+1} = h^k + \alpha  \frac{1}{n}\sum \limits_{i=1}^n \cC_i^k(\nabla f_i(w^k) - h_i^k)$
		\STATE Perform update step\\ 
		\qquad $y^{k+1} = \prox_{\eta\psi}(x^k - \eta g^k)$ \alglabel{line:yp}
		\STATE $z^{k+1} = \beta z^k + (1-\beta)x^k + \frac{\gamma}{\eta} (y^{k+1} - x^k)$
		\STATE $w^{k+1} = \begin{cases}
		y^k, &\text {with probability } p\\
		w^k, &\text {with probability } 1-p
		\end{cases}$ \alglabel{line:prob}
		\ENDFOR
	\end{algorithmic}
\end{algorithm}

\subsection{Convergence theory}

Our main convergence result for ADIANA (Algorithm \ref{alg:ADIANA})  is formulated in Theorem \ref{thm:accdcgd}. We focus on the strongly convex setting.

\begin{theorem}\label{thm:accdcgd}
	Suppose $f$ is $\mu$-strongly convex and that the functions $f_i$ have $L$-Lipschitz continuous gradient for all $i$. Further, let  the compression operator $\calC$ satisfy \eqref{eq:comp}. Choose the ADIANA (Algorithm \ref{alg:ADIANA}) parameters as follows: 
	\begin{align*}
	&\eta =\min \left\{\frac{1}{2L}, \frac{n}{64\omega(2p(\omega+1)+1)^2L} \right\}, \\
	&\theta_1 = \min\left\{\frac{1}{4},~ \sqrt{\frac{\eta\mu}{p}} \right\}, \quad \theta_2 = \frac{1}{2},\\
	&\alpha = \frac{1}{\omega+1}, \quad  \beta = 1-\gamma \mu, 
			\quad \gamma=\frac{\eta}{2(\theta_1 + \eta\mu)}, \\
	&p=\min\left\{1, \frac{\max\{1,\sqrt{\frac{n}{32\omega}}-1\}}{2(1+\omega)}\right\}.
	\end{align*}
	Then the number of iterations performed by ADIANA  to find an $\epsilon$-solution such that \[\EE[\ns{z^k-x^*}]\leq \epsilon\]  is at most
	\begin{equation*}
	k=
	\begin{cases}
	O\left(\left[\omega + \omega\sqrt{\frac{L}{n \mu}}\ \right]\log \frac{1}{\epsilon}  \right),& n \leq \omega,\\
	O\left(\left[\omega  + \sqrt{\frac{L}{\mu}} + \sqrt{\sqrt{\frac{\omega}{n}} \frac{\omega L}{\mu}}\ \right]\log \frac{1}{\epsilon}  \right),&  n>\omega. 
	\end{cases}
	\end{equation*}
\end{theorem}

As we have explained in the introduction, the above rate is vastly superior to that of non-accelerated distributed CGD methods, including  that of DIANA.

\subsection{Proof sketch}

In the proof, we use the following notation:
\begin{align}
\cZ^k &\eqdef \norm{z^k - x^*}^2, \label{def:z}\\
\cY^k &\eqdef P(y^k) - P(x^*), \label{def:y}\\
\cW^k &\eqdef P(w^k) - P(x^*), \label{def:w}\\ 
\cH^k &\eqdef \frac{1}{n}\sum_{i=1}^n \norm{h_i^k - \nabla f_i(w^k)}^2. \label{def:h}
\end{align}

We first present a key technical lemma which plays a similar role to that of Lemma~\ref{lem:key}.

\begin{lemma}\label{lem:keyaccdcgd}
	If the parameters satisfy 
	$\eta \leq \frac{1}{2L}, \theta_1 \leq \frac{1}{4}, \theta_2 =\frac{1}{2},  \gamma=\frac{\eta}{2(\theta_1 + \eta\mu)}$  and $\beta = 1-\gamma \mu$, 
	then we have for any iteration $k$, 
	\begin{align}
	&\frac{2\gamma\beta}{\theta_1}\E{\cY^{k+1}} +  \E{\cZ^{k+1}} \notag\\
	& \leq
	(1 - \theta_1 - \theta_2)\frac{2\gamma\beta}{\theta_1}\cY^k
	+
	\beta\cZ^k  \notag\\
	&\qquad +
	2\gamma\beta\frac{\theta_2}{\theta_1}\cW^k
	+
	\frac{\gamma\eta}{\theta_1}\E{\norm{g^k - \nabla f(x^k)}^2} \label{eq:wandg}\\
	&\qquad -
	\frac{\gamma}{4Ln\theta_1}\sum_{i=1}^n\norm{\nabla f_i(w^k) - \nabla f_i(x^k)}^2 \label{eq:lem2w} \\
	&\qquad -
	\frac{\gamma}{8Ln\theta_1}\sum_{i=1}^n\norm{\nabla f_i(y^k) - \nabla f_i(x^k)}^2. \label{eq:lem2y}  
	\end{align}
\end{lemma}

Theorem~\ref{thm:accdcgd} can be proved  by combing the above lemma with three additional Lemmas: Lemma~\ref{lem:keyw}, \ref{lem:keyg} and \ref{lem:keyh}, which we present next. In view of the presence of $\cW^k$ in \eqref{eq:wandg}, the following result is useful as it allows us to add $\cW^{k+1}$ into the Lyapunov function.

\begin{lemma}\label{lem:keyw}
	According to Line \ref{line:prob} of Algorithm \ref{alg:ADIANA} and Definition \eqref{def:y}--\eqref{def:w}, we have
	\begin{equation*}
	\E{\cW^{k+1}} = (1-p)\cW^k + p\cY^k.
	\end{equation*}
\end{lemma}

To cancel the term $\E{\norm{g^k - \nabla f(x^k)}^2}$ in \eqref{eq:wandg}, we use the defining property of compression operator (i.e., \eqref{eq:comp}) :
\begin{lemma}\label{lem:keyg}
	If the compression operator $\calC$ satisfies \eqref{eq:comp}, we have
	\begin{align}
	&\E{\norm{g^k - \nabla f(x^k)}^2}  \notag\\
	&\leq
	\frac{2\omega}{n^2}\sum_{i=1}^n \norm{\nabla f_i(w^k) - \nabla f_i(x^k)}^2
	+
	\frac{2\omega}{n}\cH^k \,.\label{eq:appearh}
	\end{align}
\end{lemma}

Note that the bound on variance obtained above introduces an additional term $\cH^k$ (see \eqref{eq:appearh}). We will therefore add the terms $\cH^{k+1}$ into the Lyapunov function as well.
\begin{lemma}\label{lem:keyh}
	If $\alpha \leq \frac{1}{\omega+1}$, we have
	\begin{align*}
	\E{\cH^{k+1}} &\leq \left(1-\frac{\alpha}{2}\right)\cH^k  \notag\\
		& + \left(1 + \frac{2p}{\alpha} \right)\frac{2p}{n}\sum_{i=1}^n\norm{\nabla f_i(w^k) - \nabla f_i(x^k)}^2   \notag\\
		&+ \left(1 + \frac{2p}{\alpha} \right)\frac{2p}{n}\sum_{i=1}^n\norm{\nabla f_i(y^k) - \nabla f_i(x^k)}^2. 
	\end{align*}
\end{lemma}
Note that the terms $\sum_{i=1}^n\norm{\nabla f_i(w^k) - \nabla f_i(x^k)}^2$ and 
$\sum_{i=1}^n\norm{\nabla f_i(y^k) - \nabla f_i(x^k)}^2$ in Lemma~\ref{lem:keyh} and \eqref{eq:appearh} can be cancelled by \eqref{eq:lem2w} and \eqref{eq:lem2y} by choosing the parameters appropriately.

Finally, it is not hard to obtain the following key inequality for the Lyapunov function by plugging Lemmas \ref{lem:keyw}-\ref{lem:keyh} into our key Lemma~\ref{lem:keyaccdcgd}:
\begin{align}
&\EE\left[c_1\cY^{k+1}+ c_2\cZ^{k+1} +c_3\cW^{k+1}+c_4\cH^{k+1}\right]  \notag\\
&\leq (1-c_5) \left(c_1\cY^{k}+ c_2\cZ^{k} +c_3\cW^{k}+c_4\cH^{k}\right). \label{eq:lyapunov}
\end{align}
Above, the constants $c_1,\ldots, c_5$ are related to the algorithm parameters $\eta, \theta_1, \theta_2, \alpha, \beta, \gamma$ and $p$.
Finally, the proof of Theorem~\ref{thm:accdcgd} can be derived (i.e., plug into the specified parameters) from inequality~\eqref{eq:lyapunov}. The detailed proof can be found in the appendix.

\section{Experiments}
\label{sec:exp}

In this section, we demonstrate the performance of our accelerated method ADIANA (Algorithm \ref{alg:ADIANA}) and previous methods with different compression operators on the regularized logistic regression problem, 
\begin{equation*}
\min_{x\in\R^d} \bigg\{\frac{1}{n}\sum \limits_{i=1}^{n}\log \big(1 + \exp(-b_i a_i^\top x)\big) + \frac{\lambda}{2}\|x\|^2\bigg\},
\end{equation*}
where $\{a_i, b_i\}_{i\in[n]}$ are data samples. 

\paragraph{Data sets.} In our experiments we use four standard datasets, namely, \texttt{a5a}, \texttt{mushrooms}, \texttt{a9a} and \texttt{w6a} from the LIBSVM library. Some of the experiments are provided in the appendix.

\paragraph{Compression operators.} We use three different compression operators: random sparsification (see e.g.~\citep{stich2018sparsified}), random dithering (see e.g.~\citep{alistarh2017qsgd}), and natural compression (see e.g.~\citep{Cnat}). For random-$r$ sparsification, the number of communicated bits per iteration is $32r$,  and we choose $r = d/4$. For random dithering, we choose $s = \sqrt{d}$, which means the number of communicated bits per iteration is $2.8d+32$ \citep{alistarh2017qsgd}. For natural compression, the number of communicated bits  per iteration is $9d$ bits \citep{Cnat}. 

\paragraph{Parameter setting.} 
In our experiments, we use the theoretical stepsize and parameters for all the three algorithms: vanilla distributed compressed gradient descent (DCGD), DIANA \citep{DIANA}, and our ADIANA (Algorithm \ref{alg:ADIANA}). 
The default number of nodes/machines is $20$ and the regularization parameter $\lambda=10^{-3}$. The numerical results for different number of nodes can be found in the appendix. 
For the figures, we plot the relation of the optimality loss gap $f(x^k) - f(x^*)$ and the number of accumulated transmitted bits. 
The optimal value $f(x^*)$ for each case is obtained by getting the minimum of three uncompressed versions of ADIANA, DIANA, and DCGD for $100000$ iterations.

\subsection{Comparison with DIANA and DCGD}
\label{sec:compare-diana-dcgd}
In this subsection, we compare our ADIANA with DIANA and DCGD with three compression operators: random sparsification, random dithering, and natural compression in Figures~\ref{fig:a5a111} and \ref{fig:mushrooms111}. 

The experimental results indeed show that our ADIANA converges fastest for all three compressors, and natural compression uses the fewest communication bits than random dithering and random sparsification.
Moreover, because the compression error of vanilla DCGD is nonzero in general, DCGD can only converge to the neighborhood of the optimal solution while DIANA and ADIANA can converge to the optimal solution.

\subsection{Communication efficiency}
\label{sec:compare-comp-noncomp}
\begin{figure*}[h]
	\vspace{0cm}
	\centering
	\begin{tabular}{cccc}
		\includegraphics[width=0.32\linewidth]{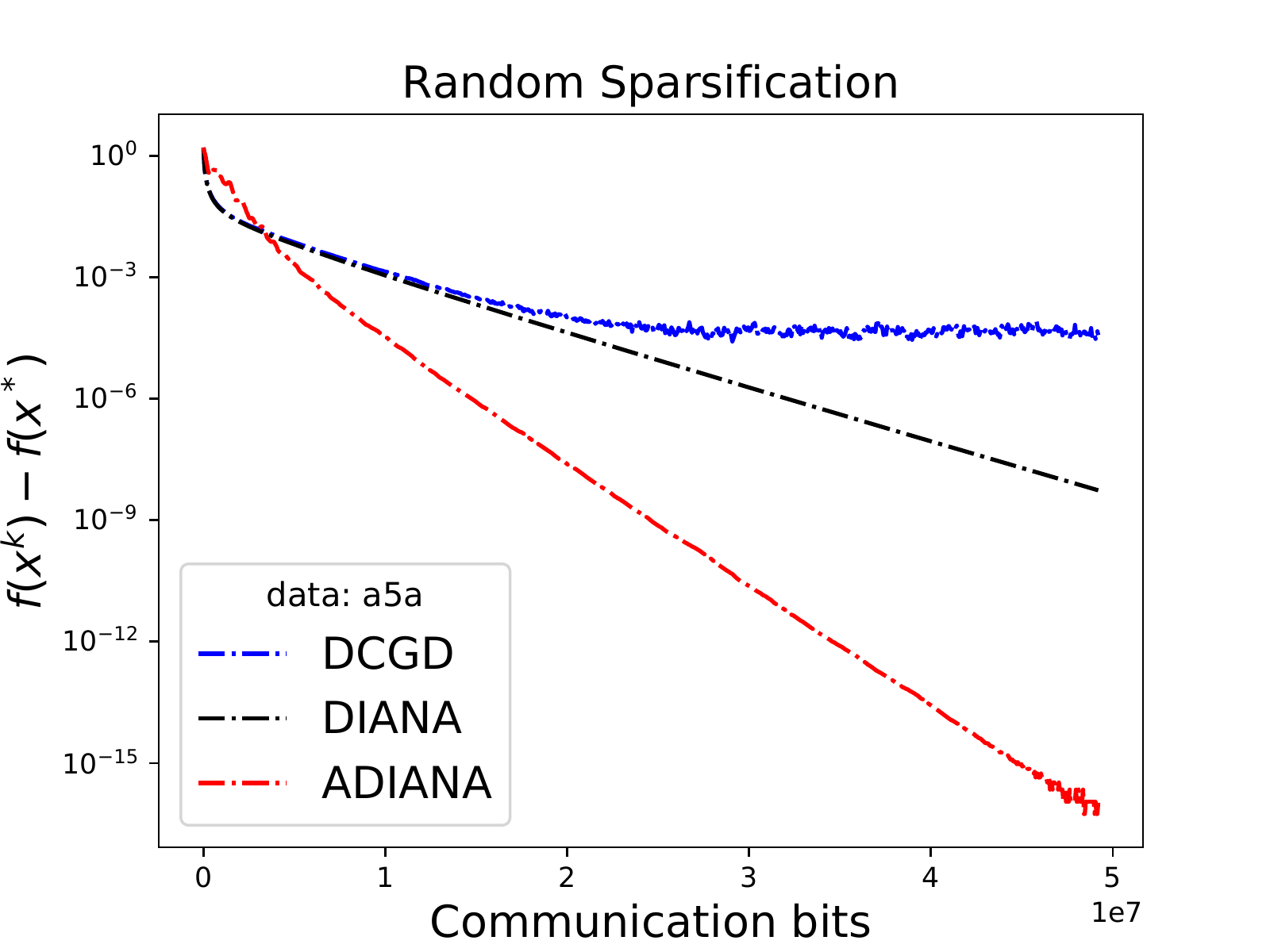}&
		\includegraphics[width=0.32\linewidth]{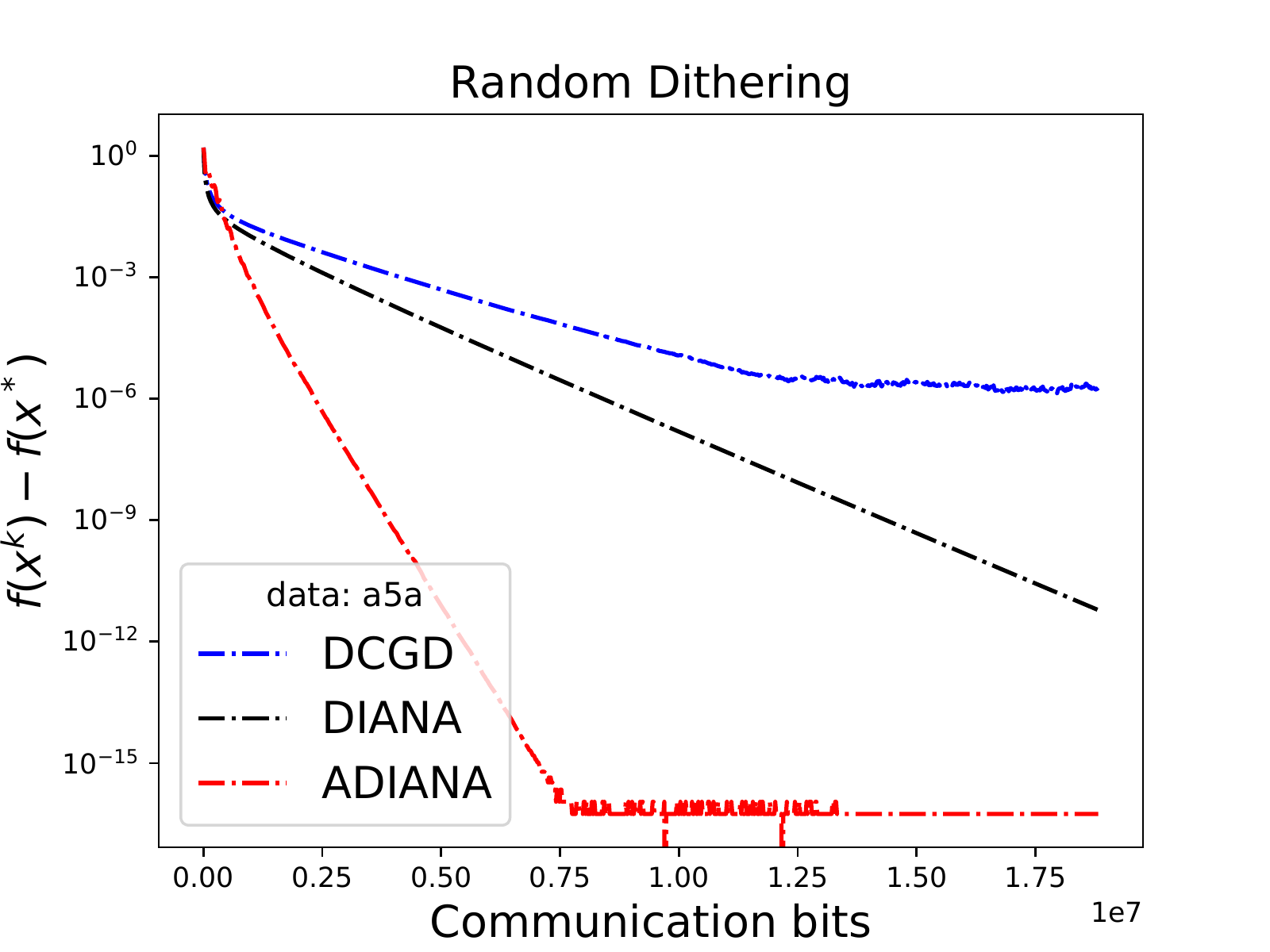}&
		\includegraphics[width=0.32\linewidth]{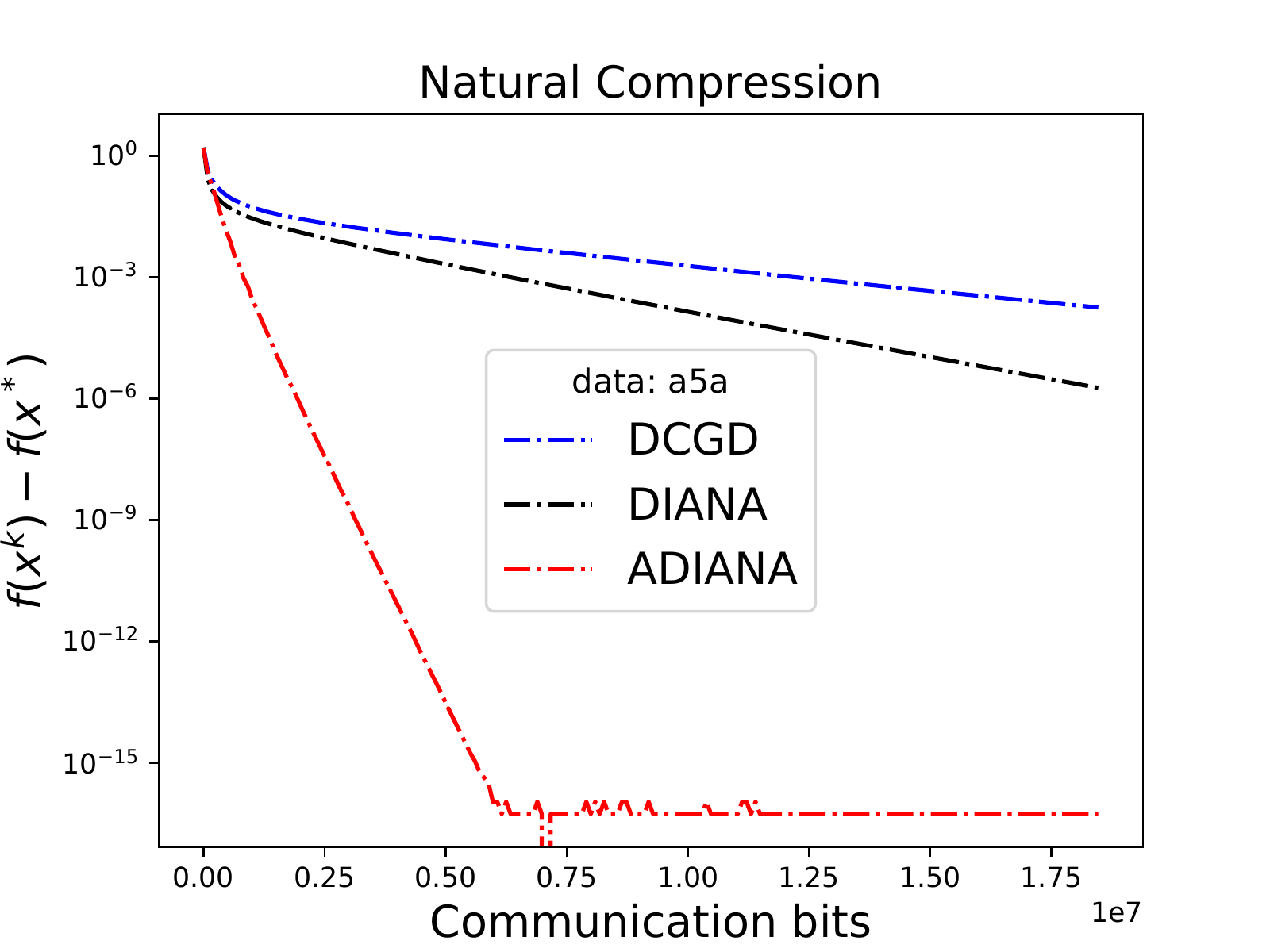}&
	\end{tabular}
	
	\caption{The communication complexity of different methods for three different compressors (random sparsification, random dithering and natural compression) on the \texttt{a5a} dataset.}
	\label{fig:a5a111}
	
	%
	\vspace{3mm}
	\centering
	\begin{tabular}{cccc}
		\includegraphics[width=0.32\linewidth]{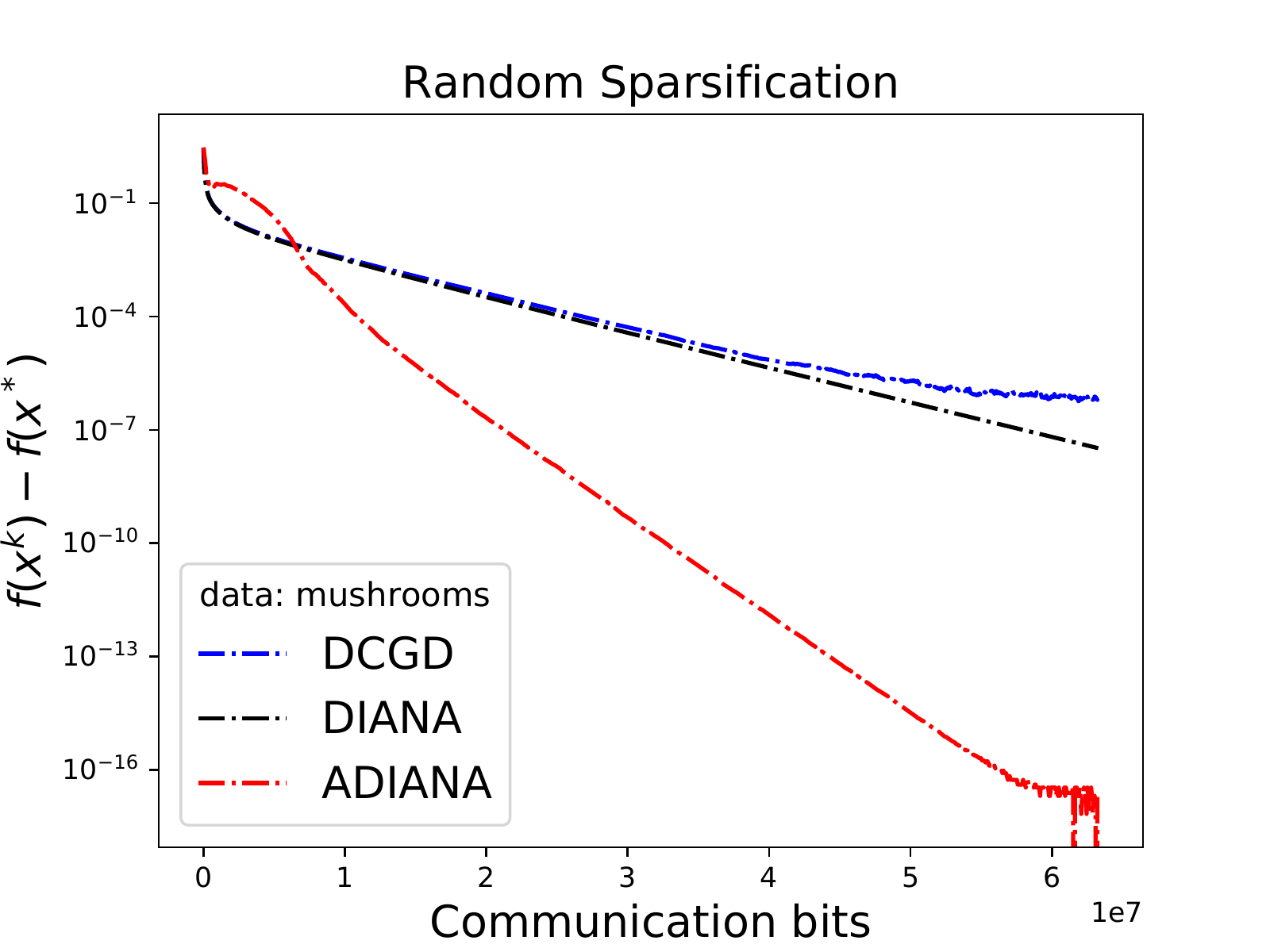}&
		\includegraphics[width=0.32\linewidth]{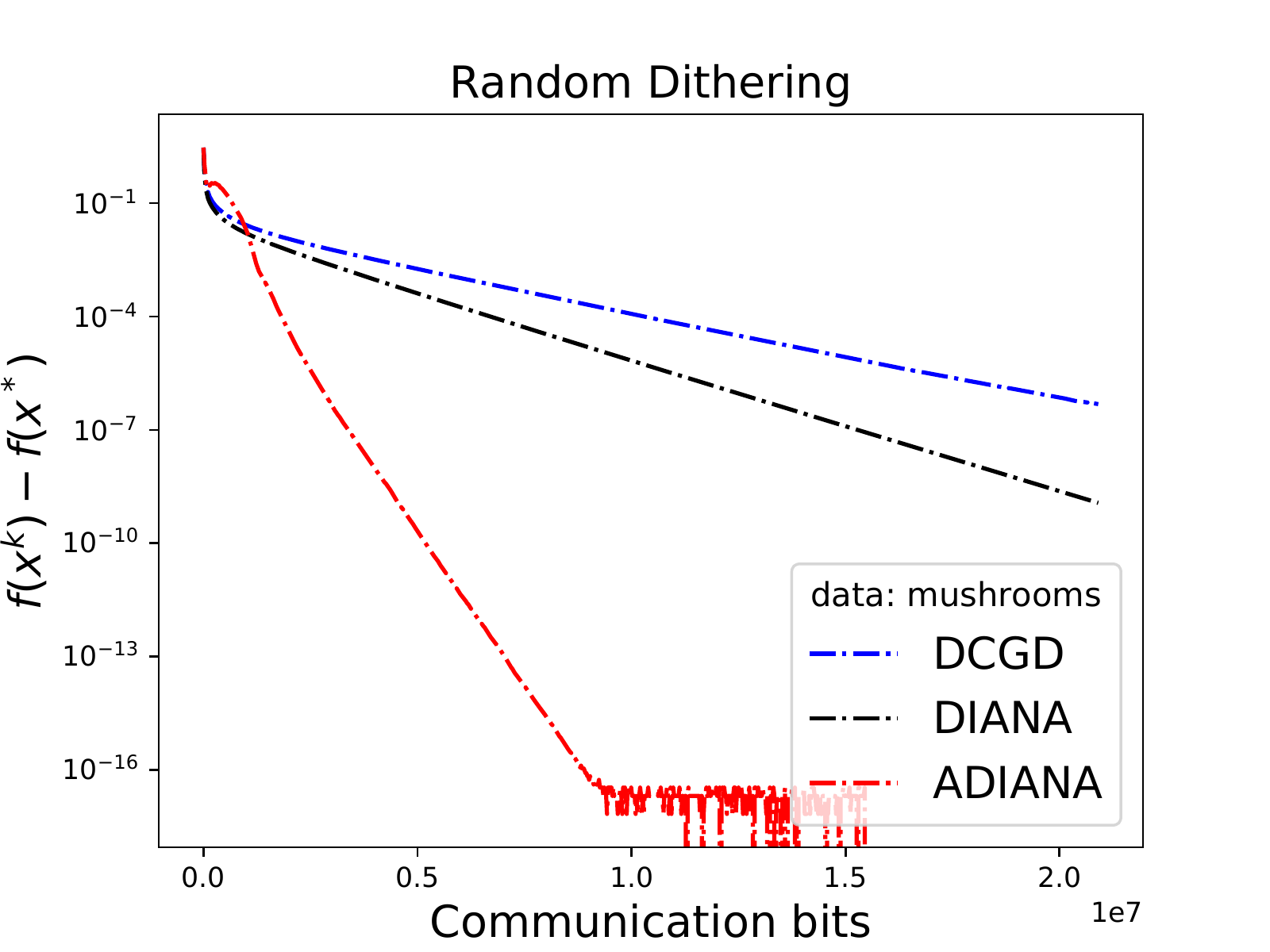}&
		\includegraphics[width=0.32\linewidth]{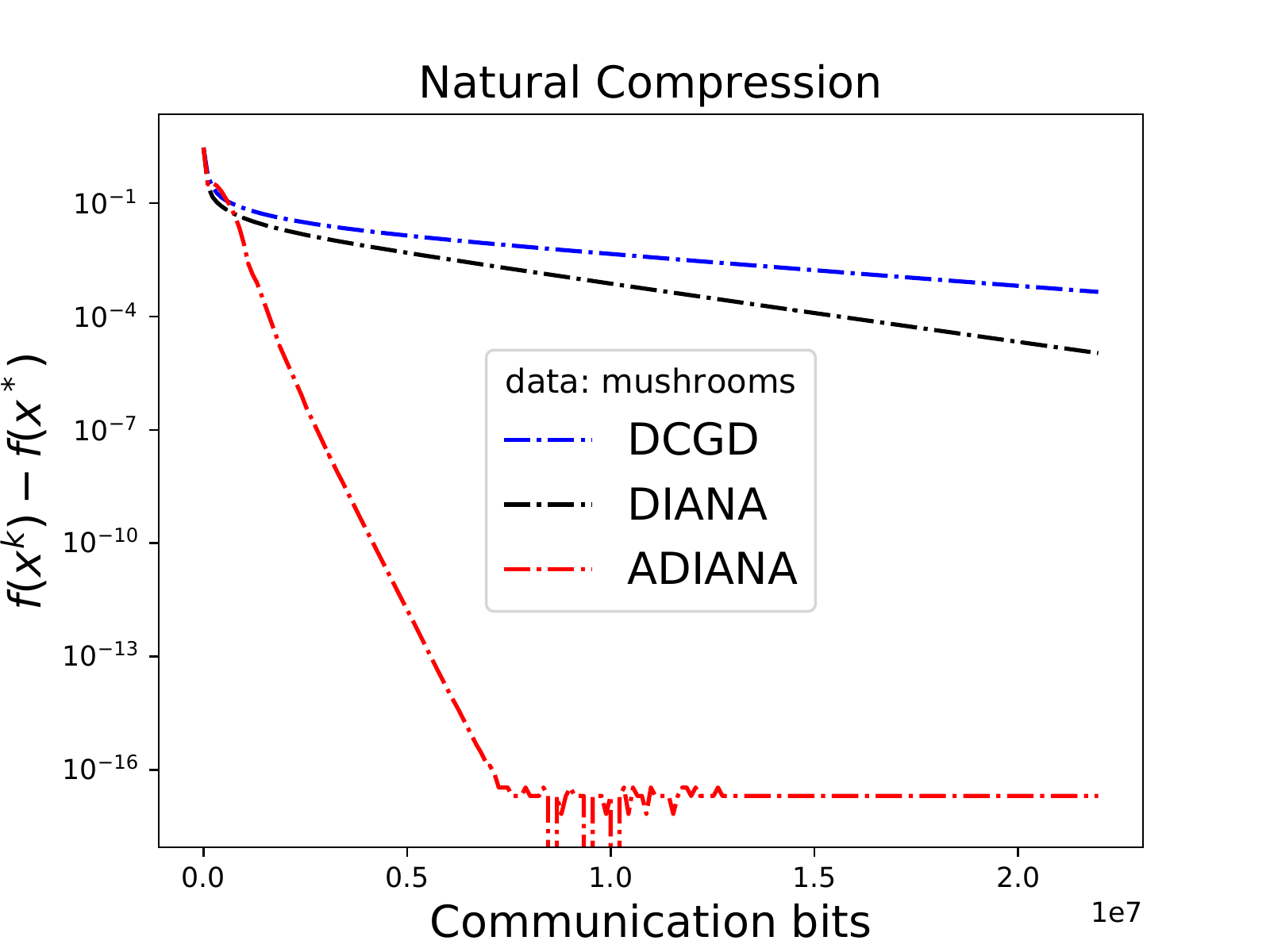}&
	\end{tabular}
	
	\caption{The communication complexity of different methods for three different compressors (random sparsification, random dithering and natural compression) on the  \texttt{mushrooms} dataset.}
	\label{fig:mushrooms111}
	
%
	\vspace{10mm}
	\centering
	\begin{tabular}{cccc}
		\includegraphics[width=0.32\linewidth]{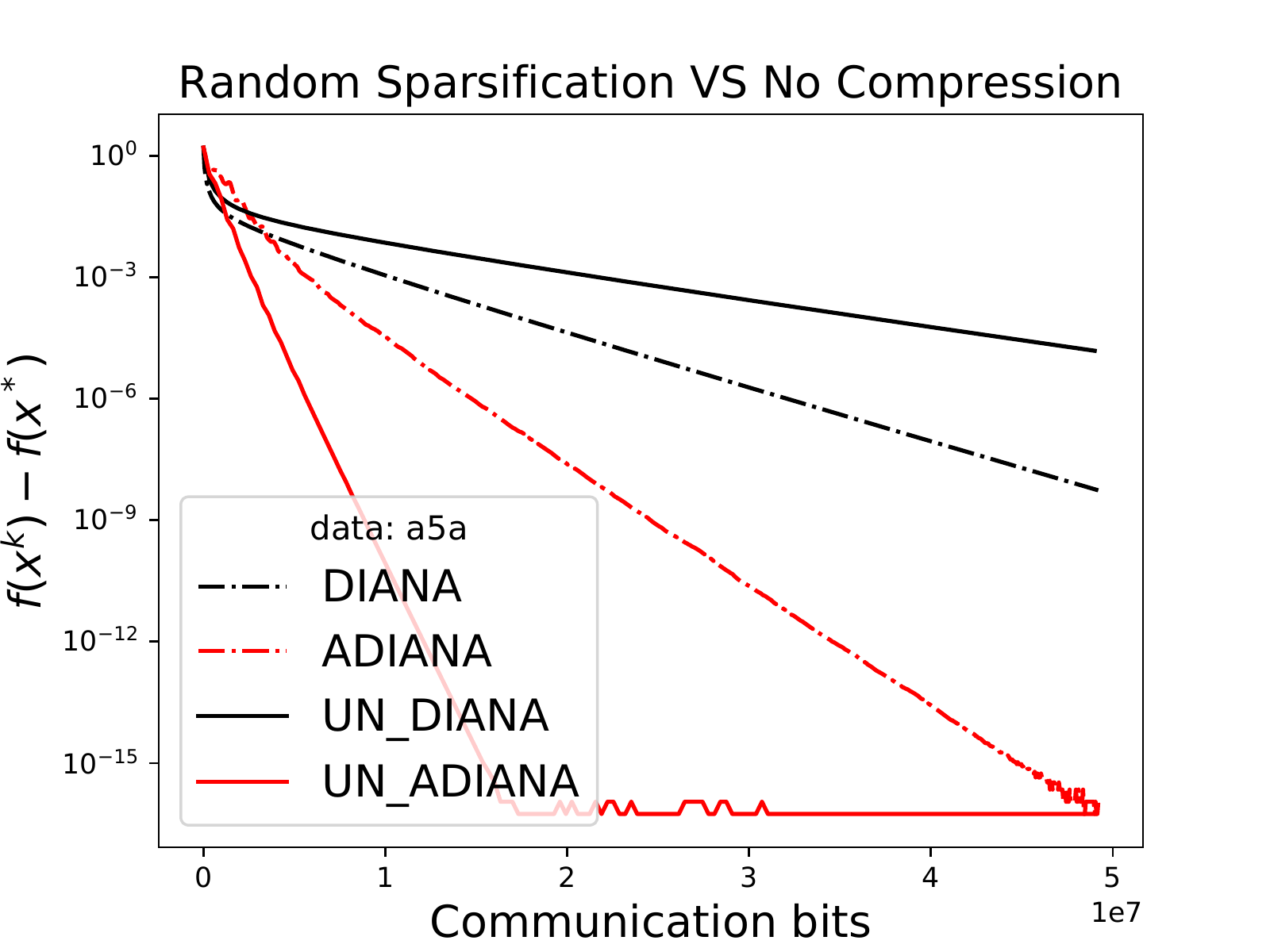}&
		\includegraphics[width=0.32\linewidth]{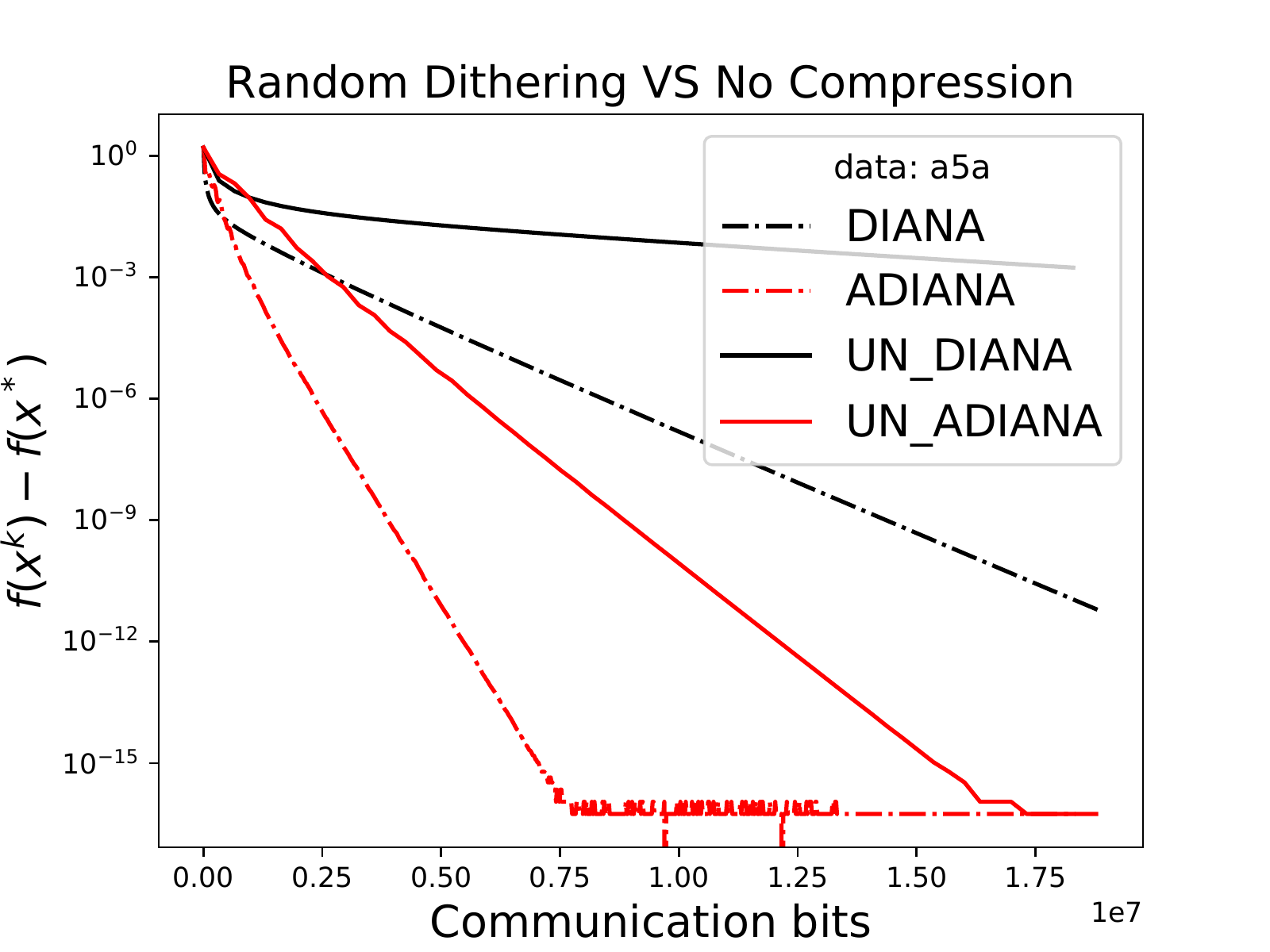}&
		\includegraphics[width=0.32\linewidth]{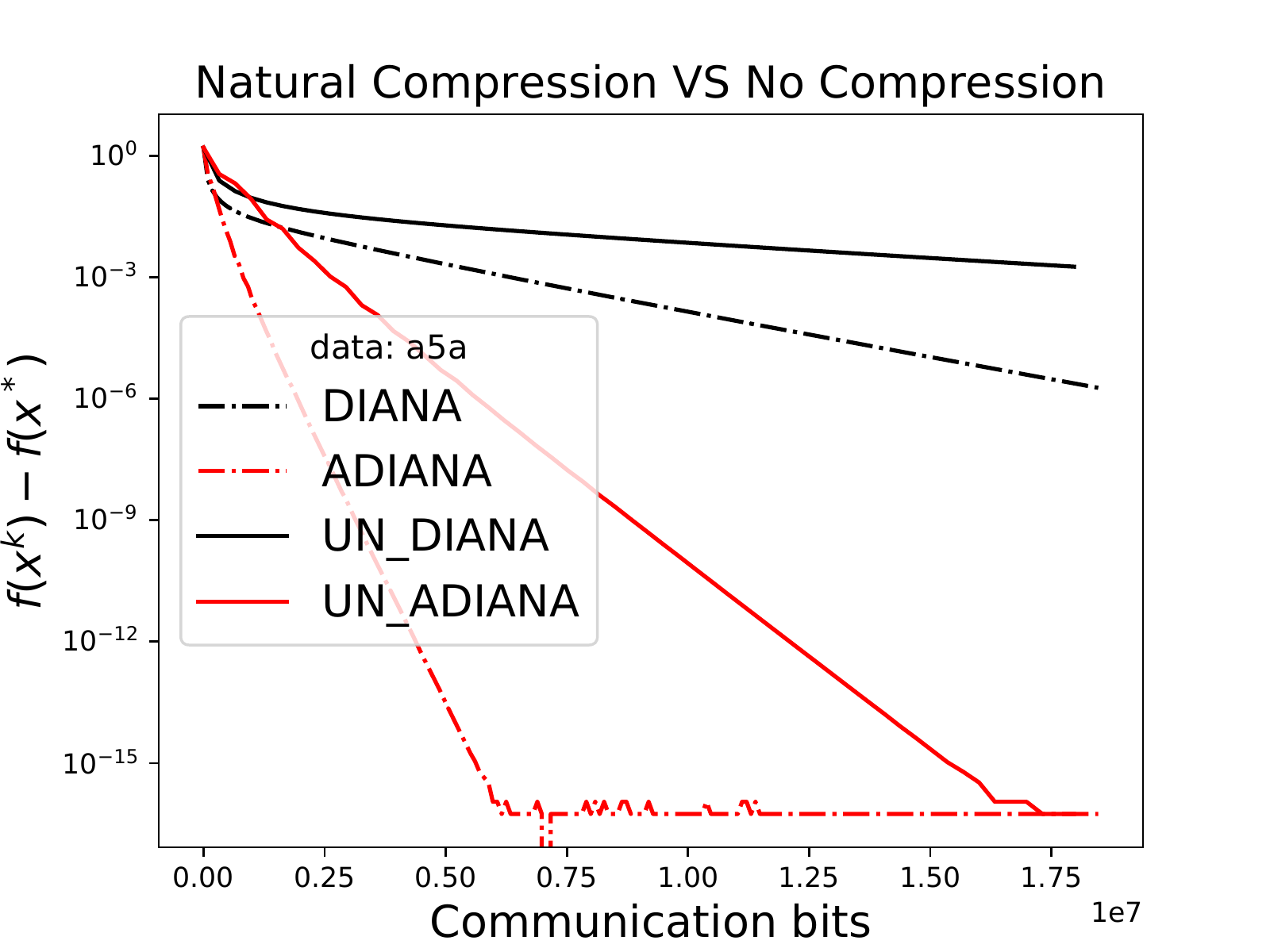}&
	\end{tabular}

	\caption{The communication complexity of DIANA and ADIANA with and without compression on the \texttt{a5a} dataset.}
	\label{fig:a5a222}
	
%
%
	\vspace{3mm}
	\centering
	\begin{tabular}{cccc}
		\includegraphics[width=0.32\linewidth]{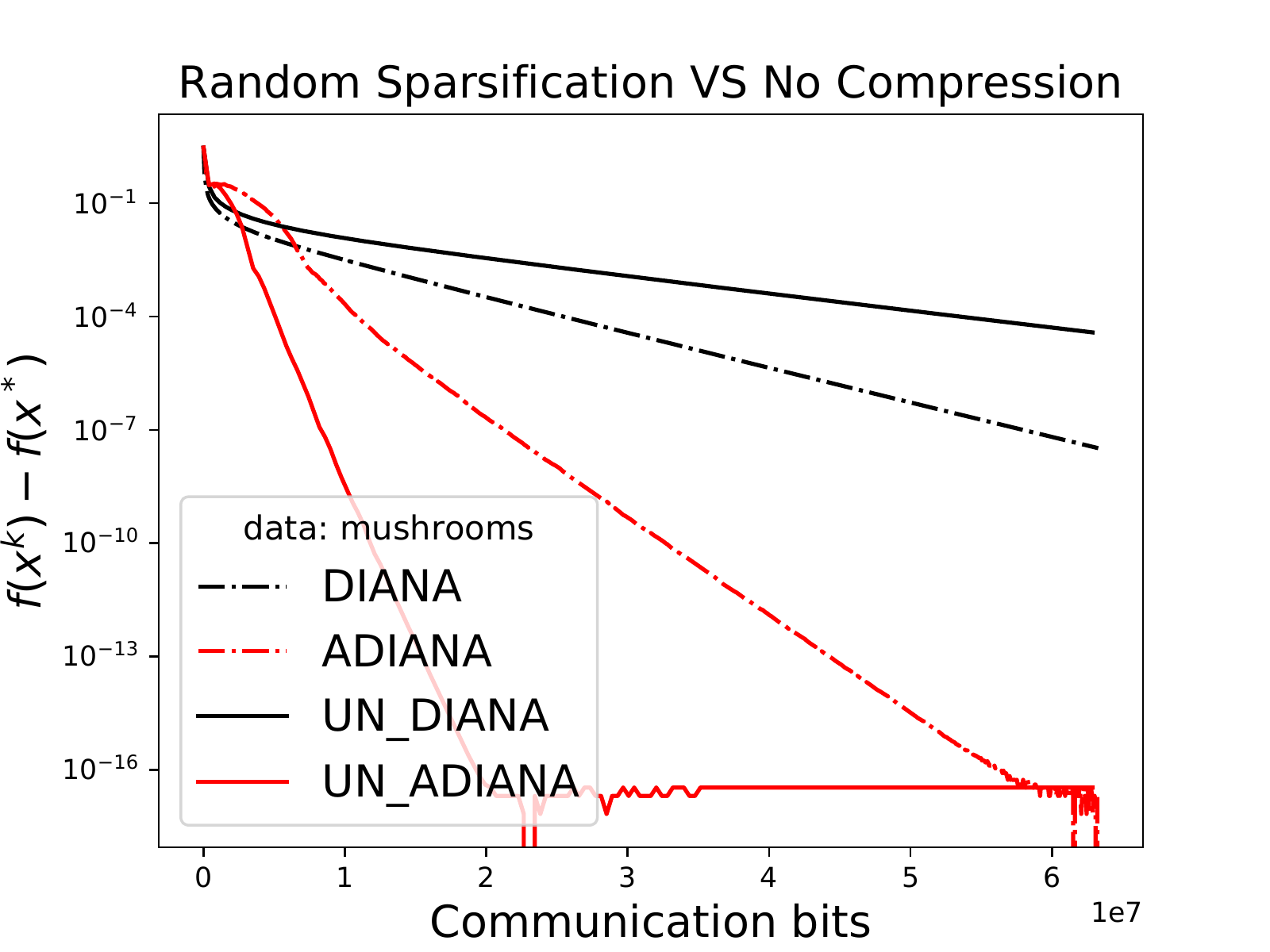}&
		\includegraphics[width=0.32\linewidth]{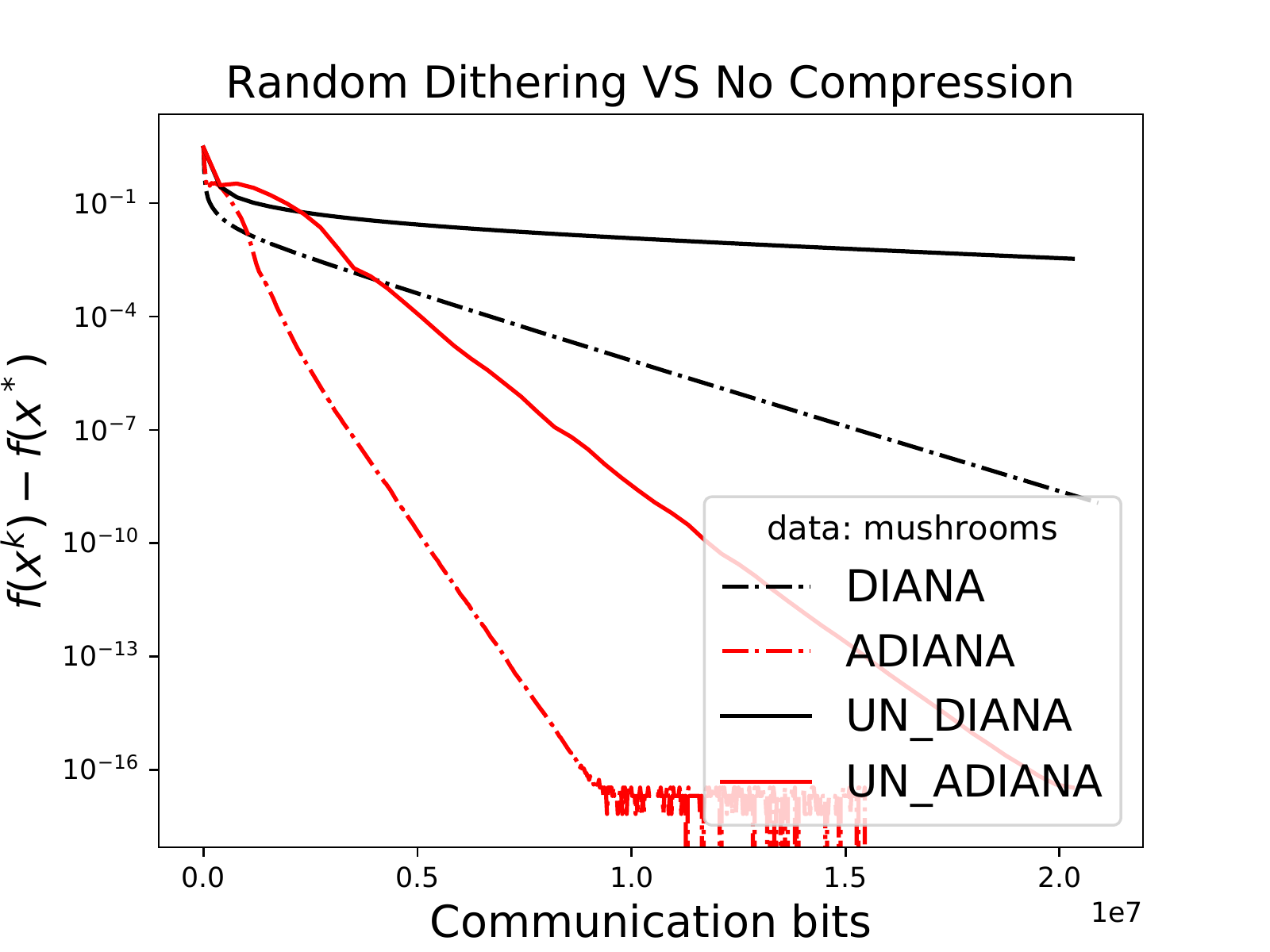}&
		\includegraphics[width=0.32\linewidth]{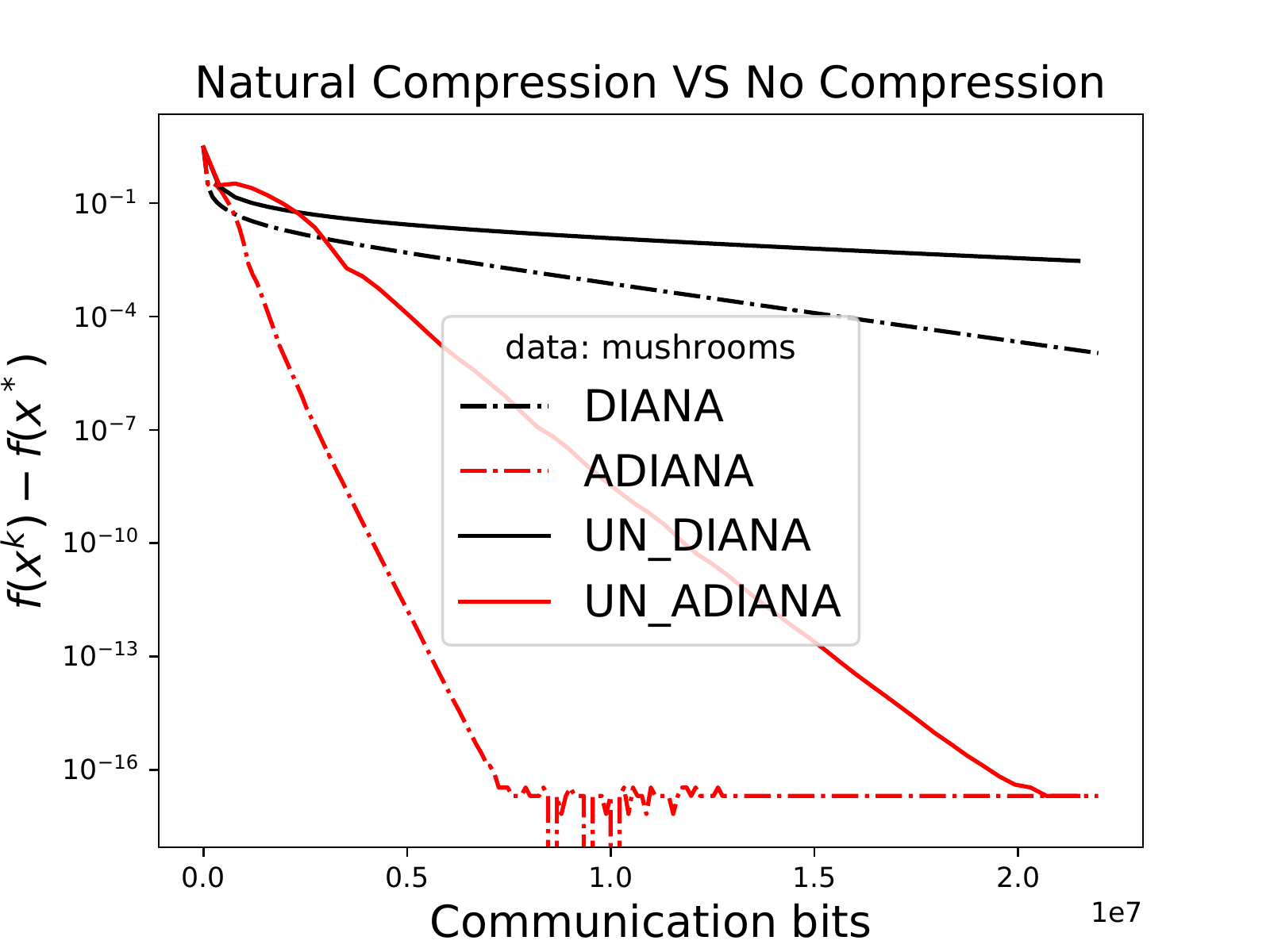}&
	\end{tabular}

	\caption{ The communication complexity of DIANA and ADIANA with and without compression on the \texttt{mushrooms} dataset.}
	\label{fig:mushrooms222}
	
\end{figure*}

Now, we compare our ADIANA and DIANA, with and without compression to show the communication efficiency of our accelerated method ADIANA in Figures~\ref{fig:a5a222} and \ref{fig:mushrooms222}. 

According to the left top and left bottom of Figure \ref{fig:mushrooms222}, DIANA is better than its uncompressed version if the compression operator is random sparsification. However, ADIANA behaves worse than its uncompressed version. 
For random dithering (middle figures) and natural compression (right figures), ADIANA is about twice faster than its uncompressed version, and is much faster than DIANA with/without compression.
These numerical results indicate that ADIANA (which enjoys both acceleration and compression) could be a more practical communication efficiency method, i.e., acceleration (better than non-accelerated DIANA) and compression (better than the uncompressed version), especially for random dithering and natural compression.

\clearpage
\section*{Acknowledgements}
The authors would like to acknowledge support from the KAUST Baseline Research Fund and the KAUST Visual Computing Center. 
Zhize Li and Xun Qian thank for support from the KAUST Extreme Computing Research Center.


\bibliography{acgd}
\bibliographystyle{icml2020}


\clearpage
\onecolumn
\appendix

{\bf \LARGE Appendix}

\section{Missing Proofs}
\label{sec:app1}

In this appendix, we provide the detailed proofs for our Theorems \ref{thm:cgdconvex}--\ref{thm:accdcgd}.

\subsection{Proof of Theorem \ref{thm:cgdconvex}}
\label{sec:appthm1pf}
We first restate our Theorem \ref{thm:cgdconvex} here and then provide the detailed proof.

\begingroup
\def\thetheorem{\ref{thm:cgdconvex}}
\begin{theorem}
	Suppose $f(x)$ is convex with $L$-Lipschitz continuous gradient and the compression operator $\calC(\cdot)$ satisfies \eqref{eq:comp}.  Let step size $\eta=\frac{1}{(1+\omega)L}$, then the number of iterations performed by CGD to find an $\epsilon$-solution such that $\EE[f(x^k)-f(x^*)]\leq \epsilon$ is at most $k=O(\frac{(1+\omega)L}{\epsilon})$.
\end{theorem}
\addtocounter{theorem}{-1}
\endgroup
\begin{proof}
	According to CGD update $x^{k+1} =x^k-\eta g^k$, we have
	\begin{align}
	\EB{\ns{x^{k+1}-\xs}} &= \EB{\ns{x^k-\eta g^k -\xs}} \notag\\
	&= \EB{\ns{x^k-\eta \cC(\nabla f(x^k)) -\xs}}  \notag\\
	&= \EB{\ns{x^k-\xs} -2\eta\inner{\cC(\nabla f(x^k))}{x^k-\xs} +\eta^2 \ns{\cC(\nabla f(x^k))}} \notag\\
	&\overset{\eqref{eq:comp}}{=}\ns{x^k-\xs} -2\eta\inner{\nabla f(x^k)}{x^k-\xs} +\eta^2 \EB{\ns{\cC(\nabla f(x^k))}} \notag\\
	&\overset{\eqref{eq:comp}}{=}\ns{x^k-\xs} -2\eta\inner{\nabla f(x^k)}{x^k-\xs} +\eta^2 \ns{\nabla f(x^k)}+\EB{\ns{\cC(\nabla f(x^k)) - \nabla f(x^k)}} \notag\\
	&\overset{\eqref{eq:comp}}{\leq}\ns{x^k-\xs} -2\eta\inner{\nabla f(x^k)}{x^k-\xs} +\eta^2 (1+\omega)\ns{\nabla f(x^k)} \notag \\
	&\leq \ns{x^k-\xs} -2\eta (f(x^k)-f(\xs)) +\eta^2 (1+\omega)\ns{\nabla f(x^k)}, \label{eq:relationx}
	\end{align}
	where the last inequality holds due to convexity of $f$.
	Besides, according to $L$-smoothness of $f$ (see \eqref{def_smoothness}), we have
	\begin{align}
	\EB{f(x^{k+1})-f(\xs)} 
	&\leq \EB{f(x^{k})-f(\xs) +\inner{\nabla f(x^k)}{x^{k+1}-x^k}  +\frac{L}{2}\ns{x^{k+1}-x^k}} \notag\\
	&= \EB{f(x^{k})-f(\xs) +\inner{\nabla f(x^k)}{-\eta \cC(\nabla f(x^k))}  +\frac{L\eta^2}{2}\ns{\cC(\nabla f(x^k))}} \notag\\
	&\overset{\eqref{eq:comp}}{=} \EB{f(x^{k})-f(\xs) -\eta \ns{\nabla f(x^k)} +\frac{L\eta^2}{2}\ns{\cC(\nabla f(x^k))}} \notag\\
	&\overset{\eqref{eq:comp}}{\leq} \EB{f(x^{k})-f(\xs) -\eta\left(1-\frac{L\eta(1+\omega)}{2}\right) \ns{\nabla f(x^k)}}. \label{eq:relationf}
	\end{align}
	Adding \eqref{eq:relationx}  and $\frac{\eta(1+\omega)}{1-\frac{L\eta(1+\omega)}{2}}$ times \eqref{eq:relationf} to cancel the term $\ns{\nabla f(x^k)}$, we have
	\begin{align}
	&\EB{\frac{\eta(1+\omega)}{1-\frac{L\eta(1+\omega)}{2}}\left(f(x^{k+1})-f(\xs)\right) + \ns{x^{k+1}-\xs} +2\eta \left(f(x^{k})-f(\xs)\right)}  \notag\\
	&\leq \EB{\frac{\eta(1+\omega)}{1-\frac{L\eta(1+\omega)}{2}}\left(f(x^{k})-f(\xs)\right) + \ns{x^{k}-\xs}}. \label{eq:relationfandx}
	\end{align}
	Summing up the above inequality from iteration $0$ to $k$, we have
	\begin{align}
	\EB{2\eta \sum_{i=0}^{k}\left(f(x^{i})-f(\xs)\right)}  
	&\leq \frac{\eta(1+\omega)}{1-\frac{L\eta(1+\omega)}{2}}\left(f(x^{0})-f(\xs)\right) + \ns{x^{0}-\xs} \notag\\
	&\leq \frac{\eta(1+\omega)}{1-\frac{L\eta(1+\omega)}{2}}\left(\frac{L}{2}\ns{x^0-\xs}\right) + \ns{x^{0}-\xs} \notag\\
	&=\frac{2\ns{x^{0}-\xs}}{2-L\eta(1+\omega)}, \label{eq:thm1last}
	\end{align}
	where the last inequality uses the $L$-smoothness of $f$. Finally, noting that $\EE[f(x^{i+1})-f(\xs)]\leq\EE[f(x^i)-f(\xs)]$ for all $i=0,\ldots k$ according to \eqref{eq:relationf}, then \eqref{eq:thm1last} turns to be 
	\begin{align}
	\EB{2\eta k\left(f(x^{k})-f(\xs)\right)}  
	&\leq \frac{2\ns{x^{0}-\xs}}{2-L\eta(1+\omega)} \notag \\
	\EB{f(x^{k})-f(\xs)}  
	&\leq \frac{\ns{x^{0}-\xs}}{\left(2-L\eta(1+\omega)\right)\eta k} =\frac{(1+\omega)L\ns{x^{0}-\xs}}{k},
	\end{align}
	where the last equality uses the choice of step size $\eta =\frac{1}{(1+\omega)L}$. Now the proof of Theorem \ref{thm:cgdconvex} is finished by letting the number of iteration $k=\frac{(1+\omega)L\ns{x^{0}-\xs}}{\epsilon}$, i.e.,  obtain the $\epsilon$-solution $x^k$ such that $\EE[f(x^k)-f(x^*)]\leq \epsilon$ within $O(\frac{(1+\omega)L}{\epsilon})$ iterations.
\end{proof}

\vspace{15pt}
\subsection{Proof of Theorems \ref{thm:acgdconvex} and \ref{thm:acgdstrongconvex}}
\label{sec:appthm23pf}
First, we restate our key Lemma \ref{lem:key} here. Then we use it to prove Theorems \ref{thm:acgdconvex} and \ref{thm:acgdstrongconvex}. Finally, we provide the proof for this key Lemma \ref{lem:key}.

\begingroup
\def\thelemma{\ref{lem:key}}
\begin{lemma}
	If parameters $\{\et\}, \{\at\}, \{\bt\}$, $\{\st\}$ and $p$ satisfy 
	$\at = \frac{1 - \st/p}{1-\bt\st/p}$
	$\bt \leq \min\{\frac{\mu\et}{\st p},1\}$, $p\geq \frac{(1+L\et)(1+\omega)}{2}$ and the compression operator $\calC(\cdot)$ satisfies \eqref{eq:comp},
	then we have for any iteration $k$, $\forall x \in \R^d$, 
	\begin{align}
	\frac{2\et}{\st^2}\EE[f(\xp) -f(x)]+ \EE[\ns{\vp-x}] 
	&\leq \left(1-\frac{\st}{p}\right)\frac{2\et}{\st^2}\left(f(\xt) -f(x)\right) +(1-\bt)\ns{\vt-x},  \label{eq:keyapp}
	\end{align}
	where the expectation is with respect to the randomness of compression $\calC(\cdot)$ at iteration $k$.
\end{lemma}
\addtocounter{lemma}{-1}
\endgroup

Now, we recall our Theorem \ref{thm:acgdconvex} here and are ready to prove it using Lemma \ref{lem:key}.
\begingroup
\def\thetheorem{\ref{thm:acgdconvex}}
\begin{theorem}
	Suppose $f(x)$ is convex with $L$-Lipschitz continuous gradient and the compression operator $\calC(\cdot)$ satisfies \eqref{eq:comp}.  Let parameters $\et \equiv \frac{1}{L}$, $\at =1-\frac{2}{k+2}$, $\bt\equiv 0$, $\st = \frac{2p}{k+2}$ and $p=1+\omega$, then the number of iterations performed by ACGD (Algorithm \ref{alg:acgd}) to find an $\epsilon$-solution such that $\EE[f(x^k)-f(x^*)]\leq \epsilon$ is at most $k=O\left((1+\omega)\sqrt{\frac{L}{\epsilon}}\right)$.
\end{theorem}
\addtocounter{theorem}{-1}
\endgroup
\begin{proofof}{Theorem \ref{thm:acgdconvex}}
	First, we know $\mu=0$ in this general convex case. By choosing  step size $\et \equiv \frac{1}{L}$, $p=1+\omega$, $\bt\equiv 0$, $\st = \frac{2p}{k+2}$ and $\at =1-\frac{2}{k+2}$, Lemma \ref{lem:key} turns to be
	\begin{align*}
	\frac{(k+2)^2}{2L(1+\omega)^2}\EE[f(\xp) -f(x)]+ \EE[\ns{\vp-x}] 
	&\leq \frac{(k+1)^2-1}{2L (1+\omega)^2}\left(f(\xt) -f(x)\right) + \ns{\vt-x}. 
	\end{align*}
	Summing up the above inequality from iteration $0$ to $k-1$ and letting $x=x^*$, we have
	\begin{align}
	\frac{(k+1)^2}{2L(1+\omega)^2}\EE[f(\xt) -f(x^*)]
	&\leq \sum_{i=0}^{k-1}\frac{-1}{2L (1+\omega)^2}\left(f(y^i) -f(x^*)\right) + \ns{y^0-x^*} \leq \ns{x^0-x^*},   \\
	\EE[f(\xt) -f(x^*)] &\leq \frac{2L(1+\omega)^2\ns{x^0-x^*}}{(k+1)^2}, \notag
	\end{align}
	where the second inequality uses $f(y^i)-f(x^*)\geq 0$ and $y^0=x^0$. 
	Now the proof of Theorem \ref{thm:acgdconvex} is finished, i.e.,  we obtain the $\epsilon$-solution $y^k$ such that $\EE[f(y^k)-f(x^*)]\leq \epsilon$ within $k=O\left((1+\omega)\sqrt{\frac{L}{\epsilon}}\right)$ iterations.
\end{proofof}

\newpage
Similarly, we recall our Theorem \ref{thm:acgdstrongconvex} here and also prove it using Lemma \ref{lem:key}.
\begingroup
\def\thetheorem{\ref{thm:acgdstrongconvex}}
\begin{theorem}
	Suppose $f(x)$ is $\mu$-strongly convex with $L$-Lipschitz continuous gradient and the compression operator $\calC(\cdot)$ satisfies \eqref{eq:comp}. Let parameters $\et \equiv \frac{1}{L}$, $\at \equiv \frac{p}{p+\sqrt{\mu/L}}$, $\bt \equiv \frac{\sqrt{\mu/L}}{p}$, $\st \equiv \sqrt{\frac{\mu}{L}}$ and $p=1+\omega$, then the number of iterations performed by ACGD (Algorithm \ref{alg:acgd}) to find an $\epsilon$-solution such that $\EE[f(x^k)-f(x^*)]\leq \epsilon$ (or $\EE[\ns{x^k-x^*}]\leq \epsilon$)  is at most $k=O\left((1+\omega)\sqrt{\frac{L}{\mu}}\log \frac{1}{\epsilon}\right)$.
\end{theorem}
\addtocounter{theorem}{-1}
\endgroup
\begin{proofof}{Theorem \ref{thm:acgdstrongconvex}}
	By choosing  step size $\et \equiv \frac{1}{L}$, $p=1+\omega$, $\st \equiv \sqrt{\frac{\mu}{L}}$, $\bt \equiv \frac{\sqrt{\mu/L}}{p}$ and $\at \equiv \frac{1}{1+\sqrt{\mu/p^2L}}$, Lemma \ref{lem:key} turns to be
	\begin{align*}
	\frac{2}{\mu}\EE[f(\xp) -f(x)] + \EE[\ns{\vp-x}] 
	&\leq \left(1-\frac{\sqrt{\mu/L}}{1+\omega}\right)\frac{2}{\mu}\left(f(\xt) -f(x)\right) 
	+\left(1-\frac{\sqrt{\mu/L}}{1+\omega}\right)\ns{\vt-x} \\
	&= \left(1-\frac{\sqrt{\mu/L}}{1+\omega}\right)\left(\frac{2}{\mu}\left(f(\xt) -f(x)\right) + \ns{\vt-x}\right).
	\end{align*}
	Telescoping the above inequality from iteration $0$ to $k$ and letting $x=x^*$, we have
	\begin{align}
	\frac{2}{\mu}\EE[f(\xt) -f(x^*)] + \EE[\ns{\vt-x^*}] 
	&= \left(1-\frac{\sqrt{\mu/L}}{1+\omega}\right)^k\left(\frac{2}{\mu}\left(f(y^0) -f(x^*))\right) + \ns{z^0-x^*}\right) \notag\\
	&\leq \left(1-\frac{\sqrt{\mu/L}}{1+\omega}\right)^k\left(\frac{4}{\mu}\left(f(x^0) -f(x^*)\right)\right),  \label{eq:appmu} \\
	\EE[f(\xt) -f(x^*)] &\leq \left(1-\frac{\sqrt{\mu/L}}{1+\omega}\right)^k\left(f(x^0) -f(x^*)\right), \notag\\
	\EE[\ns{\vt-x^*}] &\leq \left(1-\frac{\sqrt{\mu/L}}{1+\omega}\right)^k\left(\frac{4}{\mu}\left(f(x^0) -f(x^*)\right)\right), \notag
	\end{align}
	where the first inequality uses $\mu$-strongly convex of $f$ (see \eqref{def_strongconvex}) and $y^0=z^0=x^0$. 
	Now the proof of Theorem \ref{thm:acgdstrongconvex} is finished, i.e.,  we obtain the $\epsilon$-solution $y^k$ (or $z^k$) such that $\EE[f(y^k)-f(x^*)]\leq \epsilon$ (or $\EE[\ns{\vt-x^*}]\leq \epsilon$) within $k=O\left((1+\omega)\sqrt{\frac{L}{\mu}}\log \frac{1}{\epsilon}\right)$ iterations.
\end{proofof}

\vspace{15pt}
\head{Key Lemma} Now, the only remaining thing is to prove the key Lemma \ref{lem:key}.

\begin{proofof}{Lemma \ref{lem:key}}
	First, we get the following equality (its proof is deferred to the end):
	\begin{align}
	\forall x\in \R^d,\quad  \EE[\ns{\vp-x}]
	&=  (1-\bt)\ns{\vt-x} -\bt(1-\bt)\ns{\yt-\vt} +\bt\ns{\yt-x} +\frac{\et^2}{\st^2 p^2}\EE\ns{\gt}  \notag \\
	&\qquad\qquad	 
	+\frac{2\et}{\st^2}\innerb{\nabla f(\yt)}{\left(1-\frac{\st}{p} \right)(\xt-\yt)+\frac{\st}{p}(x-\yt)}. \label{eq:lemz1}
	\end{align}
	Then, to cancel the last inner product in \eqref{eq:lemz1}, we use the  property (smoothness and/or strong convexity) of $f$:
	\begin{align}
	\EE[f(\xp)] &\leq \EB{f(\yt) +\inner{\nabla f(\yt)}{\xp-\yt} + \frac{L}{2}\ns{\xp-\yt}} \label{eq:lsmooth} \\
	&= \EB{f(\yt) -\frac{\et}{p}\inner{\nabla f(\yt)}{\gt} + \frac{L\et^2}{2p^2}\ns{\gt}} \notag\\
	&= f(\yt) -\frac{\et}{p} \ns{\nabla f(\yt)} + \frac{L\et^2}{2p^2}\EE\ns{\gt} \notag\\
	&\leq f(u) -\inner{\nabla f(\yt)}{u-\yt}-\frac{\mu}{2}\ns{u-\yt} -\frac{\et}{p} \ns{\nabla f(\yt)} + \frac{L\et^2}{2p^2}\EE\ns{\gt}, \quad \forall u\in \R^d   \label{eq:mustrong}
	\end{align} 
	where \eqref{eq:lsmooth} and \eqref{eq:mustrong} uses $L$-smoothness (see \eqref{def_smoothness}) and $\mu$-strong convexity (see \eqref{def_strongconvex}), respectively.
	By adding $\left(1-\frac{\st}{p}\right)$ times \eqref{eq:mustrong} (where $u=\xt$) and $\frac{\st}{p}$ times \eqref{eq:mustrong} (where $u=x$), we have 
	\begin{align}
	\EE[f(\xp)] &\leq \left(1-\frac{\st}{p}\right) f(\xt) + \frac{\st}{p} f(x) 
	-\innerb{\nabla f(\yt)}{\left(1-\frac{\st}{p}\right)(\xt-\yt) +\frac{\st}{p}(x-\yt)}  \notag \\
	&	\qquad \qquad - \left(1-\frac{\st}{p}\right)\frac{\mu}{2}\ns{\xt-\yt} -\frac{\st \mu}{2p}\ns{x-\yt} -\frac{\et}{p} \ns{\nabla f(\yt)} + \frac{L\et^2}{2p^2}\EE\ns{\gt}. \label{eq:getinner}
	\end{align}
	Now, this key lemma (i.e., \eqref{eq:keyapp}) is proved as follows by adding $\frac{2\et}{\st^2}$ times \eqref{eq:getinner} and \eqref{eq:lemz1}:
	\begin{align}
	&\frac{2\et}{\st^2}\EE[f(\xp) -f(x)] + \EE[\ns{\vp-x}]  \notag\\
	&\leq \left(1-\frac{\st}{p}\right)\frac{2\et}{\st^2}\left(f(\xt) -f(x)\right) +(1-\bt)\ns{\vt-x}
	-\bt(1-\bt)\ns{\yt-\vt} -\left(\frac{\mu\et}{\st p}-\bt \right)\ns{\yt-x} \notag\\
	&\qquad \qquad -\frac{2\et^2}{\st^2 p}\ns{\nabla f(\yt)} +\frac{L\et^3+\et^2}{\st^2 p^2}\EE\ns{g_k} \notag\\
	&\leq \left(1-\frac{\st}{p}\right)\frac{2\et}{\st^2}\left(f(\xt) -f(x)\right) +(1-\bt)\ns{\vt-x} 
	-\frac{2\et^2}{\st^2 p}\ns{\nabla f(\yt)} +\frac{L\et^3+\et^2}{\st^2 p^2}\EE\ns{g_k} \label{eq:betacond}\\
	&= \left(1-\frac{\st}{p}\right)\frac{2\et}{\st^2}\left(f(\xt) -f(x)\right) +(1-\bt)\ns{\vt-x}
	-\frac{2\et^2}{\st^2 p}\ns{\nabla f(\yt)} +\frac{L\et^3+\et^2}{\st^2 p^2}\EE\ns{\calC(\nabla f(\yt)} \notag \\
	&= \left(1-\frac{\st}{p}\right)\frac{2\et}{\st^2}\left(f(\xt) -f(x)\right) +(1-\bt)\ns{\vt-x}
	-\frac{2\et^2p-L\et^3-\et^2}{\st^2 p^2}\ns{\nabla f(\yt)} \notag\\
	&\qquad \qquad +\frac{L\et^3+\et^2}{\st^2 p^2}\EE[\ns{\calC(\nabla f(\yt))-\nabla f(\yt)}]  \label{eq:exp-var}\\
	&\leq \left(1-\frac{\st}{p}\right)\frac{2\et}{\st^2}\left(f(\xt) -f(x)\right) +(1-\bt)\ns{\vt-x}
	-\frac{2\et^2p-(1+\omega)(L\et^3+\et^2)}{\st^2 p^2}\ns{\nabla f(\yt)} \label{eq:usecomp}\\
	&\leq \left(1-\frac{\st}{p}\right)\frac{2\et}{\st^2}\left(f(\xt) -f(x)\right) +(1-\bt)\ns{\vt-x}, 
	\label{eq:pcond}
	\end{align}
	where \eqref{eq:betacond} holds due to condition $\bt \leq \min\{\frac{\mu\et}{\st p},1\}$, \eqref{eq:exp-var} and \eqref{eq:usecomp} use the property of compression \eqref{eq:comp}, and \eqref{eq:pcond} uses the condition $p\geq \frac{(1+L\et)(1+\omega)}{2}$.
	Now, the only remaining thing is to prove \eqref{eq:lemz1}. For any $ x\in\R^d$, we have
	\begin{align}
	\EE[\ns{\vp-x}] 
	&= \E{\nsB{\frac{1}{\st}\xp + \left(\frac{1}{p}-\frac{1}{\st} \right)\xt 
			+ \left(1-\frac{1}{p}\right)(1-\bt)\vt
			+\left(1-\frac{1}{p}\right)\bt \yt -x}} \notag\\
	&= \E{\nsB{\left(\frac{1}{\st}+\left(1-\frac{1}{p}\right)\bt\right) \yt+ \left(\frac{1}{p}-\frac{1}{\st}\right)\xt 
			+ \left(1-\frac{1}{p}\right)(1-\bt)\vt -x
			- \frac{\et}{\st p} \gt}} \notag\\
	&=\nsB{\left(\frac{1}{\st}+\left(1-\frac{1}{p}\right)\bt\right) \yt+ \left(\frac{1}{p}-\frac{1}{\st}\right)\xt 
		+ \left(1-\frac{1}{p}\right)(1-\bt)\vt -x} +\frac{\et^2}{\st^2 p^2}\EE\ns{\gt} \notag \\
	&\qquad\qquad  
	+\frac{2\et}{\st p}\innerB{\EE[\gt]}{x-\left(\frac{1}{\st}+\left(1-\frac{1}{p}\right)\bt\right) \yt- \left(\frac{1}{p}-\frac{1}{\st}\right)\xt 
		- \left(1-\frac{1}{p}\right)(1-\bt)\vt} \notag \\
	&=(1-\bt)\ns{\vt-x} -\bt(1-\bt)\ns{\yt-\vt} +\bt\ns{\yt-x} +\frac{\et^2}{\st^2 p^2}\EE\ns{\gt}  \label{eq:lemz11} \\
	&\qquad\qquad	 
	+\frac{2\et}{\st p}\innerB{\EE[\gt]}{x-\left(\frac{1}{\st}+\left(1-\frac{1}{p}\right)\bt\right) \yt- \left(\frac{1}{p}-\frac{1}{\st}\right)\xt 
		- \left(1-\frac{1}{p}\right)(1-\bt)\vt} \notag \\
	&=(1-\bt)\ns{\vt-x} -\bt(1-\bt)\ns{\yt-\vt} +\bt\ns{\yt-x} +\frac{\et^2}{\st^2 p^2}\EE\ns{\gt}  \notag \\
	&\qquad\qquad	 
	+\frac{2\et}{\st^2}\innerb{\nabla f(\yt)}{\left(1-\frac{\st}{p}\right)(\xt-\yt)+\frac{\st}{p}(x-\yt)}, \label{eq:lemz12} 
	\end{align}
	where \eqref{eq:lemz11} and \eqref{eq:lemz12} use equalities \eqref{eq:lemz11eq} and \eqref{eq:lemz12eq}, respectively. Further,
	\begin{align}
	& \nsB{\left(\frac{1}{\st}+\left(1-\frac{1}{p}\right)\bt\right) \yt+ \left(\frac{1}{p}-\frac{1}{\st}\right)\xt 
		+ \left(1-\frac{1}{p}\right)(1-\bt)\vt -x} \notag\\
	&=\nsB{\left(\frac{1}{\st}-\left(1-\frac{1}{p}\right)\bt\right) \yt
		- \frac{1}{\st}\left(1-\frac{\st}{p}\right)\xt 
		+ \left(1-\frac{1}{p}\right)(1-\bt)\vt -x} \notag\\
	&=\nsB{\left(\frac{1}{\st}-\left(1-\frac{1}{p}\right)\bt\right) \yt
		- \frac{1}{\st}\left(1-\frac{\bt\st}{p}\right)\at\xt 
		+ \left(1-\frac{1}{p}\right)(1-\bt)\vt -x} \label{eq:useat1} \\
	&=\nsB{\left(\frac{1}{\st}-\left(1-\frac{1}{p}\right)\bt\right) \yt
		- \frac{1}{\st}\left(1-\frac{\bt\st}{p}\right)(\yt-(1-\at)\vt) 
		+ \left(1-\frac{1}{p}\right)(1-\bt)\vt -x} \label{eq:usexyz} \\
	&=\nsB{\bt \yt
		+ \frac{1}{\st}\left(1-\frac{\bt\st}{p}\right)(1-\at)\vt
		+ \left(1-\frac{1}{p}\right)(1-\bt)\vt -x} \notag \\
	&=\nsB{\bt \yt
		+ \frac{1}{p}(1-\bt)\vt
		+ \left(1-\frac{1}{p}\right)(1-\bt)\vt -x} \label{eq:useat2} \\
	&=\ns{\vt-x+\bt(\yt-\vt)} \notag\\
	&=\ns{\vt-x}+\bt^2\ns{\yt-\vt}+2\bt\inner{\yt-\vt}{\vt-x} \notag\\
	&=(1-\bt)\ns{\vt-x} -\bt(1-\bt)\ns{\yt-\vt} +\bt\ns{\yt-x}, \label{eq:lemz11eq} 
	\end{align}
	where \eqref{eq:useat1} and  \eqref{eq:useat2} hold due to the condition $\at = \frac{1 - \st/p}{1-\bt\st/p}$, and \eqref{eq:usexyz} holds due to the relation $\yt=\at \xt +(1-\at)\vt$ (Line \ref{line:y} in Algorithm \ref{alg:acgd}). 
	Now, we finish the proof for the inner product term.
	\begin{align}
	&\frac{2\et}{\st p}\innerB{\nabla f(\yt)}{x-\left(\frac{1}{\st}+\left(1-\frac{1}{p}\right)\bt\right) \yt- \left(\frac{1}{p}-\frac{1}{\st}\right)\xt 
		- \left(1-\frac{1}{p}\right)(1-\bt)\vt}  \notag\\
	&=\frac{2\et}{\st^2}\innerB{\nabla f(\yt)}{\frac{\st}{p}x-\left(\frac{1}{p}+\left(1-\frac{1}{p}\right)\frac{\bt\st}{p}\right) \yt+ \frac{1}{p}\left(1-\frac{\st}{p}\right)\xt 
		- \left(1-\frac{1}{p}\right)(1-\bt)\frac{\st}{p}\vt}  \notag\\
	&=\frac{2\et}{\st^2}\innerB{\nabla f(\yt)}{\frac{\st}{p}x-\left(\frac{1}{p}+\left(1-\frac{1}{p}\right)\frac{\bt\st}{p}\right) \yt+ \frac{1}{p}\left(1-\frac{\st}{p}\right)\xt 
		- \left(1-\frac{1}{p}\right)\left(1-\frac{\bt\st}{p}\right)(1-\at)\vt}   \label{eq:useat3}\\
	&=\frac{2\et}{\st^2}\innerB{\nabla f(\yt)}{\frac{\st}{p}x-\left(\frac{1}{p}+\left(1-\frac{1}{p}\right)\frac{\bt\st}{p}\right) \yt+ \frac{1}{p}\left(1-\frac{\st}{p}\right)\xt 
		- \left(1-\frac{1}{p}\right)\left(1-\frac{\bt\st}{p}\right)(\yt-\at\xt)}  \label{eq:usexyz1} \\
	&=\frac{2\et}{\st^2}\innerB{\nabla f(\yt)}{\frac{\st}{p}x - \yt+ \frac{1}{p}\left(1-\frac{\st}{p}\right)\xt 
		+ \left(1-\frac{1}{p}\right)\left(1-\frac{\bt\st}{p}\right)\at\xt}  \notag\\
	&=\frac{2\et}{\st^2}\innerB{\nabla f(\yt)}{\frac{\st}{p}x - \yt+ \frac{1}{p}\left(1-\frac{\st}{p}\right)\xt 
		+ \left(1-\frac{1}{p}\right)\left(1-\frac{\st}{p}\right)\xt}  \label{eq:useat4}\\
	&=\frac{2\et}{\st^2}\innerb{\nabla f(\yt)}{\left(1-\frac{\st}{p}\right)(\xt-\yt)+\frac{\st}{p}(x-\yt)}, \label{eq:lemz12eq} 
	\end{align}
	where \eqref{eq:useat3} and  \eqref{eq:useat4} hold due to the condition $\at = \frac{1 - \st/p}{1-\bt\st/p}$, and \eqref{eq:usexyz1} holds due to the relation $\yt=\at \xt +(1-\at)\vt$ (Line \ref{line:y} in Algorithm \ref{alg:acgd}). 
\end{proofof}

\newpage
\subsection{Proof of Theorem \ref{thm:accdcgd}}
\label{sec:appthm4pf}
In this section, we provide the detailed proof for accelerated result in the distributed case.
Similar to previous Section \ref{sec:appthm23pf}, we first restate our Lemmas \ref{lem:keyaccdcgd}--\ref{lem:keyh} here. Then we use them to prove Theorem \ref{thm:accdcgd}. Finally, we provide the proof for these Lemmas.
Before restating lemmas, we recall the following notation:
\begin{align}
\cZ^k &\eqdef \norm{z^k - x^*}^2, \label{def:zapp}\\
\cY^k &\eqdef P(y^k) - P(x^*), \label{def:yapp}\\
\cW^k &\eqdef P(w^k) - P(x^*), \label{def:wapp}\\ 
\cH^k &\eqdef \frac{1}{n}\sum_{i=1}^n \norm{h_i^k - \nabla f_i(w^k)}^2. \label{def:happ}
\end{align}

\begingroup
\def\thelemma{\ref{lem:keyaccdcgd}}
\begin{lemma}
	If parameters satisfy 
	$\eta \leq \frac{1}{2L}, \theta_1 \leq \frac{1}{4}, \theta_2 =\frac{1}{2},  \gamma=\frac{\eta}{2(\theta_1 + \eta\mu)}$  and $\beta = 1-\gamma \mu$, 
	then we have for any iteration $k$, 
	\begin{align}
	\frac{2\gamma\beta}{\theta_1}\E{\cY^{k+1}} +  \E{\cZ^{k+1}} 
	&\leq
	(1 - \theta_1 - \theta_2)\frac{2\gamma\beta}{\theta_1}\cY^k
	+
	\beta\cZ^k  +
	2\gamma\beta\frac{\theta_2}{\theta_1}\cW^k
	+
	\frac{\gamma\eta}{\theta_1}\E{\norm{g^k - \nabla f(x^k)}^2} \notag\\ 
	&\qquad -
	\frac{\gamma}{4Ln\theta_1}\sum_{i=1}^n\norm{\nabla f_i(w^k) - \nabla f_i(x^k)}^2 
	- \frac{\gamma}{8Ln\theta_1}\sum_{i=1}^n\norm{\nabla f_i(y^k) - \nabla f_i(x^k)}^2. \label{eq:lem2app}  
	\end{align}
\end{lemma}
\addtocounter{lemma}{-1}
\endgroup

\begingroup
\def\thelemma{\ref{lem:keyw}}
\begin{lemma}
	According to Line \ref{line:prob} of Algorithm \ref{alg:ADIANA} and Definition \eqref{def:yapp}--\eqref{def:wapp}, we have
	\begin{equation}
	\E{\cW^{k+1}} = (1-p)\cW^k + p\cY^k. \label{eq:wapp}
	\end{equation}
\end{lemma}
\addtocounter{lemma}{-1}
\endgroup

\begingroup
\def\thelemma{\ref{lem:keyg}}
\begin{lemma}
	If the compression operator $\calC(\cdot)$ satisfies \eqref{eq:comp}, we have
	\begin{align}
	\E{\norm{g^k - \nabla f(x^k)}^2}  \leq
	\frac{2\omega}{n^2}\sum_{i=1}^n \norm{\nabla f_i(w^k) - \nabla f_i(x^k)}^2
	+
	\frac{2\omega}{n}\cH^k \label{eq:gapp}
	\end{align}
\end{lemma}
\addtocounter{lemma}{-1}
\endgroup

\begingroup
\def\thelemma{\ref{lem:keyh}}
\begin{lemma}
	If $\alpha \leq 1/(1+\omega)$, we have
	\begin{align}
	\E{\cH^{k+1}} &\leq (1-\frac{\alpha}{2})\cH^k   + (1 + \frac{2p}{\alpha})\frac{2p}{n}\left(\sum_{i=1}^n\norm{\nabla f_i(w^k) - \nabla f_i(x^k)}^2  
	+ \sum_{i=1}^n\norm{\nabla f_i(y^k) - \nabla f_i(x^k)}^2\right). \label{eq:happ}
	\end{align}
\end{lemma}
\addtocounter{lemma}{-1}
\endgroup

Now, we recall our Theorem \ref{thm:accdcgd} here and are ready to prove it using Lemmas \ref{lem:keyaccdcgd}--\ref{lem:keyh}.
\begingroup
\def\thetheorem{\ref{thm:accdcgd}}
\begin{theorem}
	Suppose $f(x)$ is $\mu$-strongly convex and all $f_is$ have $L$-Lipschitz continuous gradients, and the compression operator $\calC(\cdot)$ satisfies \eqref{eq:comp}. Let parameters $\alpha = \frac{1}{\omega+1}$, $\eta =\min\{\frac{1}{2L}, \frac{n}{64\omega(2p(\omega+1)+1)^2L}\}$, $\theta_1 = \min\left\{\frac{1}{4}, \sqrt{\frac{\eta\mu}{p}} \right\}$, $\theta_2 = \frac{1}{2}$, 
	$\gamma=\frac{\eta}{2(\theta_1 + \eta\mu)}$, $\beta = 1-\gamma \mu$, and $p=\min\big\{1, \frac{\max\{1,\sqrt{n/32\omega}-1\}}{2(1+\omega)}\big\}$, then the number of iterations performed by ADIANA (Algorithm \ref{alg:ADIANA}) to find an $\epsilon$-solution such that $\EE[\ns{z^k-x^*}]\leq \epsilon$  is at most
	\begin{equation*}
	k=
	\begin{cases}
	O\left(\left[\omega + \omega\sqrt{\frac{L}{n \mu}}\ \right]\log \frac{1}{\epsilon}  \right),& n \leq \omega,\\
	O\left(\left[\omega  + \sqrt{\frac{L}{\mu}} + \sqrt{\sqrt{\frac{\omega}{n}} \frac{\omega L}{\mu}}\ \right]\log \frac{1}{\epsilon}  \right),&  n>\omega. 
	\end{cases}
	\end{equation*}
\end{theorem}
\addtocounter{theorem}{-1}
\endgroup

\begin{proofof}{Theorem \ref{thm:accdcgd}}
	We define the following Lyapunov function $\Psi$ and induce it as follows:
	\begin{align}
	\E{\Psi^{k+1}}
	&\eqdef
	\E{\cZ^{k+1} + \frac{2\gamma\beta}{\theta_1}\cY^{k+1} + 2\gamma\beta \frac{\theta_2(1+\theta_1)}{\theta_1p}\cW^{k+1}
		+\frac{8\gamma\eta\omega}{\alpha\theta_1 n}\cH^{k+1}}    \notag \\
	&\leq
	\beta\cZ^k  +
	(1 - \theta_1 - \theta_2)\frac{2\gamma\beta}{\theta_1}\cY^k
	+
	2\gamma\beta\frac{\theta_2}{\theta_1}\cW^k
	+
	\frac{\gamma\eta}{\theta_1}\E{\norm{g^k - \nabla f(x^k)}^2} \notag\\ 
	&\qquad -
	\frac{\gamma}{4Ln\theta_1}\sum_{i=1}^n\norm{\nabla f_i(w^k) - \nabla f_i(x^k)}^2 
	- \frac{\gamma}{8Ln\theta_1}\sum_{i=1}^n\norm{\nabla f_i(y^k) - \nabla f_i(x^k)}^2 \notag\\
	&\qquad  
	+ \E{2\gamma\beta \frac{\theta_2(1+\theta_1)}{\theta_1p}\cW^{k+1}
		+\frac{8\gamma\eta\omega}{\alpha\theta_1 n}\cH^{k+1}} \label{eq:uselm2}\\
	&=
	\beta\cZ^k  +
	(1 - \theta_1 - \theta_2)\frac{2\gamma\beta}{\theta_1}\cY^k
	+
	2\gamma\beta\frac{\theta_2}{\theta_1}\cW^k
	+
	\frac{\gamma\eta}{\theta_1}\E{\norm{g^k - \nabla f(x^k)}^2} \notag\\ 
	&\qquad -
	\frac{\gamma}{4Ln\theta_1}\sum_{i=1}^n\norm{\nabla f_i(w^k) - \nabla f_i(x^k)}^2 
	- \frac{\gamma}{8Ln\theta_1}\sum_{i=1}^n\norm{\nabla f_i(y^k) - \nabla f_i(x^k)}^2 \notag\\
	&\qquad  
	+ 2\gamma\beta \frac{\theta_2(1+\theta_1)}{\theta_1p}(1-p)\cW^{k} 
	+2\gamma\beta \frac{\theta_2(1+\theta_1)}{\theta_1}\cY^k
	+\E{\frac{8\gamma\eta\omega}{\alpha\theta_1 n}\cH^{k+1}} \label{eq:uselm3}\\
	&\leq
	\beta\cZ^k
	+
	\left(1 - \frac{\theta_1}{2}\right)\frac{2\gamma\beta}{\theta_1}\cY^k
	+
	\left(1 - \frac{\theta_1p}{2}\right)2\gamma\beta \frac{\theta_2(1+\theta_1)}{\theta_1p}\cW^k
	\notag\\
	&\qquad -
	\frac{\gamma}{4Ln\theta_1}\sum_{i=1}^n\norm{\nabla f_i(w^k) - \nabla f_i(x^k)}^2 
	- \frac{\gamma}{8Ln\theta_1}\sum_{i=1}^n\norm{\nabla f_i(y^k) - \nabla f_i(x^k)}^2 \notag\\
	&\qquad  
	+
	\frac{\gamma\eta}{\theta_1}\E{\norm{g^k - \nabla f(x^k)}^2}
	+\E{\frac{8\gamma\eta\omega}{\alpha\theta_1 n}\cH^{k+1}} \label{eq:usetheta12}\\
	&\leq
	\beta\cZ^k
	+
	\left(1 - \frac{\theta_1}{2}\right)\frac{2\gamma\beta}{\theta_1}\cY^k
	+
	\left(1 - \frac{\theta_1p}{2}\right)2\gamma\beta \frac{\theta_2(1+\theta_1)}{\theta_1p}\cW^k
	\notag\\
	&\qquad -
	\frac{\gamma}{4Ln\theta_1}\sum_{i=1}^n\norm{\nabla f_i(w^k) - \nabla f_i(x^k)}^2 
	- \frac{\gamma}{8Ln\theta_1}\sum_{i=1}^n\norm{\nabla f_i(y^k) - \nabla f_i(x^k)}^2 \notag\\
	&\qquad  
	+
	\frac{2\gamma\eta\omega}{\theta_1 n^2}\norm{\nabla f_i(w^k) - \nabla f_i(x^k)}^2 
	+\frac{2\gamma\eta\omega}{\theta_1 n}\cH^k
	+\E{\frac{8\gamma\eta\omega}{\alpha\theta_1 n}\cH^{k+1}} \label{eq:uselm4}\\
	&\leq
	\beta\cZ^k
	+
	\left(1 - \frac{\theta_1}{2}\right)\frac{2\gamma\beta}{\theta_1}\cY^k
	+
	\left(1 - \frac{\theta_1p}{2}\right)2\gamma\beta \frac{\theta_2(1+\theta_1)}{\theta_1p}\cW^k
	\notag\\
	&\qquad -
	\frac{\gamma}{4Ln\theta_1}\sum_{i=1}^n\norm{\nabla f_i(w^k) - \nabla f_i(x^k)}^2 
	- \frac{\gamma}{8Ln\theta_1}\sum_{i=1}^n\norm{\nabla f_i(y^k) - \nabla f_i(x^k)}^2 \notag\\
	&\qquad  
	+
	\frac{2\gamma\eta\omega}{\theta_1 n^2}\norm{\nabla f_i(w^k) - \nabla f_i(x^k)}^2 
	+\frac{2\gamma\eta\omega}{\theta_1 n}\cH^k
	+\frac{8\gamma\eta\omega}{\alpha\theta_1 n}\left(1-\frac{\alpha}{2}\right)\cH^{k}
	\notag\\
	&\qquad  
	+ \left(1 + \frac{2p}{\alpha}\right)\frac{16\gamma\eta\omega p}{\alpha\theta_1 n^2}\left(\sum_{i=1}^n\norm{\nabla f_i(w^k) - \nabla f_i(x^k)}^2  
	+ \sum_{i=1}^n\norm{\nabla f_i(y^k) - \nabla f_i(x^k)}^2\right)  \label{eq:uselm5}\\
	&=
	\beta\cZ^k
	+
	\left(1 - \frac{\theta_1}{2}\right)\frac{2\gamma\beta}{\theta_1}\cY^k
	+
	\left(1 - \frac{\theta_1p}{2}\right)2\gamma\beta \frac{\theta_2(1+\theta_1)}{\theta_1p}\cW^k
	+
	\left(1 - \frac{\alpha}{4}\right)\frac{8\gamma\eta\omega}{\alpha\theta_1 n}\cH^{k}
	\notag\\
	&\qquad -
	\frac{\gamma}{n\theta_1}\left(\frac{1}{8L}-\frac{2\eta \omega}{n}\right)\sum_{i=1}^n\norm{\nabla f_i(w^k) - \nabla f_i(x^k)}^2
	\notag\\
	&\qquad  
	-\frac{\gamma}{n\theta_1}\left(\frac{1}{8L}
	-(1 + \frac{2p}{\alpha})\frac{16\eta\omega p}{\alpha n}\right)
	\left(\sum_{i=1}^n\norm{\nabla f_i(w^k) - \nabla f_i(x^k)}^2  
	+ \sum_{i=1}^n\norm{\nabla f_i(y^k) - \nabla f_i(x^k)}^2\right) \notag\\
	&\leq  \beta\cZ^k
	+
	\left(1 - \frac{\theta_1}{2}\right)\frac{2\gamma\beta}{\theta_1}\cY^k
	+
	\left(1 - \frac{\theta_1p}{2}\right)2\gamma\beta \frac{\theta_2(1+\theta_1)}{\theta_1p}\cW^k
	+
	\left(1 - \frac{\alpha}{4}\right)\frac{8\gamma\eta\omega}{\alpha\theta_1 n}\cH^{k} \label{eq:usestepsize} 
	\end{align}
	\begin{align}
	&\leq  (1-\frac{\eta \mu}{4 \theta_1})\cZ^k
	+
	\left(1 - \frac{\theta_1}{2}\right)\frac{2\gamma\beta}{\theta_1}\cY^k
	+
	\left(1 - \frac{\theta_1p}{2}\right)2\gamma\beta \frac{\theta_2(1+\theta_1)}{\theta_1p}\cW^k
	+
	\left(1 - \frac{\alpha}{4}\right)\frac{8\gamma\eta\omega}{\alpha\theta_1 n}\cH^{k} \label{eq:usebeta} \\
	&\leq  \left(1-\min\left\{\frac{\alpha}{4},\frac{p}{8},\frac{\sqrt{\eta \mu p}}{4}\right\}\right)\Psi^k, \label{eq:usetheta1}
	\end{align}
	where \eqref{eq:uselm2} uses Lemma \ref{lem:keyaccdcgd}, \eqref{eq:uselm3} uses Lemma \ref{lem:keyw}, \eqref{eq:usetheta12} uses $\theta_1\leq 1/4$ and $\theta_2= 1/2$, \eqref{eq:uselm4} uses Lemma \ref{lem:keyg}, \eqref{eq:uselm5} uses Lemma \ref{lem:keyh}, \eqref{eq:usestepsize} uses $\eta =\min\{\frac{1}{2L}, \frac{n}{64\omega(2p(\omega+1)+1)^2L}\}$, \eqref{eq:usebeta} uses $\gamma=\frac{\eta}{2(\theta_1 + \eta\mu)}$, $\beta = 1-\gamma \mu\leq 1-\frac{\eta \mu}{4 \theta_1}$ due to $\eta\mu \leq \theta_1$, and \eqref{eq:usetheta1} uses $\theta_1 = \min\left\{\frac{1}{4}, \sqrt{\frac{\eta\mu}{p}} \right\}$.
	
	Telescoping the above inequality \eqref{eq:usetheta1} from iteration $0$ to $k$, we have
	$\EE[\Psi^k] \leq \left(1-\min\big\{\frac{\alpha}{4},\frac{p}{8},\frac{\sqrt{\eta \mu p}}{4}\big\}\right)^k\Psi^0$. To obtain an $\epsilon$-solution $z^k$ such that $\EE[\ns{z^k-x^*}]\leq \epsilon$, the number of iterations is at most
	\begin{align}
	k&=\max\Big\{\frac{4}{\alpha},\frac{8}{p},\frac{4}{\sqrt{\eta \mu p}}\Big\} \log \frac{\Psi^0}{\epsilon} \notag\\
	&=\max\left\{4(1+\omega),\frac{8}{p},
	4\sqrt{\frac{L}{\mu}\max\left\{\frac{2}{p},\frac{64\omega(2p(\omega+1)+1)^2}{np}\right\}}\right\} 
	\log \frac{\Psi^0}{\epsilon}
	\end{align}
	By letting $p=\min\big\{1, \frac{\max\{1,\sqrt{n/32\omega}-1\}}{2(1+\omega)}\big\}$,  it is not hard to verify that the number of iterations performed by ADIANA (Algorithm \ref{alg:ADIANA}) to find an $\epsilon$-solution such that $\EE[\ns{z^k-x^*}]\leq \epsilon$ is at most
	\begin{equation*}
	k=
	\begin{cases}
	O\left(\left[\omega + \omega\sqrt{\frac{L}{n \mu}}\ \right]\log \frac{1}{\epsilon}  \right),& n \leq \omega,\\
	O\left(\left[\omega  + \sqrt{\frac{L}{\mu}} + \sqrt{\sqrt{\frac{\omega}{n}} \frac{\omega L}{\mu}}\ \right]\log \frac{1}{\epsilon}  \right),&  n>\omega,
	\end{cases}
	\end{equation*}
	where $n$ is the number of parallel machines.
\end{proofof}

\head{Key Lemmas} Now, the remaining thing is to prove Lemmas \ref{lem:keyaccdcgd}--\ref{lem:keyh}. We first prove the relatively simple Lemmas \ref{lem:keyw}--\ref{lem:keyh} and then prove the key Lemma \ref{lem:keyaccdcgd}.

\begin{proofof}{Lemma \ref{lem:keyw}}
	According to Line \ref{line:prob} of Algorithm \ref{alg:ADIANA}, i.e.,  
	$w^{k+1} = \begin{cases}
	y^k, &\text{with probability } p\\
	w^k, &\text{with probability } 1-p
	\end{cases}$,
	and definitions $\cY^k \eqdef P(y^k) - P(x^*)$ and
	$\cW^k \eqdef P(w^k) - P(x^*)$, this lemma is directly obtained, i.e.,
	\begin{equation}\label{eq:6}
	\E{\cW^{k+1}} = (1-p)\cW^k + p\cY^k.
	\end{equation}
\end{proofof}

\begin{proofof}{Lemma \ref{lem:keyg}} This lemma is proved as follows:
	\begin{align*}
	\E{\norm{g^k - \nabla f(x^k)}^2}
	&=
	\E{\nsB{\frac{1}{n}\sum_{i=1}^n\cC(\nabla f_i(x^k) - h_i^k) + h_i^k - \nabla f_i(x^k)}}\\
	&=
	\frac{1}{n^2}\sum_{i=1}^n \E{\nsB{\cC(\nabla f_i(x^k) - h_i^k) + h_i^k - \nabla f_i(x^k)}}\\
	&\leq
	\frac{\omega}{n^2}\sum_{i=1}^n \norm{\nabla f_i(x^k) - h_i^k}^2 \\
	&\leq
	\frac{2\omega}{n^2}\sum_{i=1}^n \norm{\nabla f_i(w^k) - \nabla f_i(x^k)}^2
	+
	\frac{2\omega}{n}\cH^k,
	\end{align*}
	where first inequality uses the property of $\omega$-compression operator (i.e., \eqref{eq:comp}) and the last inequality uses Cauchy-Schwarz inequality and definition $\cH^k \eqdef \frac{1}{n}\sum_{i=1}^n \norm{h_i^k - \nabla f_i(w^k)}^2$. 
\end{proofof}

\begin{proofof}{Lemma \ref{lem:keyh}}
	This lemma is proved as follows:
	\begin{align}
	\E{\cH^{k+1}}
	&=
	\frac{1}{n}\sum_{i=1}^n\E{\norm{h_i^{k+1} - \nabla f_i(w^{k+1})}^2} \notag\\
	&=
	\frac{p}{n}\sum_{i=1}^n\E{\norm{h_i^{k+1} - \nabla f_i(y^k)}^2}
	+
	\frac{1-p}{n}\sum_{i=1}^n\E{\norm{h_i^{k+1} - \nabla f_i(w^k)}^2} \label{eq:usenextw}\\
	&\leq
	\left(1 + \frac{2p}{\alpha}\right)\frac{p}{n}\sum_{i=1}^n\norm{\nabla f_i(w^k) - \nabla f_i(y^k)}^2
	+
	\left(\frac{1-p}{n} + 
	\left(1 + \frac{\alpha}{2p}\right)\frac{p}{n}\right)\sum_{i=1}^n\E{\norm{h_i^{k+1} - \nabla f_i(w^k)}^2} \label{eq:usecauchy1}\\
	&\leq
	\left(1 + \frac{2p}{\alpha}\right)\frac{p}{n}\sum_{i=1}^n\norm{\nabla f_i(w^k) - \nabla f_i(y^k)}^2
	+
	\left(\frac{1+\alpha/2}{n}\right)(1-\alpha)\sum_{i=1}^n\norm{h_i^{k} - \nabla f_i(w^k)}^2  \label{eq:usealpha}\\
	&\leq
	\left(1 + \frac{2p}{\alpha}\right)\frac{p}{n}\sum_{i=1}^n\norm{\nabla f_i(w^k) - \nabla f_i(y^k)}^2
	+
	\left(1-\frac{\alpha}{2}\right)\cH^k \label{eq:usedefh}\\
	&\leq \left(1-\frac{\alpha}{2} \right)\cH^k   + \left(1 + \frac{2p}{\alpha} \right)\frac{2p}{n}\left(\sum_{i=1}^n\norm{\nabla f_i(w^k) - \nabla f_i(x^k)}^2  
	+ \sum_{i=1}^n\norm{\nabla f_i(y^k) - \nabla f_i(x^k)}^2\right), \label{eq:usecauchy2}
	\end{align}
	where \eqref{eq:usenextw} uses the definition of $w^{k+1}$ (see Line \ref{line:prob} of Algorithm \ref{alg:ADIANA}), \eqref{eq:usecauchy1} uses Cauchy-Schwarz inequality, 
	\eqref{eq:usealpha} uses $h_i^{k+1}=h_i^k+\alpha \cC(\nabla f_i(w^k) - h_i^k)$, the property of $\omega$-compression operator (i.e., \eqref{eq:comp}) and $\alpha\leq 1/(1+\omega)$, 
	and \eqref{eq:usecauchy2} uses Cauchy-Schwarz inequality. 
\end{proofof}

\begin{proofof}{Lemma \ref{lem:keyaccdcgd}}
	Similar to the proof of key Lemma \ref{lem:key} in the last section, we first get the following equality:
	\begin{align}
	\E{\cZ^{k+1}}
	&=
	\EB{\norm{\beta z^k + (1-\beta)x^k - x^* + \frac{\gamma}{\eta} (y^{k+1} - x^k) }^2} \notag\\
	&=
	\norm{\beta (z^k-x^*) + (1-\beta)(x^k - x^*)}^2
	+
	\EB{\frac{2\gamma}{\eta} \dotprod{y^{k+1} - x^k}{\beta (z^k-x^k) + x^k - x^*}}
	+
	\frac{\gamma^2}{\eta^2}\E{\norm{y^{k+1} - x^k}^2} \notag\\
	&\leq
	\beta\cZ^k
	+
	(1-\beta)\norm{x^k - x^*}^2 
	-
	\beta(1-\beta)\ns{x^k-z^k}
	+
	\frac{\gamma^2}{\eta^2}\E{\norm{y^{k+1} - x^k}^2} \notag\\
	&\quad
	+\EB{
		\frac{2\gamma}{\eta}\dotprod{y^{k+1} - x^k}{x^k - x^*}
		+
		\frac{2\gamma\beta}{\eta}\dotprod{y^{k+1} - x^k}{z^k - x^k}} \notag\\
	&=
	\beta\cZ^k
	+
	(1-\beta)\norm{x^k - x^*}^2 
	-
	\beta(1-\beta)\ns{x^k-z^k}
	+
	\frac{\gamma^2}{\eta^2}\E{\norm{y^{k+1} - x^k}^2} \notag\\
	&\quad +
	\EB{\frac{2\gamma}{\eta}\dotprod{x^k - y^{k+1}}{x^* - x^k}
		+
		\frac{2\gamma\beta \theta_2}{\eta \theta_1}\dotprod{x^k - y^{k+1}}{w^k - x^k}
		+
		\frac{2\gamma\beta (1 - \theta_1 - \theta_2)}{\eta \theta_1}\dotprod{x^k - y^{k+1}}{y^k - x^k}}, \label{eq:lemz1d}
	\end{align}
	where \eqref{eq:lemz1d} uses $x^k = \theta_1 z^k + \theta_2 w^k + (1 - \theta_1 - \theta_2)y^k$ (see Line \ref{line:xk} in Algorithm \ref{alg:ADIANA}).
	To cancel the inner products in \eqref{eq:lemz1d}, we first use the  property of $f$:
	\begin{align}
	\EE[f(\xp)] &\leq \EB{f(\yt) +\inner{\nabla f(\yt)}{\xp-\yt} + \frac{L}{2}\ns{\xp-\yt}} \notag \\
	&\leq \EB{\inner{\nabla f(\yt)}{\xp-\yt} + \frac{L}{2}\ns{\xp-\yt}}  \notag\\
	&\qquad +
	f(u) -\inner{\nabla f(\yt)}{u-\yt}
	-\max\left\{\frac{\mu}{2}\ns{u-\yt},\frac{1}{2Ln}\sum_{i=1}^n\ns{\nabla f_i(u) - \nabla f_i(x^k)} \right\}, \ \  \forall u\in \R^d   \label{eq:mustrongd}
	\end{align} 
	where the inequalities hold uses the $L$-smoothness and $\mu$-strong convexity of $f$ and $f(x)\eqdef\frac{1}{n}\sum_{i=1}^{n}f_i(x)$.
	
	Then, according the definition of $y^{k+1}$ (see Line \ref{line:yp} of Algorithm \ref{alg:ADIANA}), we have
	\begin{equation}\label{eq:defyp}
	y^{k+1} = x^k - \eta g^k - \eta \Delta,
	\end{equation}
	where $\Delta \in \partial \psi(y^{k+1})$. Besides, according to convexity of $\psi$, we have
	\begin{align}
	\EE[\psi(y^{k+1})] &\leq \EB{\psi(u)-\inner{\Delta}{u-y^{k+1}}}
	=\EB{\psi(u)-\inner{\Delta}{u-x^{k}}  +\inner{\Delta}{y^{k+1}-x^{k}}},  \ \   \forall u\in \R^d. \label{eq:cvxpsi} 
	\end{align} 
	Adding \eqref{eq:mustrongd} and \eqref{eq:cvxpsi}, we have
	\begin{align}
	\forall u\in \R^d,\quad	\EE[P(y^{k+1})] 
	&\leq \EB{P(u)
		-\inner{\Delta+g^k}{u-x^{k}}  +\inner{\Delta+\nabla f(x^k)}{y^{k+1}-x^{k}} +\frac{L}{2}\ns{\xp-\yt}} \notag\\
	&\qquad\quad 
	-\max\Big\{\frac{\mu}{2}\ns{u-\yt},\frac{1}{2Ln}\sum_{i=1}^n\ns{\nabla f_i(u) - \nabla f_i(x^k)} \Big\} \notag\\
	&= \EB{P(u)
		-\frac{1}{\eta}\inner{x^k-y^{k+1}}{u-x^{k}}  +\inner{\nabla f(x^k) -g^k}{y^{k+1}-x^{k}} -\frac{1}{\eta}\ns{\xp-\yt}} \notag\\
	&\qquad\quad +\EB{\frac{L}{2}\ns{\xp-\yt}}
	-\max\Big\{\frac{\mu}{2}\ns{u-\yt},\frac{1}{2Ln}\sum_{i=1}^n\ns{\nabla f_i(u) - \nabla f_i(x^k)} \Big\} \label{eq:useyp}\\
	&\leq \EB{P(u)
		-\frac{1}{\eta}\inner{x^k-y^{k+1}}{u-x^{k}}  +\frac{\eta}{2}\ns{\nabla f(x^k) -g^k} -\frac{1}{2\eta}\ns{y^{k+1}-x^{k}}} \notag\\
	&\qquad\quad +\EB{\frac{L}{2}\ns{\xp-\yt}}
	-\max\Big\{\frac{\mu}{2}\ns{u-\yt},\frac{1}{2Ln}\sum_{i=1}^n\ns{\nabla f_i(u) - \nabla f_i(x^k)} \Big\} \label{eq:useyoung}\\
	&\leq \EB{P(u)
		-\frac{1}{\eta}\inner{x^k-y^{k+1}}{u-x^{k}}  +\frac{\eta}{2}\ns{\nabla f(x^k) -g^k} -\frac{1}{4\eta}\ns{y^{k+1}-x^{k}}} \notag\\
	&\qquad\quad 
	-\max\Big\{\frac{\mu}{2}\ns{u-\yt},\frac{1}{2Ln}\sum_{i=1}^n\ns{\nabla f_i(u) - \nabla f_i(x^k)} \Big\}, \label{eq:inner}
	\end{align} 
	where \eqref{eq:useyp} uses \eqref{eq:defyp}, \eqref{eq:useyoung} uses Young's inequality, and \eqref{eq:inner} uses the condition $\eta\leq 1/2L$.
	
	Now, we are ready to prove this lemma by canceling the inner products in \eqref{eq:lemz1d} using \eqref{eq:inner}. 
	By plugging $2\gamma$ times \eqref{eq:inner} (where $u=x^*$),  
	$\frac{2\gamma\beta \theta_2}{\theta_1}$ times \eqref{eq:inner} (where $u=w^k$), and $\frac{2\gamma\beta (1-\theta_1-\theta_2)}{\theta_1}$ times \eqref{eq:inner} (where $u=y^k$) into \eqref{eq:lemz1d}, we have
	\begin{align}
	\E{\cZ^k}
	&\leq
	\beta\cZ^k
	+
	(1-\beta)\norm{x^k - x^*}^2 
	-
	\beta(1-\beta)\ns{x^k-z^k}
	+
	\frac{\gamma^2}{\eta^2}\E{\norm{y^{k+1} - x^k}^2} \notag\\
	&\qquad +
	2\gamma
	\EB{P(x^*) - P(y^{k+1}) + \frac{\eta}{2}\norm{\nabla f(x^k) -g^k}^2 - \frac{1}{4\eta}\norm{y^{k+1} - x^k}^2  - \frac{\mu}{2}\norm{x^k - x^*}^2 } \notag\\
	&\qquad+
	\frac{2\gamma\beta\theta_2}{\theta_1}
	\EE\Big[P(w^k) - P(y^{k+1}) + \frac{\eta}{2}\norm{\nabla f(x^k) -g^k}^2 - \frac{1}{4\eta}\norm{y^{k+1} - x^k}^2  \notag\\
	&\qquad\qquad\qquad\qquad\qquad\qquad\qquad\qquad\qquad\qquad\qquad\qquad\qquad\qquad
	- \frac{1}{2Ln}\sum_{i=1}^n\norm{\nabla f_i(w^k) - \nabla f_i(x^k)}^2 \Big]\notag\\
	&\qquad +
	\frac{2\gamma\beta(1 - \theta_1 - \theta_2)}{\theta_1}
	\EE\Big[P(y^k) - P(y^{k+1}) + \frac{\eta}{2}\norm{\nabla f(x^k) -g^k}^2  - \frac{1}{4\eta}\norm{y^{k+1} - x^k}^2]  \notag\\
	&\qquad\qquad\qquad\qquad\qquad\qquad\qquad\qquad\qquad\qquad\qquad\qquad\qquad\qquad
	- \frac{1}{2Ln}\sum_{i=1}^n\norm{\nabla f_i(y^k) - \nabla f_i(x^k)}^2\Big] \notag
	\end{align}
	\begin{align}
	&=
	\beta\cZ^k
	+
	(1-\beta - \gamma\mu)\norm{x^k - x^*}^2
	-
	\beta(1-\beta)\ns{x^k-z^k}
	+
	\left(\frac{\gamma^2}{\eta^2} -  \frac{\gamma\beta}{2\eta\theta_1} \right)
	\E{\norm{y^{k+1} - x^k}^2} \notag\\
	&\qquad+
	\frac{2\gamma\beta\theta_2}{\theta_1}\cW^k
	+
	\frac{2\gamma\beta(1 - \theta_1 - \theta_2)}{\theta_1}\cY^k
	-
	\frac{2\gamma\beta}{\theta_1}\E{\cY^{k+1}} \notag\\
	&\qquad-
	2\gamma(1-\beta)\E{\cY^{k+1}}
	+
	\left(\gamma\eta  + \gamma\beta\eta\frac{1-\theta_1}{\theta_1} \right)\E{\norm{\nabla f(x^k) -g^k}^2}\notag\\
	&\quad-
	\frac{\gamma\beta\theta_2}{Ln\theta_1}\sum_{i=1}^n\norm{\nabla f_i(w^k) - \nabla f_i(x^k)}^2
	-
	\frac{\gamma\beta(1 - \theta_1 - \theta_2)}{Ln\theta_1}\sum_{i=1}^n\norm{\nabla f_i(y^k) - \nabla f_i(x^k)}^2 \notag\\
	&\leq
	\beta\cZ^k
	+
	\frac{2\gamma\beta(1 - \theta_1 - \theta_2)}{\theta_1}\cY^k
	+
	\frac{2\gamma\beta\theta_2}{\theta_1}\cW^k
	-
	\frac{2\gamma\beta}{\theta_1}\E{\cY^{k+1}}
	+
	\frac{\gamma\eta}{\theta_1}\E{\norm{\nabla f(x^k) -g^k}^2}  \notag\\
	&\qquad -
	\frac{\gamma\theta_2}{2Ln\theta_1}\sum_{i=1}^n\norm{\nabla f_i(w^k) - \nabla f_i(x^k)}^2
	-
	\frac{\gamma(1-\theta_1-\theta_2)}{2Ln\theta_1}\sum_{i=1}^n\norm{\nabla f_i(y^k) - \nabla f_i(x^k)}^2 \label{eq:usebetad}\\
	&\leq
	\beta\cZ^k
	+
	\frac{2\gamma\beta(1 - \theta_1 - \theta_2)}{\theta_1}\cY^k
	+
	\frac{2\gamma\beta\theta_2}{\theta_1}\cW^k
	-
	\frac{2\gamma\beta}{\theta_1}\E{\cY^{k+1}}
	+
	\frac{\gamma\eta}{\theta_1}\E{\norm{\nabla f(x^k) -g^k}^2}  \notag\\
	&\qquad -
	\frac{\gamma}{4Ln\theta_1}\sum_{i=1}^n\norm{\nabla f_i(w^k) - \nabla f_i(x^k)}^2 
	- \frac{\gamma}{8Ln\theta_1}\sum_{i=1}^n\norm{\nabla f_i(y^k) - \nabla f_i(x^k)}^2,
	\label{eq:usetheta1d}
	\end{align}
	where \eqref{eq:usebetad} uses $\gamma=\frac{\eta}{2(\theta_1 + \eta\mu)}$ and $\beta = 1-\gamma \mu$, and \eqref{eq:usetheta1d} uses $\theta_1\leq 1/4$ and $\theta_2= 1/2$.
\end{proofof}

\newpage
\section{Extra Experiments}

In this section, we conduct more experiments for these methods with different compression operators for the same regularized logistic regression problem on other two datasets \texttt{a9a}  and \texttt{w6a}. 
Similar to Section \ref{sec:exp}, for all methods, we use the theoretical stepsize and parameters.

\subsection{Comparison with DIANA and DCGD}

\begin{figure*}[h]
	\vspace{0cm}
	\centering
	\begin{tabular}{cccc}
		\includegraphics[width=0.32\linewidth]{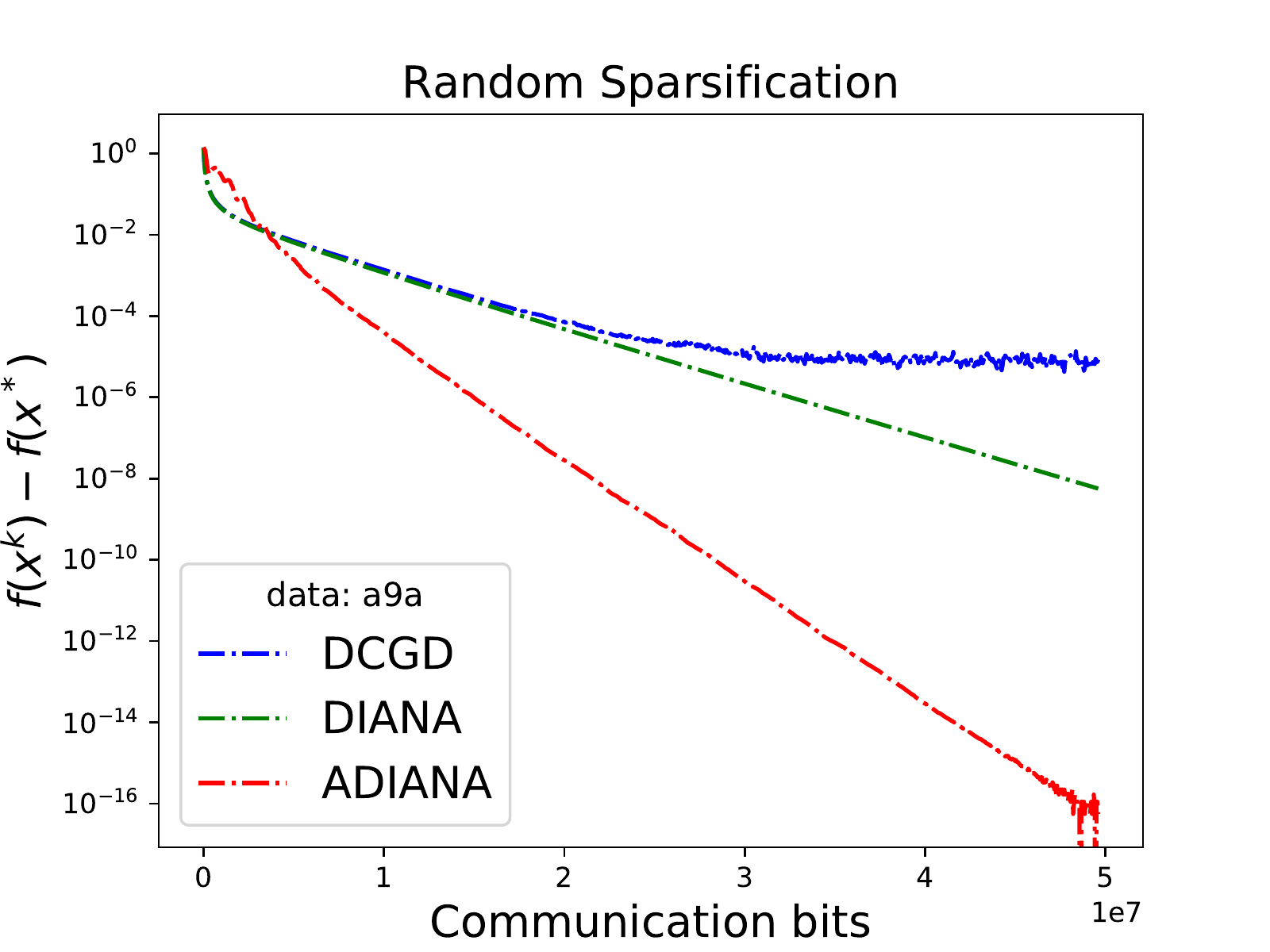}&
		\includegraphics[width=0.32\linewidth]{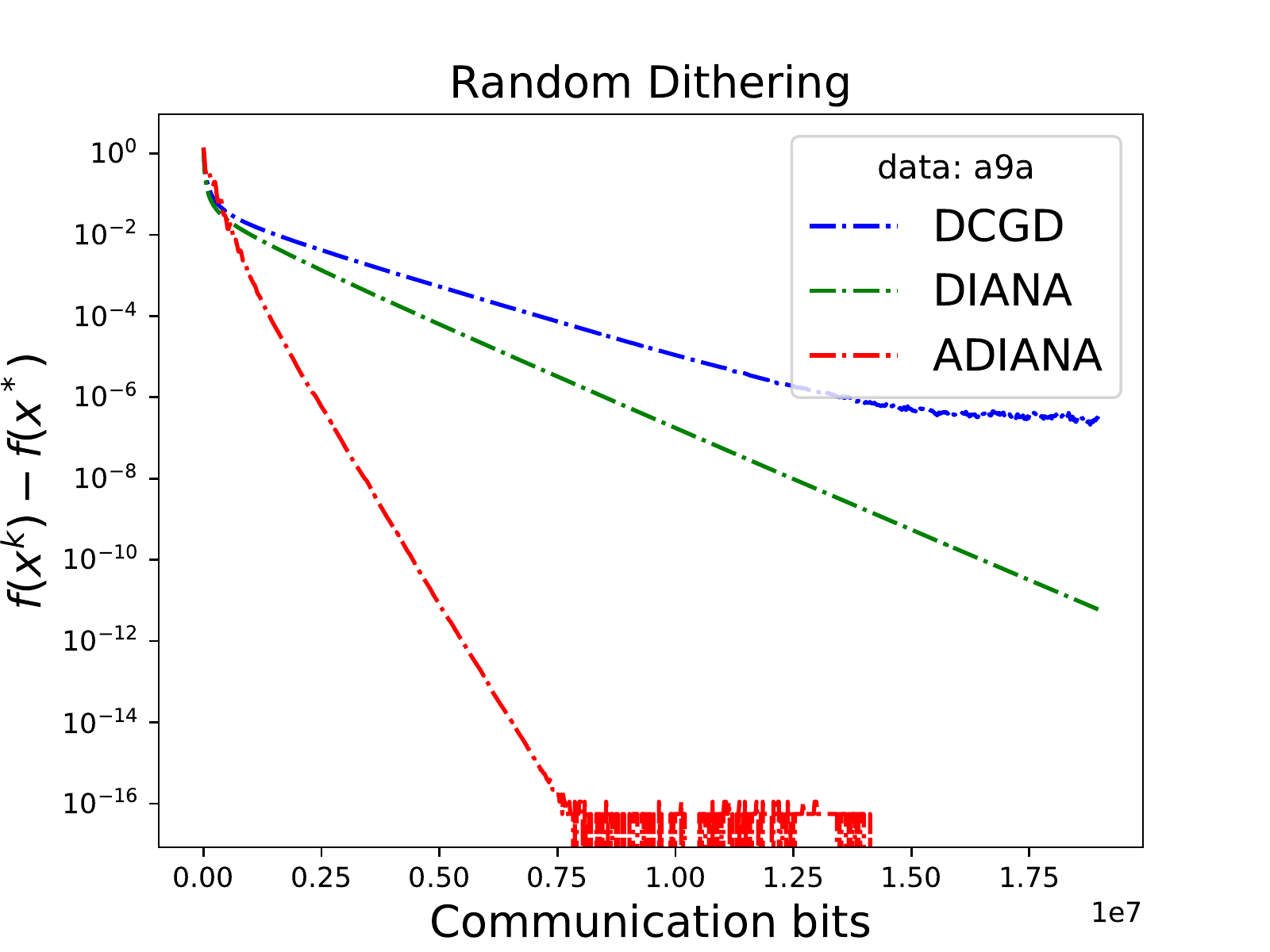}&
		\includegraphics[width=0.32\linewidth]{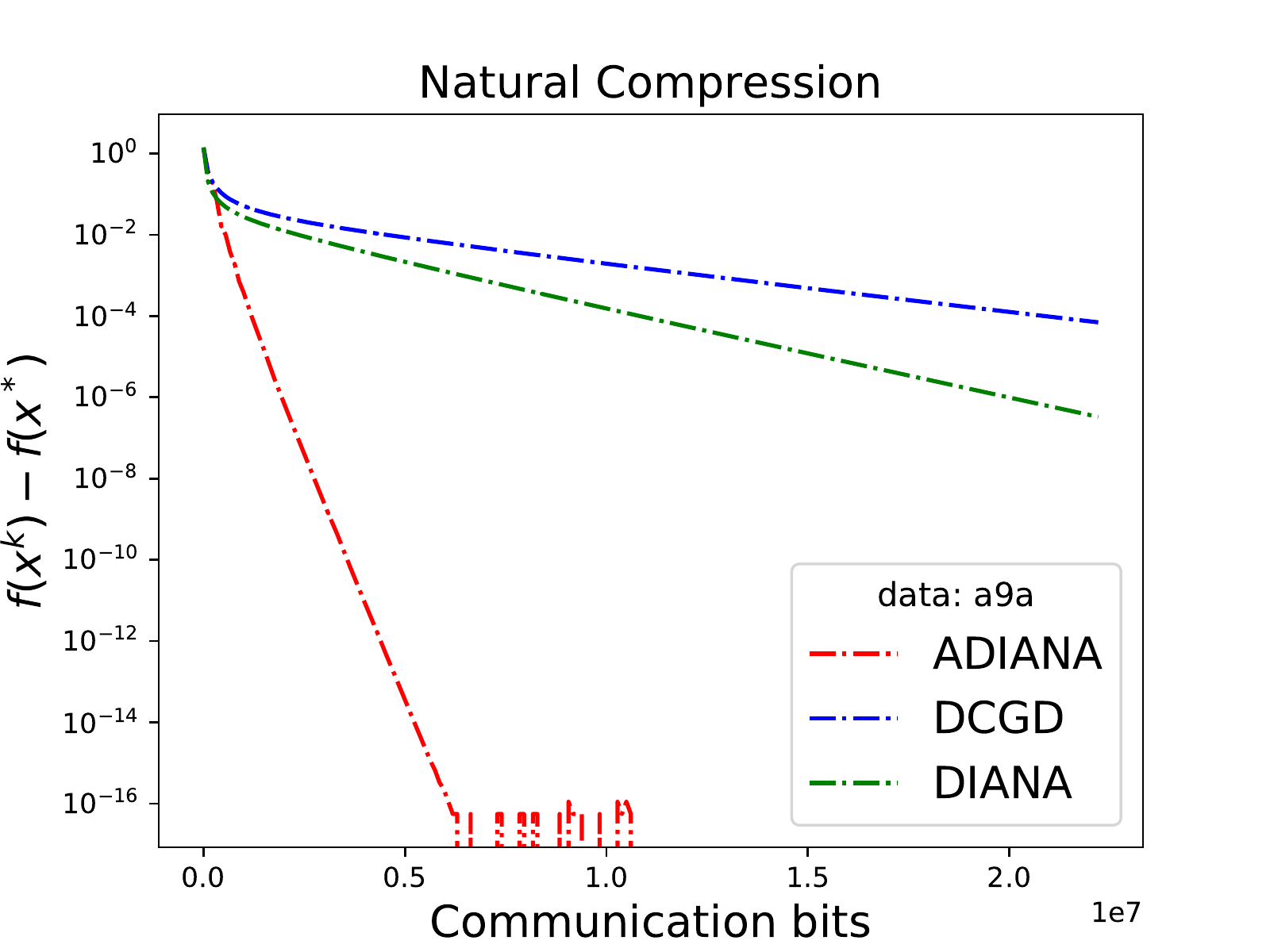}&	
	\end{tabular}

	\caption{The communication complexity of different methods for three different compressors (random sparsification, random dithering and natural compression) on the \texttt{a9a} dataset.}
	\label{fig:a9a111}
	
\end{figure*}

\begin{figure*}[h]
	\vspace{0.2cm}
	\centering
	\begin{tabular}{cccc}
		\includegraphics[width=0.32\linewidth]{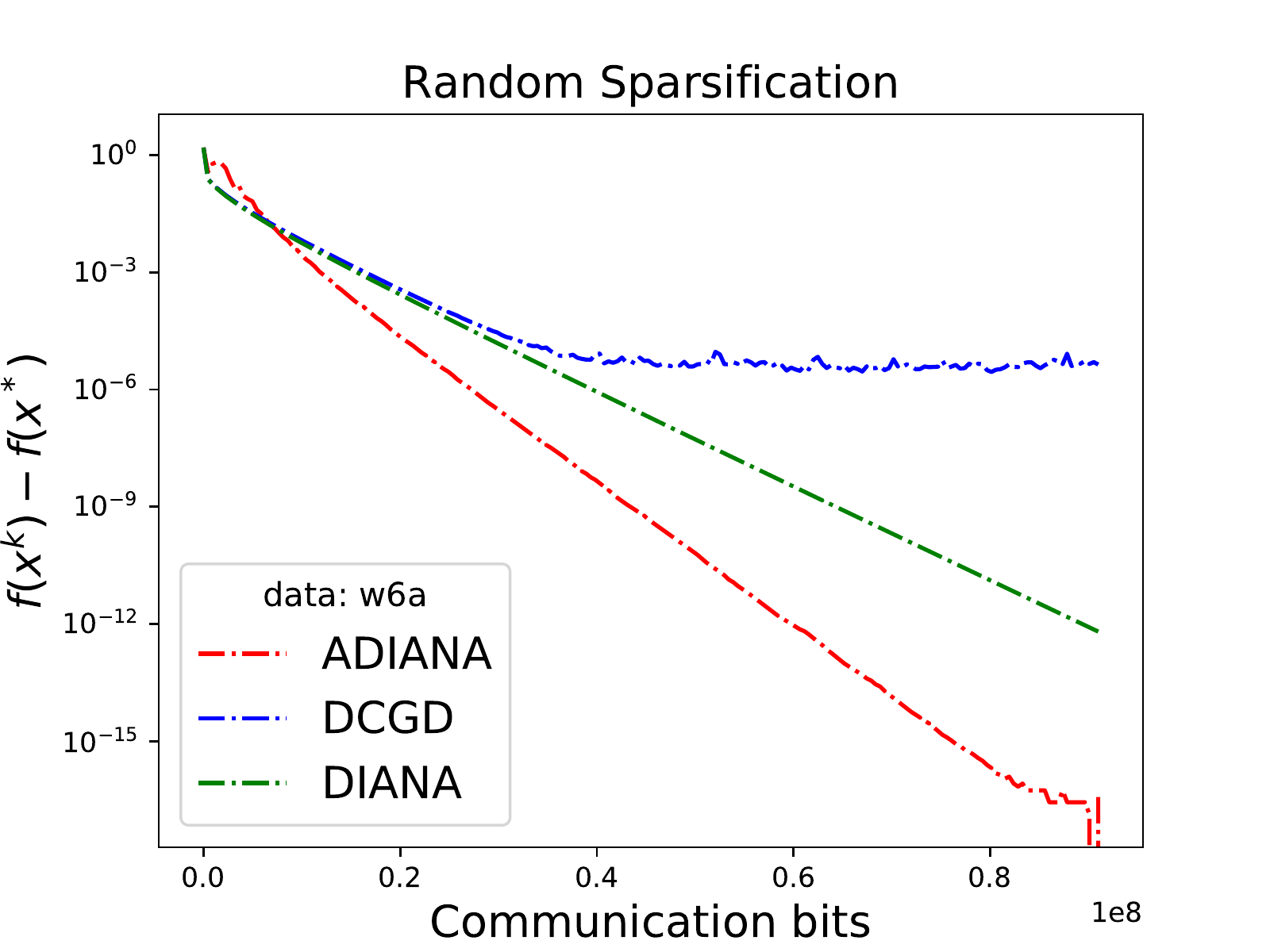}&
		\includegraphics[width=0.32\linewidth]{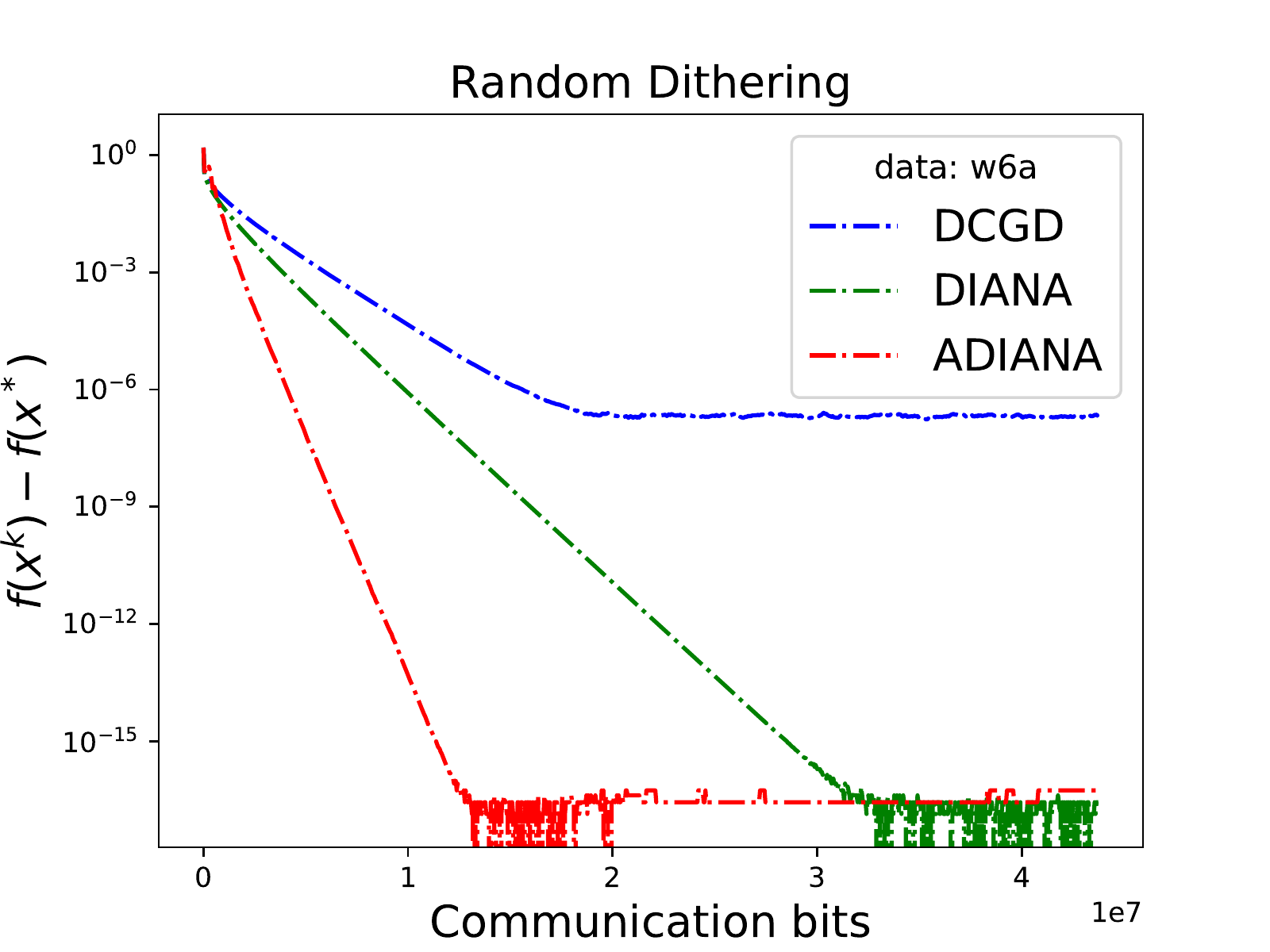}&
		\includegraphics[width=0.32\linewidth]{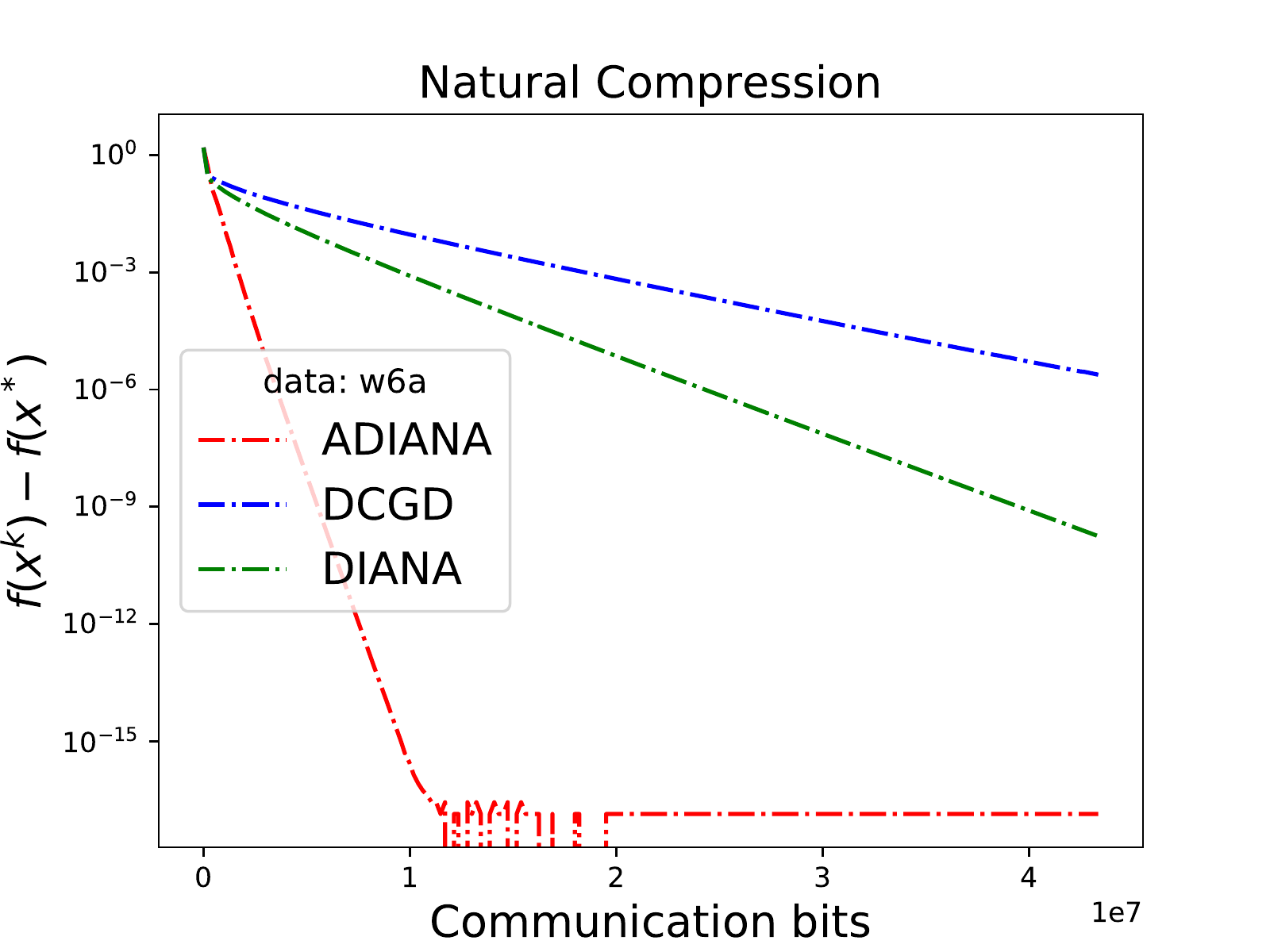}&
	\end{tabular}

	\caption{The communication complexity of different methods for three different compressors (random sparsification, random dithering and natural compression) on the \texttt{w6a} dataset.}
	\label{fig:w6a111}
	
\end{figure*}

Similar to Section \ref{sec:compare-diana-dcgd}, we compare our ADIANA with DIANA and DCGD with three compression operators: random sparsification, random dithering, and natural compression on different datasets in Figures~\ref{fig:a9a111} and \ref{fig:w6a111}. 
Similar to Figures~\ref{fig:a5a111} and \ref{fig:mushrooms111}, here Figures~\ref{fig:a9a111} and \ref{fig:w6a111} also show that our ADIANA converges fastest for all three compressors, and natural compression uses the fewest communication bits than random dithering and random sparsification.
Furthermore, because the compression error of vanilla DCGD is nonzero in general, DCGD can only converge to the neighborhood of the optimal solution while DIANA and ADIANA can converge to the optimal solution.

\subsection{Communication efficiency}

\begin{figure*}[h]
	\vspace{0cm}
	\centering
	\begin{tabular}{cccc}
		\includegraphics[width=0.32\linewidth]{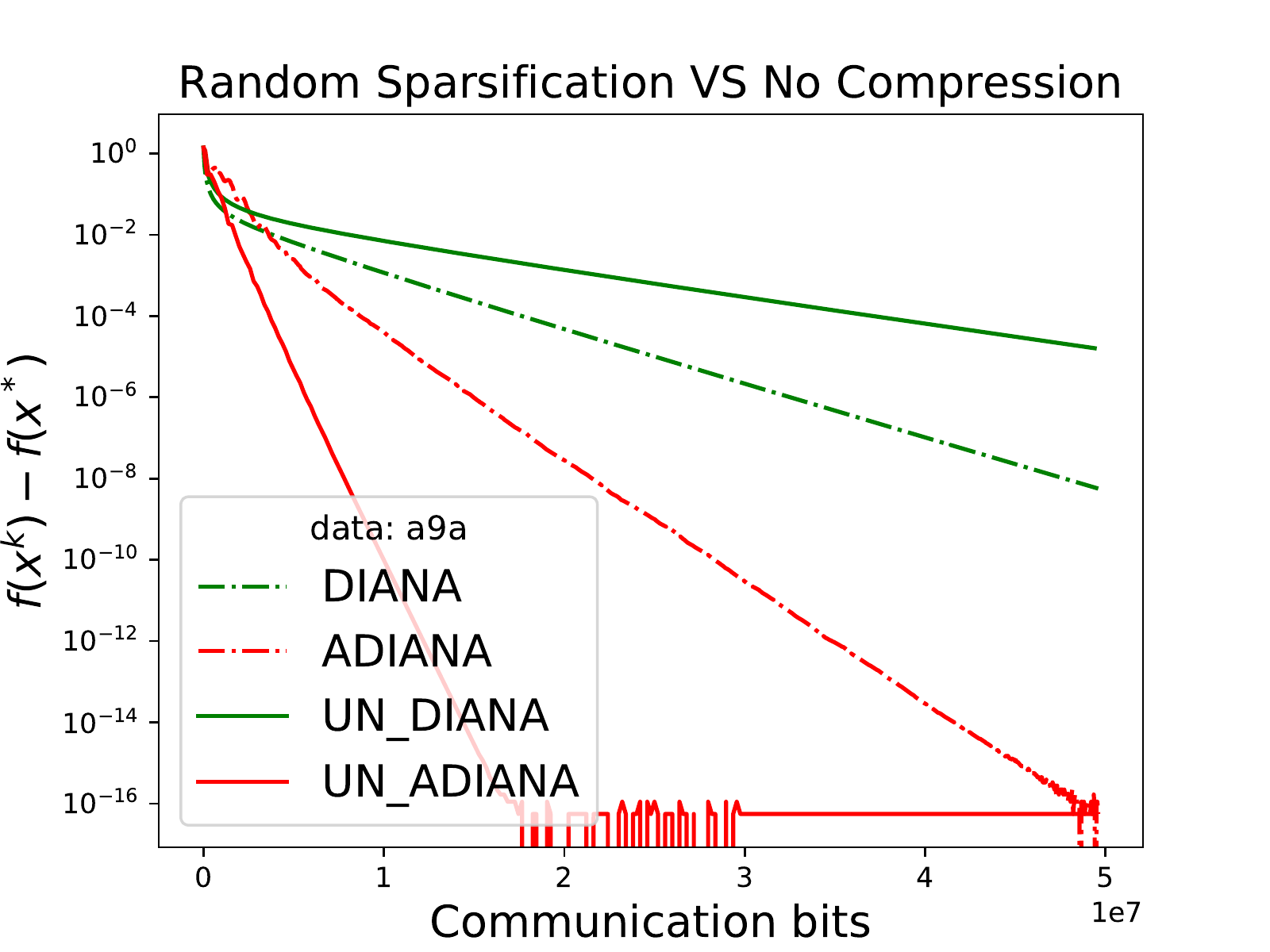}&
		\includegraphics[width=0.32\linewidth]{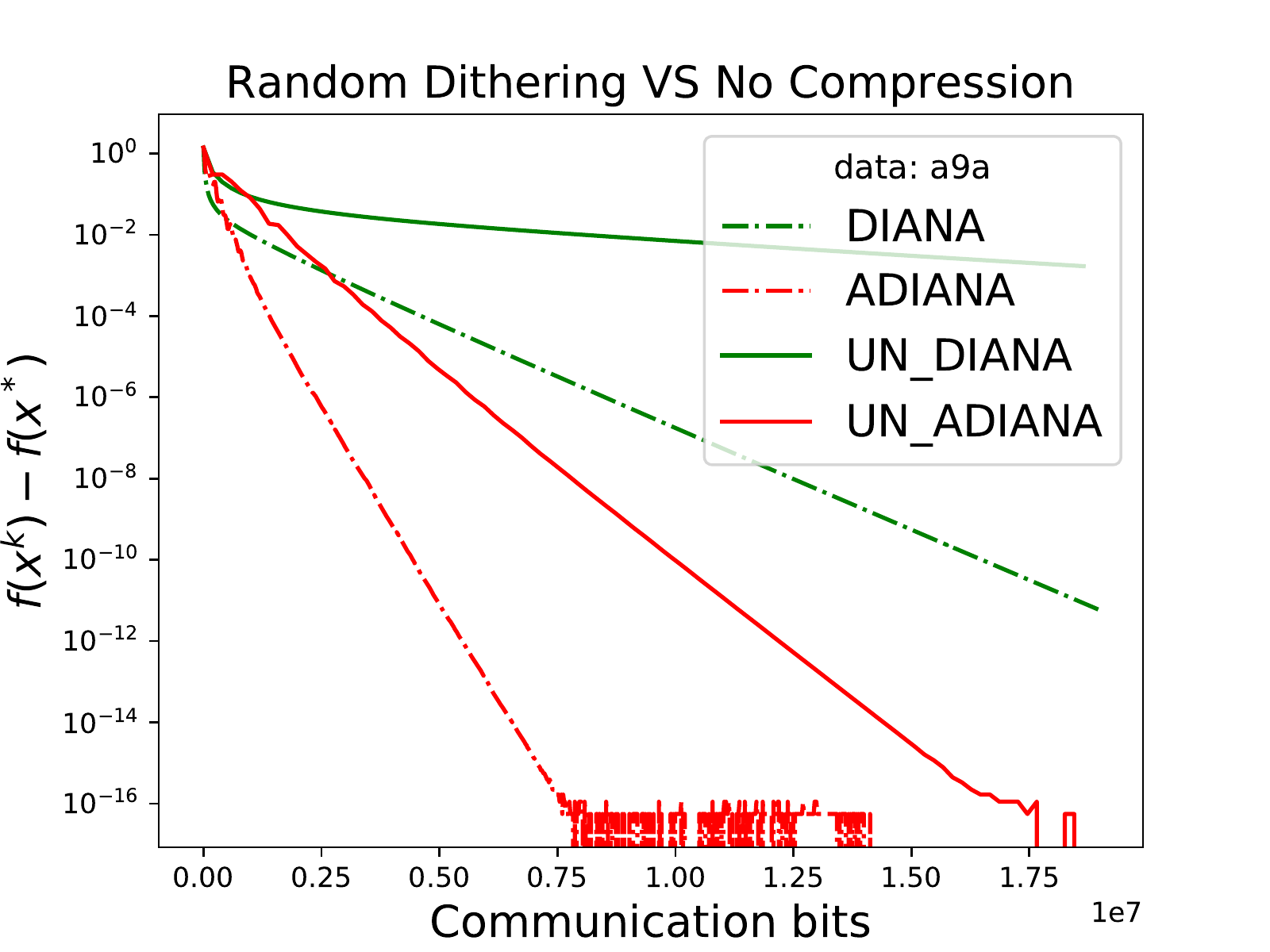}&
		\includegraphics[width=0.32\linewidth]{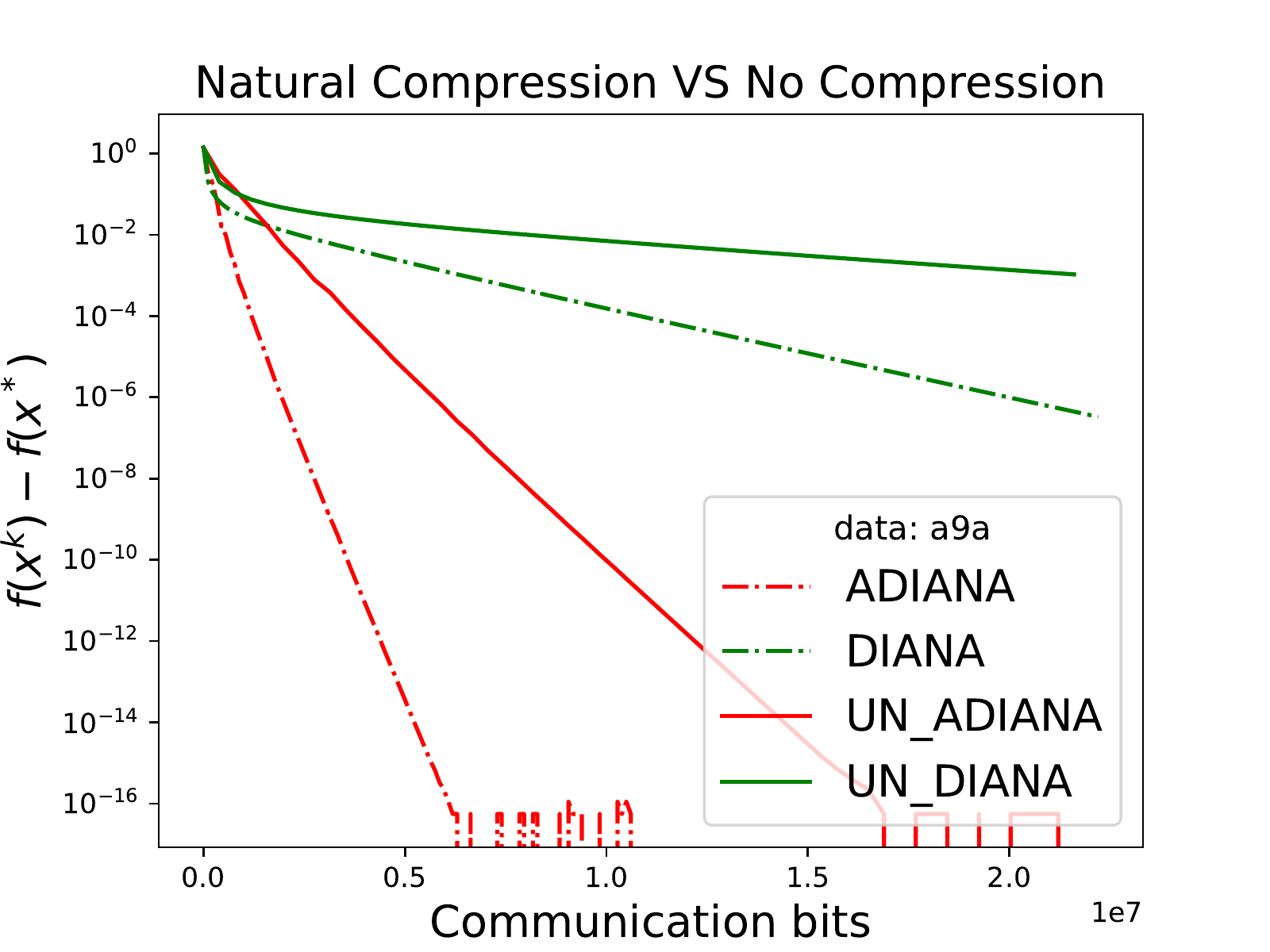}&
		
	\end{tabular}
	\vskip -0.2cm
	\caption{The communication complexity of DIANA and ADIANA with and without compression on the \texttt{a9a} dataset.}
	\label{fig:a9a222}
	
\end{figure*}

\begin{figure*}[h]
	\vspace{0cm}
	\centering
	\begin{tabular}{cccc}
		\includegraphics[width=0.32\linewidth]{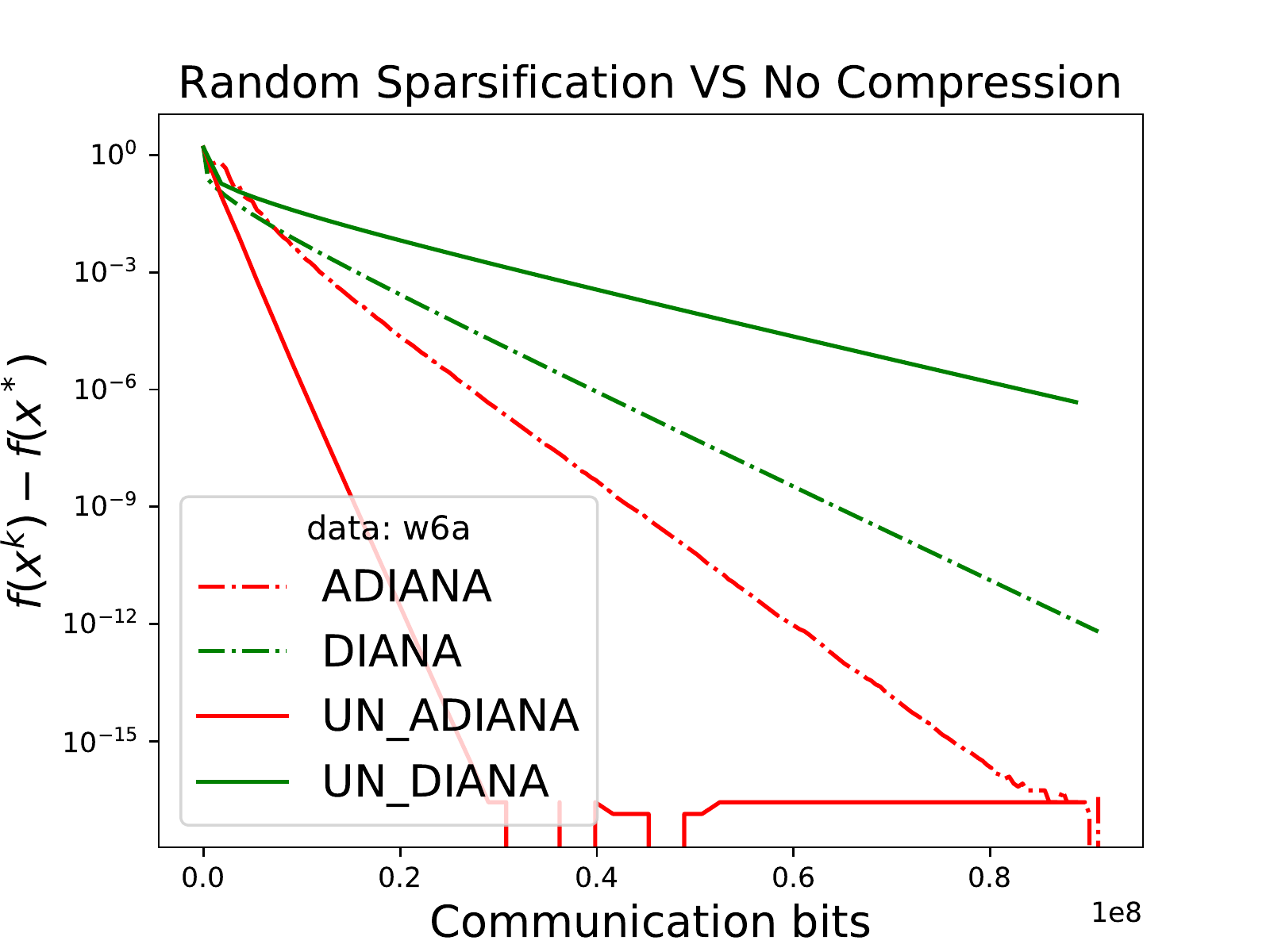}&
		\includegraphics[width=0.32\linewidth]{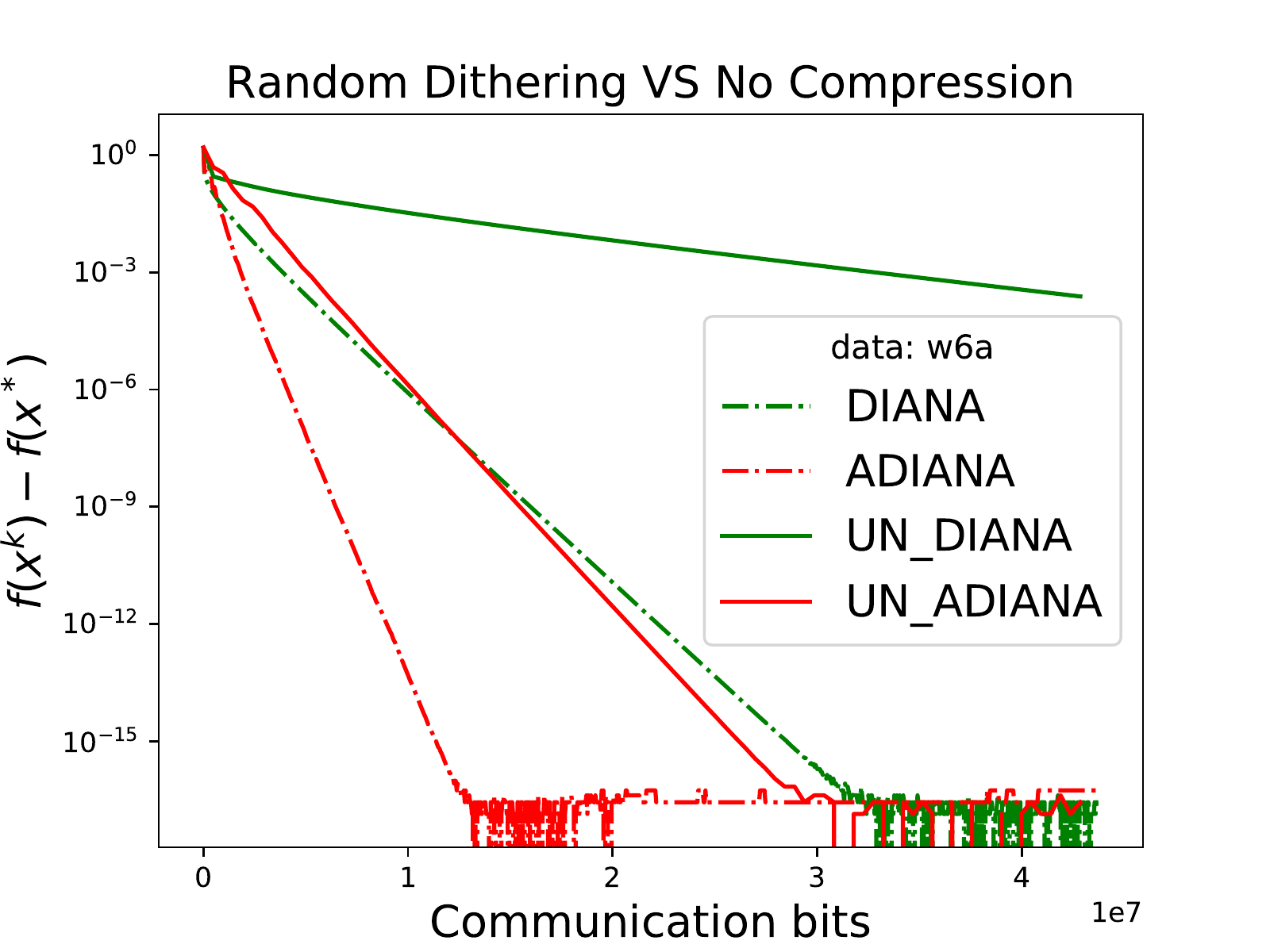}&
		\includegraphics[width=0.32\linewidth]{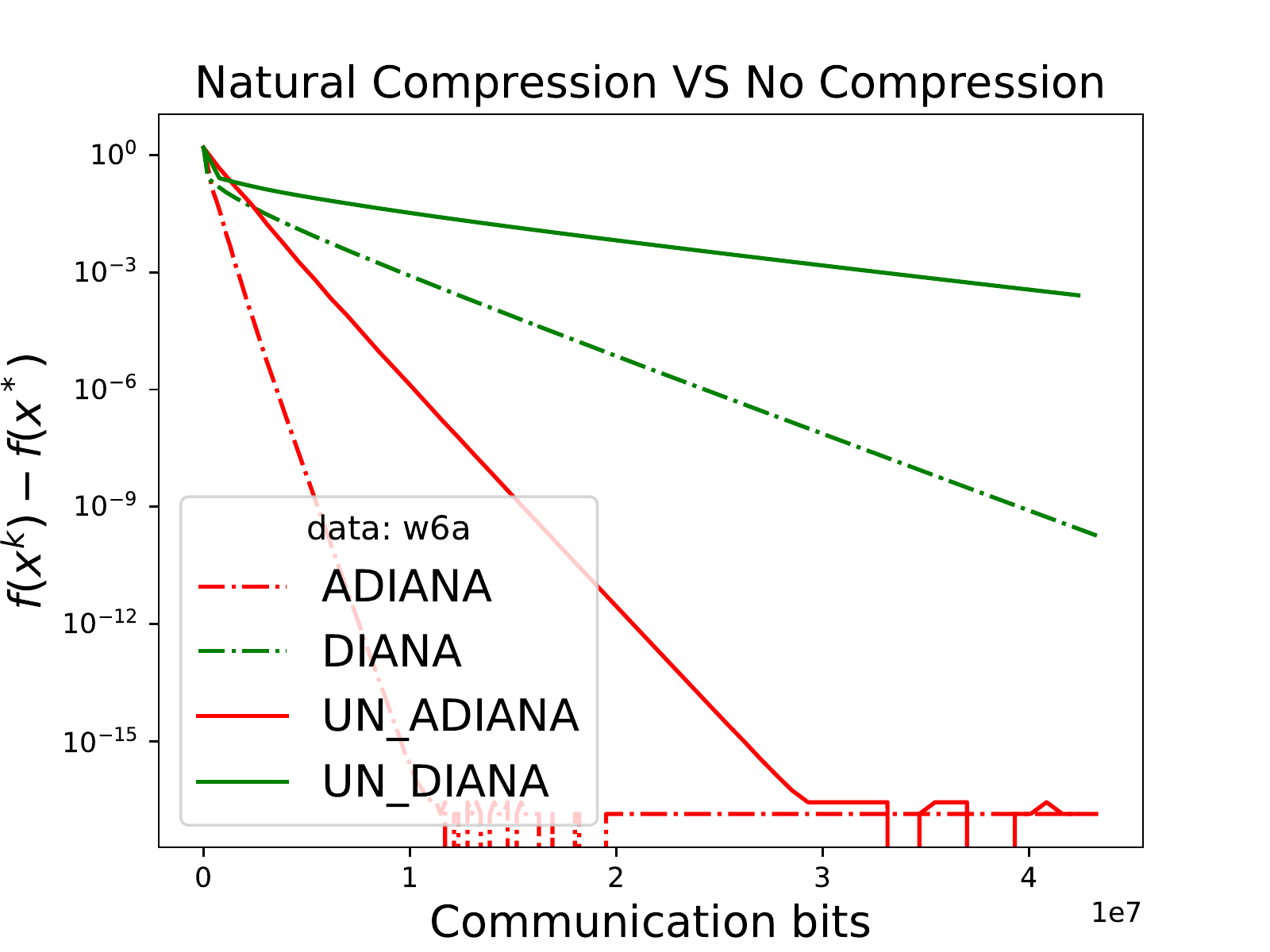}&
		
	\end{tabular}
	\vskip -0.2cm
	\caption{The communication complexity of DIANA and ADIANA with and without compression on the \texttt{w6a} dataset.}
	\label{fig:w6a222}
	
\end{figure*}

Similar to Section \ref{sec:compare-comp-noncomp}, we compare our ADIANA and DIANA, with and without compression on different datasets to show the communication efficiency of our accelerated method ADIANA in Figures~\ref{fig:a9a222} and  \ref{fig:w6a222}.
Similarly in Figures~\ref{fig:a9a111} and \ref{fig:w6a111}, DIANA is better than its uncompressed version if the compression operator is random sparsification. However, ADIANA behaves worse than its uncompressed version. 
For random dithering (middle figures) and natural compression (right figures), ADIANA is about twice faster than its uncompressed version, and is much faster than DIANA with/without compression.
These numerical results indicate that ADIANA (which enjoys both acceleration and compression) could be a more practical communication efficiency method, i.e., acceleration (better than non-accelerated DIANA) and compression (better than the uncompressed version), especially for random dithering and natural compression.

\newpage
\subsection{Different number of nodes}
\label{sec:diffnode}
\begin{figure*}[!h]
	\centering
	\begin{tabular}{cc}
		\includegraphics[width=0.35\linewidth]{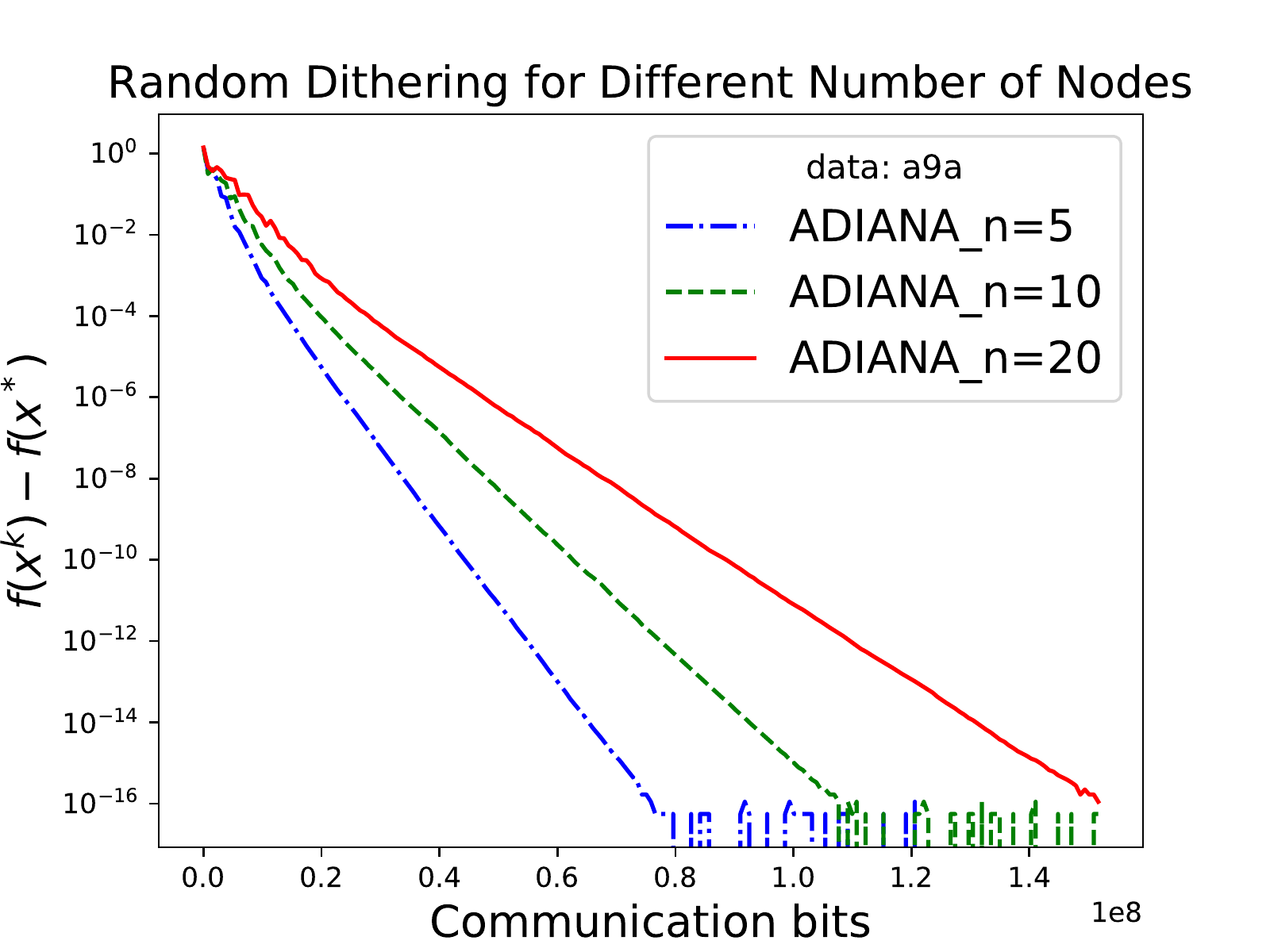}&
			\includegraphics[width=0.35\linewidth]{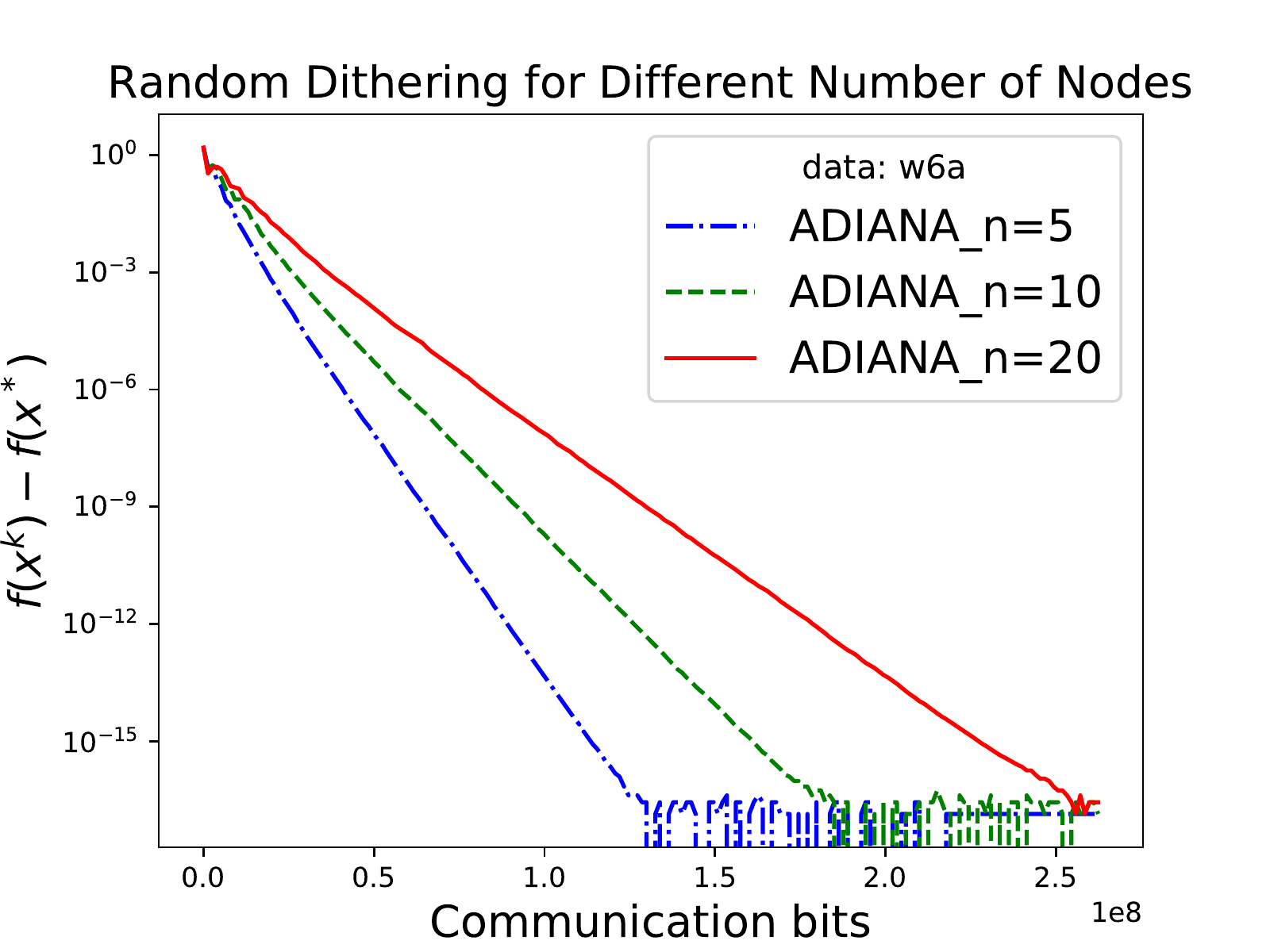}
	\end{tabular}
%
	\vspace{-2mm}
	\centering
	\begin{tabular}{cc}
		\includegraphics[width=0.35\linewidth]{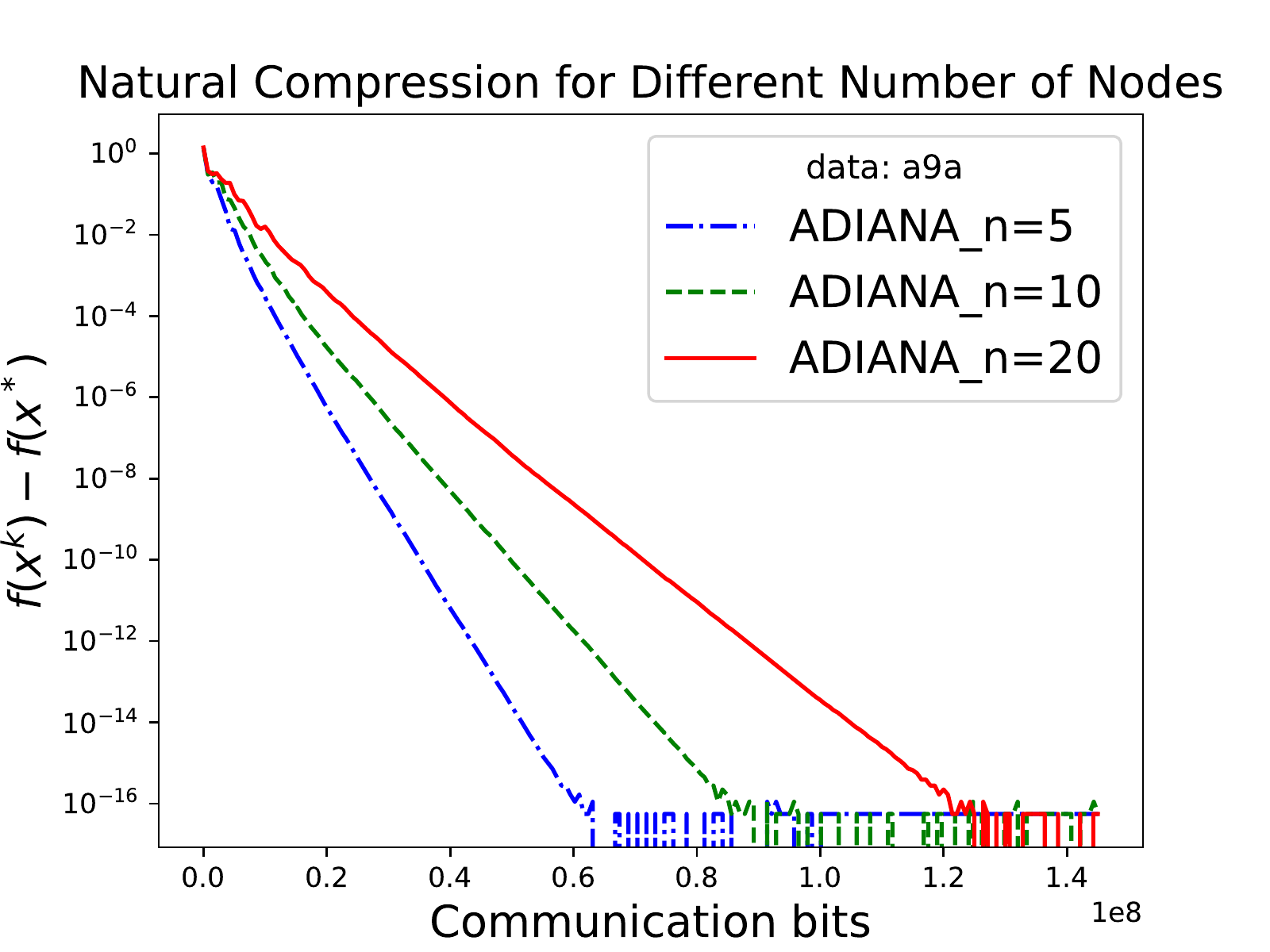}&
		\includegraphics[width=0.35\linewidth]{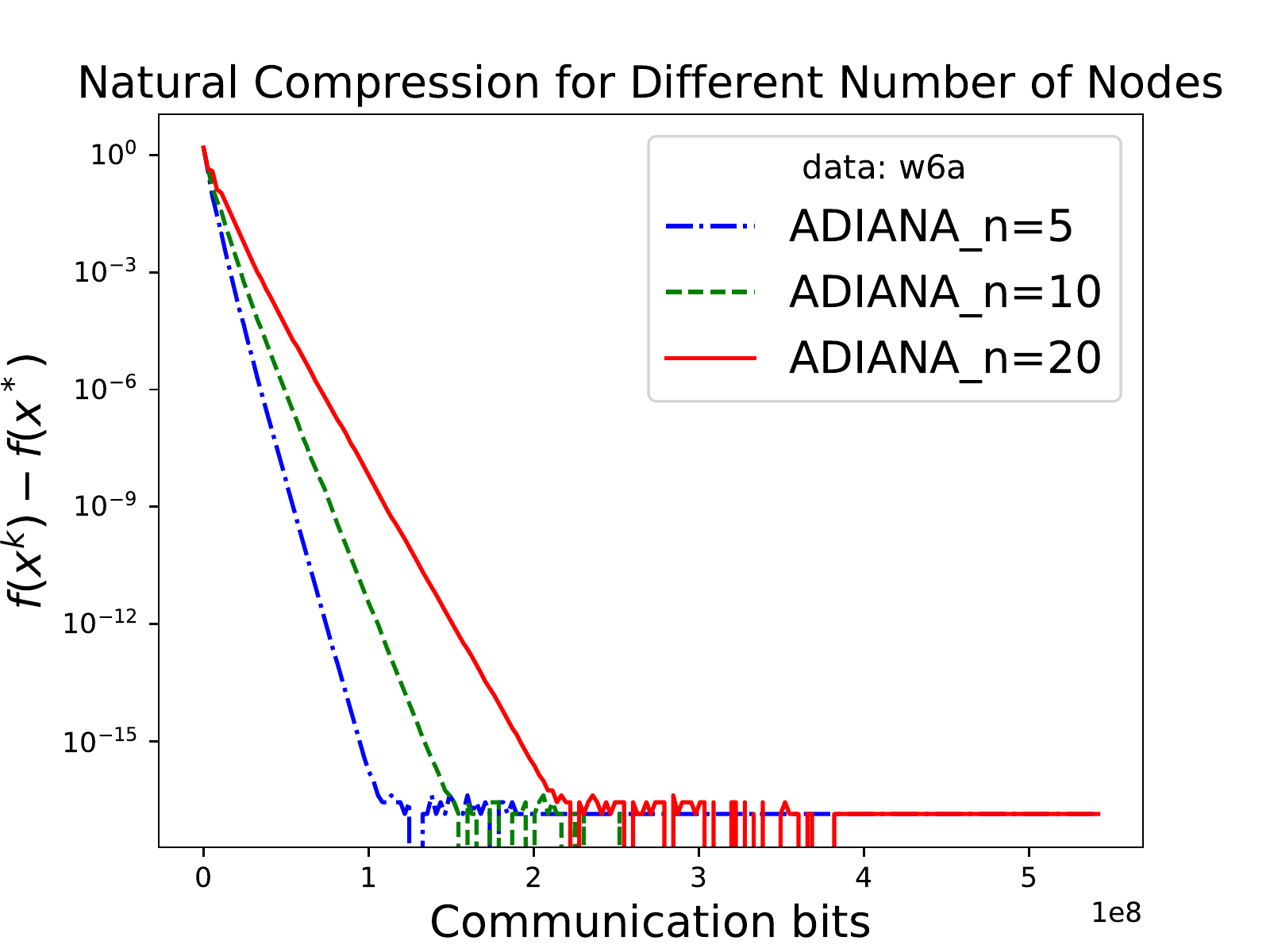}
	\end{tabular}
	
	\caption{The communication complexity of ADIANA with different number of nodes on the \texttt{a9a} (left) and \texttt{w6a} (right) datasets.}
	\label{fig:w6a333}
		\vspace{-3mm}
\end{figure*}

In this subsection, we demonstrate the performance of our ADIANA for different number of nodes/machines with random dithering and natural compression on \texttt{a9a} and \texttt{w6a} datasets in Figure~\ref{fig:w6a333}. 
In the previous numerical results, the number of communication bits is not multiplied by the number of nodes, since the number of nodes is the same for all methods. However, here the number of communication bits is multiplied by the number of nodes since they are different now.

\end{document}